\documentclass[a4paper,11pt]{amsart}

\usepackage{graphicx}
\usepackage{amssymb} 

\newcommand{\ie}{i\textup.e\textup.}
\newcommand{\eg}{e\textup.g\textup.}

\newcommand{\rank}{\operatorname{rk}}
\newcommand{\cl}{\operatorname{cl}}
\newcommand{\cw}{\operatorname{cw}}
\newcommand{\vtxlnk}{\operatorname{Link}}
\newcommand{\pathlbl}{\ell}
\newcommand{\cntr}{\partial}
\newcommand{\ccntr}{\partial^\star}

\newcommand{\cntrcclof}[1]{\langle\cntr#1\rangle}
\newcommand{\ccntrcclof}[1]{\langle\ccntr#1\rangle}

\providecommand{\abs}[1]{\lvert#1\rvert}
\providecommand{\norm}[1]{\lVert#1\rVert}

\newcommand{\GP}[2]{$\langle\,#1\,\Vert\,#2\,\rangle$}

\newtheorem{theorem}{Theorem}
\newtheorem{proposition}[theorem]{Proposition}

\newtheorem{lemma}[theorem]{Lemma}
\newtheorem{estimating_lemma}[theorem]{Estimating Lemma}
\newtheorem{inductive_lemma}[theorem]{Inductive Lemma}
\newtheorem{named_lemma:p.hall}[theorem]{Lemma}
\newtheorem{corollary.named_lemma:p.hall}[theorem]{Corollary}
\theoremstyle{remark}
\newtheorem{remark}[theorem]{Remark}
\theoremstyle{definition}
\newtheorem{definition}[theorem]{Definition}

\begin{document}



\title[Simple groups of infinite commutator width]
    {Finitely generated infinite simple groups 
    of infinite commutator width}

\author{Alexey Muranov}

\address{Institut Camille Jordan\\
Universit\'e Claude Bernard Lyon 1\\
43 blvd du 11 novembre 1918\\
69622 Villeurbanne Cedex\\
France}

\curraddr{Institut de Math\'ematiques de Toulouse\\
Universit\'e Paul Sabatier Toulouse 3\\
118 route de Narbonne\\
F--31062 Toulouse Cedex 9\\
France}

\email{muranov@math.univ-toulouse.fr}

\date{\today}


\thanks{This work was supported in part by the NSF grant DMS 0245600
of Alexander Ol'shanskii and Mark Sapir.
Currently, the author is supported by 
Chateaubriand Fellowship of French government.}





\begin{abstract}
It is shown that
there exists a finitely generated infinite simple group 
of infinite commutator width, 
and that the commutator width of a finitely generated infinite 
\emph{boundedly\/} simple group can be arbitrarily large.
Besides, such groups can be constructed with decidable word
and conjugacy problems.
\end{abstract}


\maketitle

\tableofcontents



\section{Introduction}
\label{section:introduction}

In 1951, Oystein Ore (see \cite{Ore:1951:src}) conjectured that
all elements in every non-abelian finite simple group 
are commutators.
In terms of \emph{commutator width}, the question is
whether the commutator width of every 
non-abelian finite simple group is~$1$.
This question still remains open.
However, using the Classification of Finite Simple Groups,
it was shown by John Wilson that there exists a (not found explicitly) 
common upper bound on the commutator widths of all finite simple groups
(see~\cite{Wilson:1996:fogt}).

In 1977, Martin Isaacs (see \cite{Isaacs:1977:ccs})
noted that no simple group, finite or infinite,
was known to have commutator width greater than~$1$.
In 1999, Valerij Bardakov posed the following question
(see Problem 14.13 in \cite{kn14:1999:knupgt-eng}):
\begin{quote}
Does there exist a (finitely presented) simple group of infinite
commutator width?
\end{quote}

Simple groups of infinite commutator width, 
realised as groups of certain surface diffeomorphisms,
appeared in 
\cite{BargeGhys:1992:cEM-fr,GambaudoGhys:2004:cds} in 1992.
The infinity of the commutator width is established there by 
constructing \emph{nontrivial homogeneous quasi-morphisms}.
However, the question of existence of \emph{finitely generated\/} 
simple groups of commutator width greater than $1$ seems to have been
open until now.

In this paper it is shown that presentations 
(by generators and defining relations)
of finitely generated infinite simple groups 
of infinite commutator width, as well as of large finite commutator width, 
can be constructed using methods of small-cancellation theory.
This approach is rather flexible and can yield groups with various
additional properties.%
\footnote{In particular, groups can be constructed so that to admit
no nontrivial homogeneous quasi-morphisms
(it suffices to make the \emph{stable commutator length\/} of
each element of a group equal $0$).
This is not actually done in the paper.}

\begin{definition}
The \emph{commutator\/} of two group elements $x$ and $y$,
denoted $[x,y]$, is $xyx^{-1}y^{-1}$.
The \emph{commutator length\/} of an element $g$ of the derived
subgroup of a group $G$,
denoted $\cl_G(g)$, is the minimal $n$ such that there exist elements
$x_1$, \dots, $x_n$, $y_1$, \dots, $y_n$ in $G$ such that
$g=[x_1,y_1]\dots[x_n,y_n]$.
The commutator length of the identity element is $0$.
The \emph{commutator width\/} of a group $G$, denoted $\cw(G)$,
is the maximum of the commutator lengths of the elements of its
derived subgroup $[G,G]$.
\end{definition}

\begin{definition}
The \emph{conjugate\/} of a group element $g$ by a group element $h$,
denoted $g^h$, is $hgh^{-1}$.
A nontrivial group $G$ is called \emph{$n$-boundedly simple\/} if
for every two nontrivial elements $g,h\in G$,
the element $h$ is the product of $n$ or fewer
conjugates of $g^{\pm1}$, \ie,
$$
(\exists m\le n)\,
(\exists \sigma_1,\dots,\sigma_m\in\{\pm1\})\,
(\exists x_1,\dots,x_m\in G)\,
(g=(h^{\sigma_1})^{x_1}\dots(h^{\sigma_m})^{x_m}).
$$
A group $G$ is called \emph{boundedly simple\/} if
it is $n$-boundedly simple for some natural~$n$.
\end{definition}


Every boundedly simple group is simple, but the converse
is not generally true
(\eg, for an infinite alternating group).

\begin{remark}
A group is boundedly simple if and only if
each of its ultrapowers is simple.
If a group is $n$-boundedly simple, then all its ultrapowers are
$n$-boundedly simple.
\end{remark}

\begin{theorem}
\label{theorem:bsglcw}
For every natural\/ $n$\textup,
there exists a torsion-free\/ $2$-generated simple group\/ $G$
with a rank-\/$2$ free subgroup\/ $H$ such that\/\textup:
\begin{enumerate}
\item
	for every\/ $g\in G$ and every\/ $x\in G\setminus\{1\}$\textup,
	there exist\/ $y_1$\textup{, \dots,} $y_{2n+2}$ in\/ $G$ such that\/
	$g=x^{y_1}\dots x^{y_{2n+2}}$\textup; and
\item
	for every\/ $h\in H\setminus\{1\}$ and for every\/ $m\ge2n$\textup, 
	$\cl_G(h^m)>n$
\end{enumerate}
\textup(in particular\textup, $G$ is\/ $(2n+2)$-boundedly simple\textup,
and\/ $n+1\le\cw(G)\le2n+2$\textup{).}
Moreover\textup, there exists such a group\/ $G$ with decidable
word and conjugacy problems\textup.
\end{theorem}

Note that Theorem~\ref{theorem:bsglcw} improves the result of
Theorem~2 in~\cite{Muranov:2005:dsmcbgbsg}.

\begin{theorem}
\label{theorem:sgicw}
There exists a torsion-free\/ $2$-generated simple group\/ $G$
with a rank-\/$2$ free subgroup\/ $H$ such that
for every\/ $h\in H\setminus\{1\}$\textup,
$$
\lim_{n\to+\infty}\cl_G(h^n)=+\infty
$$
\textup(in particular\textup, $G$ has 
infinite commutator width\/\textup{).}
Moreover\textup, there exists such a group\/ $G$ with decidable
word and conjugacy problems\textup.
\end{theorem}

The theorems are proved by providing examples of groups
which satisfy the required properties.
These groups are presented by generators and defining relations
in Section~\ref{section:group_presentations}.
(The constructed presentations are recursive,
as follows from the proof of Proposition~\ref{proposition:Gn_G_swcp}.)

The properties of simplicity or bounded simplicity for the constructed
groups follow directly from the imposed relations, but the
existence of free non-cyclic subgroups and estimates on the commutator
lengths of their elements are obtained through a nontrivial
analysis of van Kampen diagrams on spheres with handles.

To show that the commutator length of a given element $g$ of
a constructed group is greater than $n$, it is proved that
if $\Delta$ is a van Kampen diagram on a sphere with handles and a hole
such that some group word representing the element $g$ ``reads''
on the boundary of $\Delta$,
then the number of handles is greater than $n$.
This is done by assuming that the number of handles
is not greater than $n$, which gives a lower bound on the
Euler characteristic, and coming to a contradiction.

The contradiction is obtained as follows.
The hole in the diagram is covered with an extra face
so as to make the diagram closed.
Some arcs of the diagram are selected and distributed
among the faces.
This is done in such a manner that the sum of the lengths
of the arcs associated to each face is small,
significantly less than half of the perimeter of that face,
but almost all edges of the diagram lie on selected arcs.
This eventually leads to a contradiction with the fact that
in a closed diagram the number of edges is half the sum
of the perimeters of faces.

Certainly, the least obvious part of the proof is
distributing ``almost all'' edges among the faces,
while associating ``few'' edges to each face.
For this purpose the group presentations are constructed
with small-cancellation-type conditions.
These conditions allow one to choose a system of selected arcs
in a specific way.
The selected arcs together cover almost all edges of the diagram,
and are short relative to the perimeters of incident faces.
Using the bound on the Euler characteristic of the diagram
(determined by the number of handles),
and one combinatorial lemma by Philip Hall,
the selected arcs can be distributed among the faces so as
to have ``very few'' relatively ``short'' and
``a few'' relatively ``very short''
arcs associated to each face.

The approach used in this paper is similar to that of
\cite{Muranov:2005:dsmcbgbsg}.
There are improvements and generalizations which allow one
to obtain better estimates and
to deal with diagrams on arbitrary surfaces.
\par



\section{Construction of the groups}
\label{section:group_presentations}




Let $\mathfrak A=\{a,b\}$ be a two-letter alphabet.
Choose recursive sequences of positive numbers
$\{\lambda_n\}_{n\in\mathbb N}$ and $\{\mu_n\}_{n\in\mathbb N}$
so that for every $n\in\mathbb N$,
\begin{equation}
\label{display:main_inequality}
2\lambda_n+(14n+8)\mu_n+\frac{2n+1}{4n+4}<\frac{1}{2}.
\end{equation}
(This inequality is to be used in the proofs of some of the properties.)
Let for definiteness
\begin{equation}
\label{display:lambda_and_mu.def}
\lambda_n=\frac{1}{20n+20}\qquad\text{and}\qquad
\mu_n=\frac{1}{(14n+8)(8n+8)}.
\end{equation}

\begin{remark}
The values of $\lambda_n$ and $\mu_n$ satisfying the inequality
\thetag{\ref{display:main_inequality}}
can be chosen by the \emph{lowest parameter principle\/}
of Alexander Ol'shanskii, \eg, $1/n\succ\mu_n\succ\lambda_n$,
which means that for every $n$, if $\mu_n$ is sufficiently small,
then the desired inequality holds for all sufficiently small $\lambda_n$
(see \cite{Olshanskii:1989:gosg-rus,Olshanskii:1991:gdrg-eng}).
\end{remark}



\subsection{Boundedly simple group of large commutator width}
\label{subsection:bsglcw.group_presentations}

Take an arbitrary natural number $n$.

Let $(v_1,w_1)$, $(v_2,w_2)$, $(v_3,w_3)$, \dots\ 
be a list of all ordered pairs of reduced group words
over $\mathfrak A$.
Moreover, let the function $i\mapsto(v_i,w_i)$ be recursive.

Let $\{u_{ij}\}_{i\in\mathbb N; j=1,\dots,2n+2}$ be
a recursive \emph{indexed family\/} of reduced group words,
and $z_1$, $z_2$ be two cyclically reduced group words
over $\mathfrak A$ such that:
\begin{enumerate}
\item
    for every $i\in\mathbb N$,
    \begin{enumerate}
    \item
        $\abs{u_{i,1}}=\abs{u_{i,2}}=\dots=\abs{u_{i,2n+2}}\ge i$, and
    \item
        $\lambda_n(4n+4)\abs{u_{i,1}}\ge\abs{v_i}+(2n+2)\abs{w_i}$;
    \end{enumerate}
\item
    the family $\{u_{ij}\}_{i\in\mathbb N; j=1,\dots,2n+2}$
    satisfies the following small-can\-cellation condition:
    if $u_{i_1j_1}^{\sigma_1}=p_1sq_1$ and
    $u_{i_2j_2}^{\sigma_2}=p_2sq_2$
    (here $\sigma_1,\sigma_2\in\{\pm1\}$),
    then either
    $$
    (i_1,j_1,\sigma_1,p_1,q_1)=(i_2,j_2,\sigma_2,p_2,q_2),
    $$
    or
    $$
    \mu_n(4n+4)\abs{u_{i_1j_1}}\ge\abs{s}\le\mu_n(4n+4)\abs{u_{i_2j_2}};
    $$
\item
    $z_1$ starts and ends with $a^{+1}$, and
    $z_2$ starts and ends with $b^{+1}$
    (hence, if $t(x,y)$ is an arbitrary reduced group word
    over $\{x,y\}$, then substituting $z_1$ for $x$ and $z_2$ for $y$
    yields a reduced groups word $t(z_1,z_2)$ over $\mathfrak A$);
\item
    if $s$ is a common subword of $u_{ij}$
    and of the concatenation of several copies of
    $z_1^{\pm1}$ and $z_2^{\pm1}$, then
    $$
    \abs{s}\le\mu_n(4n+4)\abs{u_{ij}}.
    $$
\end{enumerate}

For example, $z_1$, $z_2$,
and $\{u_{ij}\}_{i\in\mathbb N; j=1,\dots,2n+2}$
may be defined as follows:
\begin{equation}
\begin{gathered}
z_1=a^2,\qquad z_2=b^2,\\
u_{ij}=\prod_{k=4(14n+8)(j-1)+1}^{4(14n+8)j}
a^{k}b^{2(\abs{u_{i-1,1}}+\abs{v_i}+\abs{w_i})+4(14n+8)(2n+2)+1-k},
\end{gathered}
\label{display:u_and_z.bsglcw.def}
\end{equation}
where in the case $i=0$, the summand $\abs{u_{i-1,1}}$
shall be replaced with $0$.
(Here multiplication in $\Pi$-notation is understood in the
usual left-to-right sense, \eg, $\prod_{i=1}^3A_i=A_1A_2A_3$,
and not $A_3A_2A_1$.)

For every $i\in\mathbb N$, let
$$
r_i=w_i^{u_{i,1}}\dots w_i^{u_{i,2n+2}}v_i^{-1}
$$
where $w_i^{u_{i,1}}\dots w_i^{u_{i,2n+2}}v_i^{-1}$
denotes the concatenation of the group words
$u_{i,1}$, $w_i$, $u_{i,1}^{-1}$,
$u_{i,2}$, $w_i$, $u_{i,2}^{-1}$, \dots,
$u_{i,2n+2}$, $w_i$, $u_{i,2n+2}^{-1}$, and $v_i^{-1}$
in this order.

Now inductively construct a group presentation
\GP{\mathfrak A}{\mathcal R_n}
as follows.
Start with $\mathcal R_n^{(0)}=\varnothing$.
On step number $i$ ($i\in\mathbb N$),
if the relation `$w_i=1$' is a consequence of the relations
`$r=1$', $r\in\mathcal R_n^{(i-1)}$,
then define $\mathcal R_n^{(i)}=\mathcal R_n^{(i-1)}$;
otherwise, define $\mathcal R_n^{(i)}=\mathcal R_n^{(i-1)}\cup\{r_i\}$.
Finally, let
$$
\mathcal R_n=\bigcup_{i\in\mathbb N}\mathcal R_n^{(i)}.
$$

Let $G_n$ be the group presented by \GP{\mathfrak A}{\mathcal R_n}.


\subsection{Simple group of infinite commutator width}
\label{subsection:sgicw.group_presentations}

Let $w_1$, $w_2$, $w_3$, \dots\ 
be a recursive list of all reduced group words over $\mathfrak A$.

Let $\{u_{ij}\}_{i\in\mathbb N; j=1,\dots,4i+4}$ be
a recursive indexed family of reduced group words,
and $z_1$, $z_2$ be two cyclically reduced group words
over $\mathfrak A$ such that:
\begin{enumerate}
\item
    for every $i\in\mathbb N$,
    \begin{enumerate}
    \item
        $\abs{u_{i,1}}=\abs{u_{i,2}}=\dots=\abs{u_{i,4i+4}}$,
    \item
        $\lambda_i(4i+4)\abs{u_{i,1}}\ge1+(2i+2)\abs{w_i}$,
    \item
        $\lambda_i(4i+4)\abs{u_{i,1}}
        \le\lambda_{i+1}(4(i+1)+4)\abs{u_{i+1,1}}$,
    \item
        $\mu_i(4i+4)\abs{u_{i,1}}\le\mu_{i+1}(4(i+1)+4)\abs{u_{i+1,1}}$,
    \item
        $\abs{u_{i,1}}\le\abs{u_{i+1,1}}$, and
    \item
        $\mu_i(4i+4)\abs{u_{i,1}}\ge i$;
    \end{enumerate}
\item
    the family $\{u_{ij}\}_{i\in\mathbb N; j=1,\dots,4i+4}$
    satisfies the following condition:
    if $u_{i_1j_1}^{\sigma_1}=p_1sq_1$ and
    $u_{i_2j_2}^{\sigma_2}=p_2sq_2$
    ($\sigma_1,\sigma_2\in\{\pm1\}$),
    then either
    $$
    (i_1,j_1,\sigma_1,p_1,q_1)=(i_2,j_2,\sigma_2,p_2,q_2),
    $$
    or
    $$
    \mu_{i_1}(4i_1+4)\abs{u_{i_1j_1}}\ge\abs{s}
    \le\mu_{i_2}(4i_2+4)\abs{u_{i_2j_2}};
    $$
\item
    $z_1$ starts and ends with $a^{+1}$, and
    $z_2$ starts and ends with $b^{+1}$;
\item
    if $s$ is a common subword of $u_{ij}$
    and of the concatenation of several copies of
    $z_1^{\pm1}$ and $z_2^{\pm1}$, then
    $$
    \abs{s}\le\mu_i(4i+4)\abs{u_{ij}}.
    $$
\end{enumerate}

For example, $z_1$, $z_2$,
and $\{u_{ij}\}_{i\in\mathbb N; j=1,\dots,4i+4}$
may be defined as follows:
\begin{equation}
\begin{gathered}
z_1=a^2,\qquad z_2=b^2,\\
u_{ij}=\prod_{k=4(14i+8)(j-1)+1}^{4(14i+8)j}
a^{k}b^{2(\abs{u_{i-1,1}}+\abs{w_i})+4(14i+8)(4i+4)+1-k},
\end{gathered}
\label{display:u_and_z.sgicw.def}
\end{equation}
where in the case $i=0$, the summand $\abs{u_{i-1,1}}$
shall be replaced with $0$.

For every $i\in\mathbb N$, let
$$
r_{i,1}=w_i^{u_{i,1}}\dots w_i^{u_{i,2i+2}}a^{-1},\qquad
r_{i,2}=w_i^{u_{i,2i+3}}\dots w_i^{u_{i,4i+4}}b^{-1}.
$$

Now inductively construct a group presentation
\GP{\mathfrak A}{\mathcal R_\infty}
as follows.
Start with $\mathcal R_\infty^{(0)}=\varnothing$.
On step number $i$,
if the relation `$w_i=1$' is a consequence of the relations
`$r=1$', $r\in\mathcal R_\infty^{(i-1)}$,
then define $\mathcal R_\infty^{(i)}=\mathcal R_\infty^{(i-1)}$;
otherwise, define
$\mathcal R_\infty^{(i)}
=\mathcal R_\infty^{(i-1)}\cup\{r_{i,1},r_{i,2}\}$.
Finally, let
$$
\mathcal R_\infty=\bigcup_{i\in\mathbb N}\mathcal R_\infty^{(i)}.
$$

Let $G_\infty$ be the group presented by 
\GP{\mathfrak A}{\mathcal R_\infty}.


\bigskip

Note that:


\begin{proposition}
\label{proposition:Gn_bs}
For every natural\/ $n$\textup,
if the group\/ $G_n$ is nontrivial\textup,
then it is\/ $(2n+2)$-boundedly simple\textup.
Moreover\textup, for every\/ $g\in G_n$
and every\/ $x\in G_n\setminus\{1\}$\textup,
there exist\/ $y_1,\dots,y_{2n+2}\in G_n$ such that\/
$g=x^{y_1}\dots x^{y_{2n+2}}$\textup.
\end{proposition}



\begin{proposition}
\label{proposition:G_s}
If the group\/ $G_\infty$ is nontrivial\textup, then is simple\textup.
\end{proposition}



Other properties of these groups shall be established in
Section~\ref{section:proof_theorems}.
\par




\section{Combinatorial complexes, maps, and van Kampen diagrams}
\label{section:combinatorial_complexes_van_kampen_diagrams}


\subsection{Combinatorial complexes}
\label{subsection:combinatorial_complexes.ccvkd}

The purpose of introducing \emph{combinatorial cell complexes\/}
is to model CW-complexes by combinatorial objects which
preserve most of the combinatorial theory without
need for deep topological proofs.

Combinatorial cell complexes described here are equivalent to a
particular type of \emph{cone categories\/} defined in
\cite{McCammond:2000:gsct}.
They should not be confused with \emph{cell categories\/}
defined therein because, for example, in a cell category a $2$-cell
cannot have less than $3$ corners, but in a combinatorial cell complex
$2$-cells with just $1$ or $2$ corners are allowed.

Cone categories provide a perfect algebraic alternative to
CW-com\-plex\-es.
Unfortunately, despite the beauty of their concise and purely algebraic
definition, the language of cone categories seems less appealing to
geometric intuition than the language of CW-complexes.
Some of the natural geometric operations,
such as ``cutting and pasting,'' are harder to visualise
when thinking about cone complexes
instead of their geometric realisations.
The author chooses to define and use \emph{combinatorial cell complexes\/}
whose analogy to CW-complexes is more evident, but he believes that it is
possible to rewrite all the statements, proofs, and definitions in
this paper in terms of cone categories.

The following terms shall all be used as synonyms:
an \emph{$i$-dimensional combinatorial cell complex\/} is the same as
a \emph{combinatorial\/ $i$-complex}, or just an \emph{$i$-complex}.
Combinatorial $0$-complexes and $1$-complexes shall be viewed as
particular cases of $2$-complexes.
Only $0$-, $1$-, and $2$-complexes will be used,
and only their definitions are discussed.
Thus all combinatorial complexes may be assumed $2$-dimensional.

Some terms whose meaning in relation to combinatorial complexes
can be unambiguously inferred from their meaning in relation to
CW-complexes or simplicial complexes may be used without definition
(\eg, connectedness, a link of a vertex, etc.).

The definition of combinatorial complexes here is very similar to
that in \cite{Muranov:2005:dsmcbgbsg}.

A (combinatorial) \emph{$0$-complex\/} $A$ is a
$3$-tuple $( A(0),\varnothing,\varnothing)$
where $A(0)$ is an arbitrary set.
Elements of $A(0)$ are called \emph{$0$-cells}, or \emph{vertices},
of~$A$.

A $0$-complex with exactly $2$ vertices shall be called
a \emph{combinatorial\/ $0$-sphere}.

A \emph{morphism\/} $\phi$ from a $0$-complex $A$ to a $0$-complex $B$
is a $3$-tuple $(\phi(0),\varnothing,\varnothing)$ where $\phi(0)$
is an arbitrary function $A(0)\to B(0)$.
If $A$, $B$, and $C$ are $0$-complexes, and $\phi\!:A\to B$ and
$\psi\!:B\to C$ are morphisms, then
the product $\psi\phi\!:A\to C$ is defined naturally:
$(\psi\phi)(0)=\psi(0)\circ\phi(0)$.
A morphism $\phi\!:A\to B$ is called an \emph{isomorphism\/}
of $A$ with $B$ if
there exists a morphism $\psi\!:B\to A$
such that $\psi\phi$ is the \emph{identity morphism\/}
of the complex $A$ and $\phi\psi$ is the identity morphism
of the complex~$B$.

A (combinatorial) \emph{$1$-complex\/} $A$ is a $3$-tuple
$( A(0),A(1),\varnothing)$ such that:
\begin{enumerate}
\item
    $A(0)$ is an arbitrary set;
\item
    $A(1)$ is a set of ordered triples of the form $(i,E,\alpha)$
    where $E$ is a combinatorial $0$-sphere and $\alpha$
    is a morphism of $E$ to the $0$-complex 
    $( A(0),\varnothing,\varnothing)$,
    such that the function $( i,E,\alpha)\mapsto i$ is injective
    on~$A(1)$.
\end{enumerate}
Elements of $A(0)$ are called \emph{$0$-cells},
or \emph{vertices}, of $A$,
elements of $A(1)$ are called \emph{$1$-cells}, or \emph{edges}.
The $0$-complex $( A(0),\varnothing,\varnothing)$ is called 
the \emph{$0$-skeleton\/} of $A$ and is denoted~$A^0$.

If $e=( i,E,\alpha)$ is an edge of a $1$-complex $A$, then
$i$ is called the \emph{index\/} of $e$, $E$ is called the
\emph{characteristic boundary\/} of $e$ and shall be denoted by $\dot e$,
and $\alpha$ is called the \emph{attaching morphism\/} of~$e$.
(The purpose of ``indexing'' $1$- and $2$-cell is to allow
distinct cells with identical characteristic boundaries and
attaching morphisms.)

Combinatorial $1$-complexes may sometimes be called \emph{graphs}.

A $1$-complex which ``looks like'' a circle
shall be called a \emph{combinatorial\/ $1$-sphere},
or a \emph{combinatorial circle}.
More precisely, a combinatorial circle is a
finite \emph{connected\/} non-empty graph in which the 
\emph{degree\/} of every vertex is~$2$.

If $A$ and $B$ are $1$-complexes, then a \emph{morphism\/} $\phi\!:A\to B$
is a $3$-tuple $(\phi(0),\phi(1),\varnothing)$ such that:
\begin{enumerate}
\item
    $\phi^0=(\phi(0),\varnothing,\varnothing)$ 
    is a morphism of $A^0$ to $B^0$;
\item
    $\phi(1)$ is a function on $A(1)$ which maps
    each $e=( i,E,\alpha)\in A(1)$
    to an ordered pair $( e',\xi)$ such that
    $e'=( i',E',\alpha')\in B(1)$,
    $\xi$ is an isomorphism of $E$ with $E'$,
    and $\phi^0\alpha=\alpha'\xi$.
\end{enumerate}
Multiplication of morphisms
of $1$-complexes is defined naturally.
For example, if $A$, $B$, $C$ are $1$-complexes,
$\phi$ and $\psi$ are morphisms,
$\phi\!:A\to B$, $\psi\!:B\to C$,
$e=( i,E,\alpha)$ is an edge of $A$,
$e'=( i',E',\alpha')$ is an edge of $B$,
$e''=( i'',E'',\alpha'')$ is an edge of $C$,
$\phi(1)(e)=( e',\xi)$,
$\psi(1)(e')=( e'',\zeta)$,
then $(\psi\phi)(1)(e)=( e'',\zeta\xi)$
(note that $\zeta\xi$ is an isomorphism of the characteristic boundary of
$e$ with that of~$e''$).
\emph{Isomorphisms\/} of $1$-complexes are defined in the natural way.

A (combinatorial) \emph{$2$-complex\/} $A$ is a $3$-tuple
$( A(0),A(1),A(2))$ such that:
\begin{enumerate}
\item
    $( A(0),A(1),\varnothing)$ is a $1$-complex, called the
    \emph{$1$-skeleton\/} of $A$ and denoted~$A^1$;
\item
    $A(2)$ is a set of ordered triples of the form $( i,F,\alpha)$
    where $F$ is a combinatorial $1$-sphere and $\alpha$
    is a morphism of $F$ to $A^1$,
    such that the function $( i,F,\alpha)\mapsto i$ is injective
    on~$A(2)$.
\end{enumerate}
Elements of $A(0)$ are called \emph{$0$-cells}, or \emph{vertices},
of the complex $A$,
elements of $A(1)$ are called \emph{$1$-cells}, or \emph{edges},
elements of $A(2)$ are called \emph{$2$-cells}, or \emph{faces}.

If $f=( i,F,\alpha)$ is a face of a $2$-complex $A$, then
$i$ is called the \emph{index\/} of $f$, $F$ is called the
\emph{characteristic boundary\/} of $f$ and shall be denoted by $\dot f$,
and $\alpha$ is called the \emph{attaching morphism\/} of~$f$.

If $A$ and $B$ are $2$-complexes, then a \emph{morphism\/} $\phi\!:A\to B$
is a $3$-tuple $(\phi(0),\phi(1),\phi(2))$ such that:
\begin{enumerate}
\item
    $\phi^1=(\phi(0),\phi(1),\varnothing)$ is a morphism 
    of $A^1$ to $B^1$;
\item
    $\phi(2)$ is a function on $A(2)$ which maps
    each $f=( i,F,\beta)\in A(2)$ 
    to an ordered pair $( f',\xi)$ such that
    $f'=( i',F',\beta')\in B(2)$,
    $\xi$ is an isomorphism of $F$ with $F'$,
    and $\phi^1\beta=\beta'\xi$.
\end{enumerate}
Products of morphisms of $2$-complexes are defined analogously
to the case of $1$-complexes.
The notion of \emph{isomorphism\/} for $2$-complexes is the natural one.

The \emph{empty\/} combinatorial complex, a \emph{finite\/} complex,
a \emph{subcomplex}, etc., are defined naturally.

Any combinatorial complex $C$ gives rise to a CW-complex, called
its \emph{geometric realisation}, which is unique up to
isomorphism (homeomorphism preserving the cellular structure).
A geometric realisation is constructed as follows:

Let $C$ be a combinatorial $2$-complex.
Construct a geometric realisation $X^0$ of $C^0$ by 
imposing a structure of a $0$-dimensional CW-complex on an arbitrary set 
which is in bijective correspondence with $C(0)$.
Next, construct a geometric realisation $X^1$ of $C^1$ by attaching 
$1$-cells to $X^0$ as follows.
For every $e\in C(1)$, construct a geometric realisation $S$ of $\dot e$
(a $2$-point $0$-dimensional CW-complex);
take $B$ to be a cone over $S$ viewed as a topological space 
(homeomorphic to a closed interval);
attach $B$ to $X^0$ via the continuous function $S\to X^0$ induced
by the attaching morphism of $e$;
this procedure yields a $1$-cell in $X^1$ for each $1$-cell of $C$.
Finally, construct a geometric realisation $X=X^2$ of $C$ by attaching
$2$-cells to $X^1$ in correspondence with $2$-cells of $C$ 
in the same manner as how $1$-cells were attached to~$X^0$.

The class of CW-complexes that are geometric realisations of
combinatorial complexes is quite narrow, and does not even
include all \emph{transverse\/} $2$-dimensional CW-complexes.
For example, a geometric realisation of a $2$-complex cannot
have a $2$-cell attached to a single $0$-cell.

Every morphism $\phi$ from a combinatorial $2$-complex $A$
to a combinatorial $2$-complex $B$
naturally determines functions $A(i)\to B(i)$, $i=0,1,2$.
By abuse of terminology, if $e$ is an $i$-cell of $A$ and
$e'$ is its image under the function $A(i)\to B(i)$ induced by a morphism
$\phi\!:A\to B$,
then $e'$ shall be called the \emph{image\/} of $e$ under $\phi$, or 
the \emph{$\phi$-image\/} of~$e$.

If $e=(i,E,\alpha)$ is an edge of a $2$-complex $C$,
then the vertices of the characteristic boundary $\dot e=E$ of $e$
are called the \emph{ends\/} of $e$.
The images of the ends of $e$ under the attaching morphism $\alpha$ of $e$
are called the \emph{end-vertices\/} of~$e$.
An edge $e$ is \emph{incident\/} to a vertex $v$ if
$v$ is an end-vertex of~$e$.
An edge with only $1$ end-vertex is called a \emph{loop}.

If $f$ is a face of a $2$-complex $\Phi$,
then the vertices and edges of the characteristic boundary of $f$ are
called the \emph{corners\/} and \emph{sides\/} of $f$, respectively.
The images of the corners and sides of a face $f$ under
the attaching morphism of $f$ are called
\emph{corner-vertices\/} and \emph{side-edges\/} of $f$, respectively,
and are said to be \emph{incident\/} to~$f$.


The \emph{size\/} of a face $f$ is the number of its sides 
(which is equal to the number of its corners).

If $v$ is a vertex of a complex $C$, then the number of edges incident
to $v$ ``counted with multiplicity'' is called the \emph{degree\/} of $v$.
More precisely, the degree of a vertex $v$ in a complex $C$
is the total number of ordered pairs $(e,x)$
where $e$ is an edge of $C$, $x$ is an end of $e$, and $v$ is the
image of $x$ under the attaching morphism of~$e$.

The \emph{link\/} of a vertex $v$ in a combinatorial complex $C$,
denoted $\vtxlnk_C v$, or $\vtxlnk v$, is a combinatorial complex
which can be viewed as
``the boundary of a nice small neighborhood of $v$.''
A precise definition is not given here since
the meaning of the term in the present context can be easily
inferred from its meaning in context of
simplicial complexes or CW-complexes.%
\footnote{In the language of cone categories,
the link of a vertex is the full subcategory of the co-slice
category of that vertex obtained by removing the initial object.}
Note that the link of a vertex $v$ in a $2$-complex is a
$1$-complex whose vertices are in bijective correspondence with
all the ordered pairs $(e,x)$ where $e$ is an edge incident to $v$ and
$x$ is an end of $e$ which is mapped to $v$ by the attaching morphism 
of $e$,
and whose edges are in bijective correspondence with all the
ordered pairs $(f,x)$ where $f$ is a face incident to $v$
and $x$ is a corner of $f$ which is mapped to $v$ 
by the attaching morphism of~$f$.

An \emph{orientation\/} of an edge $e$ is a function from the set of
ends of $e$ to $\mathbb Z$ which maps one of the ends to $+1$ and the
other to $-1$.
An \emph{oriented edge\/} is an edge together with an orientation.
The end-vertex of an oriented edge $e$ 
that is the image of the ``negative'' end is called the \emph{tail-vertex},
or the \emph{initial vertex}, of $e$.
The end-vertex that is the image of the ``positive''
end is called the \emph{head-vertex}, or the \emph{terminal vertex}, of
the oriented edge;
on figures it is indicated with an arrowhead.
An oriented edge \emph{exits\/} its tail-vertex and \emph{enters\/} 
its head-vertex.

Consider a combinatorial circle $C$ and a function $f$
which chooses an orientation of each edge of~$C$.
The choice of orientations $f$ is called \emph{coherent\/}
if every vertex of $C$ is the tail-vertex of exactly one 
(and the head-vertex of exactly one) 
of the oriented edges obtained from edges of $C$ by 
assigning orientations according to $f$.
A coherent choice of orientations of all edges of $C$ is called 
an \emph{orientation\/} of $C$,
and $C$ together with one of its orientations is called an
\emph{oriented combinatorial circle}.

An \emph{orientation\/} of a face $f$ is an orientation
of the characteristic boundary of~$f$.
An \emph{oriented face\/} is a face together with an orientation.

Every edge and every face has exactly two \emph{opposite\/} orientations.
Two oriented edges (or faces) with the same underlying non-oriented
edge (or face) but with opposite orientations are called
mutually \emph{inverse}.

A \emph{path\/} is a non-empty finite sequence of alternating vertices and
oriented edges in which every oriented edge is 
immediately preceded by its tail-vertex and 
immediately succeeded by its head-vertex.

The \emph{length\/} of a path
$p=(v_0\,e_1,v_1,\dots,e_n,v_n)$ is $n$;
it is denoted by $\abs{p}$.
The vertices $v_1,\dots,v_{n-1}$ of this path are called
\emph{intermediate}.
A \emph{trivial path\/} is a path of length zero.
By abuse of notation, a path of the form $(v_1,e,v_2)$
shall be denoted by $e$,
and a trivial path $(v)$ shall be denoted by~$v$.

The \emph{inverse path\/} to a path $p$ is defined naturally and
is denoted by $p^{-1}$.
If the \emph{terminal vertex\/} of a path $p_1$ coincides with
the \emph{initial vertex\/} of a path $p_2$, then the \emph{product\/}
$p_1p_2$ is defined (naturally).
A path $s$ is an \emph{initial subpath\/} of a path $p$ if
$p=sq$ for some path $q$.
A path $s$ is a \emph{terminal subpath\/} of $p$ if
$p=qs$ for some~$q$.

A \emph{cyclic path\/} is a path whose terminal and initial vertices 
coincide.
A \emph{cycle\/} is the set of all \emph{cyclic shifts\/} of a
cyclic path.
The cycle represented by a cyclic path $p$ shall be denoted
by $\langle p\rangle$.
The \emph{length\/} of a cycle $c$, denoted by $\abs{c}$,
is the length of an arbitrary representative of $c$.
A \emph{trivial cycle\/} is a cycle of length zero.
A path $p$ is a \emph{subpath\/} of a cycle $c$ if
for some representative $r$ of $c$ and for some $n\in\mathbb N$,
$p$ is a subpath of $r^n$ (\ie, of the product of $n$ copies of~$r$).

A path is \emph{reduced\/} if
it has no subpath of the form $ee^{-1}$ where $e$ is an oriented edge.
A cyclic path is \emph{cyclically reduced\/} if
it is reduced and its first and last oriented edges 
are not mutually inverse.
(Trivial paths are cyclically reduced.)
A cycle is \emph{reduced\/} if
it consists of cyclically reduced cyclic paths.
A path is \emph{simple\/} if
it is nontrivial, reduced, and
every its intermediate vertex appears in it only once.
A cycle is \emph{simple\/} if
it consists of simple cyclic paths.

An \emph{oriented arc\/} of a complex $C$ is a simple path
whose all intermediate vertices have degree $2$ in $C$.
An \emph{oriented pseudo-arc\/} of a complex $C$ is a nontrivial
reduced path all intermediate vertices of which have degree $2$ in $C$.
A (non-oriented) \emph{arc\/}
is a pair of mutually inverse oriented arcs.
A (non-oriented) \emph{pseudo-arc\/}
is a pair of mutually inverse oriented pseudo-arcs.
Sometimes edges may be viewed as arcs,
and oriented edges as oriented arcs.

Note that if $C$ is a connected $1$-complex and 
is not a combinatorial circle, then every pseudo-arc of $C$ is an arc.
All pseudo-arcs under consideration will be pseudo-arcs in
characteristic boundaries of faces.

A (pseudo-)arc or an edge $u$ \emph{lies on\/} a path $p$ if
at least one of the oriented (pseudo-)arcs or oriented edges, respectively,
associated with $u$ is a subpath of $p$.
An arc shall be called \emph{free\/} if
none of its edges is incident to any face.

Any morphism of combinatorial complexes induces maps
of paths and cycles.
An arc $u$ is \emph{incident\/} to a face $f$ if
the associated oriented arcs are the images
of some paths in $\dot f$
under the map induced by the attaching morphism of~$f$.

\begin{definition}
A \emph{c-path}, \emph{c-cycle}, \emph{c-arc},
or \emph{c-pseudo-arc\/} of a face $f$ is a path,
a cycle, an arc, or a pseudo-arc in the characteristic boundary of $f$, respectively.
A \emph{c-edge\/} of a face is the same as its side.
\end{definition}

It shall be assumed that distinct faces always have disjoint sets
of corners and disjoint sets of sides, as well as that distinct edges
have disjoint sets of ends.
This will assure that every c-path is a c-path of exactly one face.
Moreover, it is convenient to assume that no other set-theoretic
complications happen;
in particular, no vertex can be simultaneously an edge.

To use combinatorial complexes effectively, a few operations on them
need to be defined, and some properties of these operations need
to be established.

The following operations are easy to define:
\begin{description}
\item[Removing a face]
    this operation is self-explanatory.
\item[Removing an arc]
    If $u$ is a free arc (or free edge) in a complex $C$,
    then to \emph{remove\/} $u$ from $C$ means
    to remove all edges and all \emph{intermediate\/} vertices
    of $u$ from~$C$.
\item[Attaching a face (along a cyclic path)]
    This is the operation inverse to removing a face.
    If $p$ is a cyclic path in a complex $C$, then
    to \emph{attach\/} a face $f$ along $p$ means to attach
    a new face $f$ in such a way that the image of some
    simple cyclic c-paths of $f$ be mapped to $p$ by the
    attaching morphism of~$f$.
\item[Attaching an arc (at a pair of vertices)]
    this is the operation inverse to removing an arc.
\end{description}

Next, there are operations whose geometric meaning is clear,
but whose precise combinatorial definition may be complicated.
Instead of giving precise definitions,
these operations are informally described here:
\begin{description}
\item[Dividing an edge by a vertex]
    To \emph{divide\/} an edge $e$ by a vertex $v$
    geometrically means to put a new vertex $v$ inside~$e$.
\item[Dividing a face by an edge]
    To \emph{divide\/} a face $f$ by an edge $e$ through its
    (not necessarily distinct) corners $x$ and $y$
    geometrically means to connect the corners $x$ and $y$ within $f$
    by a new edge~$e$.
\item[Pulling an edge into a face]
    To \emph{pull\/} an edge $e$ into a face $f$ through its corner $x$
    geometrically means to put a new vertex $v$ inside $f$ 
    and connect it within $f$ to the corners $x$ by a new edge~$e$.
\item[Merging two edges across a vertex]
    this is the operation inverse to dividing an edge by a vertex.
\item[Merging two faces across an edge]
    this is the operation inverse to dividing a face by an edge.
\item[Pushing an edge out of a face]
    this is the operation inverse to pull\-ing an edge into a face.
\end{description}

\begin{definition}
The complex obtained from a given $2$-complex $C$ by an arbitrary
(finite) sequence of operations of dividing edges by vertices,
dividing faces by edges, and pulling edges into faces is called
a \emph{subdivision\/} of $C$.
Two combinatorial $2$-complexes are called
\emph{geometrically equivalent\/} if 
they have isomorphic subdivisions.
\end{definition}

As shown in Section~9 of 
\cite{Olshanskii:1989:gosg-rus,Olshanskii:1991:gdrg-eng},
if the topological spaces of the geometric realisations of
two combinatorial complexes $A$ and $B$
are homeomorphic $2$-dimensional surfaces with or without boundaries,
then $A$ and $B$ are geometrically equivalent in the defined above sense.
Vice versa, it is obvious that geometrically equivalent 
combinatorial complexes have homeomorphic geometric realisations.

\begin{definition}
A \emph{combinatorial surface\/} is a non-empty combinatorial complex
in which every vertex link is either a combinatorial circle
or a \emph{combinatorial segment\/}
(a finite connected $1$-complex in which $2$ vertices have degree $1$
and all the others have degree $2$).
A combinatorial surface is \emph{closed\/} if the link of every vertex
is a combinatorial circle.
\end{definition}

Only combinatorial surfaces will be discussed, rather than
topological ones.
The adjective ``combinatorial'' shall be omitted for brevity.

Alternatively, a combinatorial surface may be defined as a
combinatorial complex whose geometric realisation is
a $2$-dimensional surface, with or without boundary.
Only finite combinatorial surfaces are discussed in this paper.

Consider a combinatorial surface $S$ and a function $f$
which chooses an orientation of each face of $S$.
Then the function $f$ induces an orientation on every side of
every face of $S$.
The attaching morphisms of faces carry the orientations of sides
of these faces over to the edges that are the images of these sides
(under the attaching morphisms).
The choice of orientations $f$ is called \emph{coherent\/}
if for every edge $e$ of $S$ that is the image of two distinct face sides,
the orientations of these two sides determined by $f$
induce (via the attaching morphisms) opposite orientations of~$e$.

\begin{definition}
A coherent choice of orientations of all faces of a combinatorial surface
$S$ is called an \emph{orientation\/} of $S$,
and $S$ together with an orientation is called an
\emph{oriented combinatorial surface}.
A combinatorial surface which can be oriented is called \emph{orientable}.
\end{definition}

It can be shown that a combinatorial surface is orientable if and only if
its geometric realisation is.

An orientable connected combinatorial surface has exactly two orientations.

Let a \emph{sample combinatorial disc\/} be a $2$-complex $C$ consisting
of $1$ vertex, $1$ edge, and $1$ face whose attaching morphism
is an isomorphism of its characteristic boundary with the $1$-skeleton
of~$C$.

\begin{definition}
A \emph{combinatorial disc\/} is an arbitrary $2$-complex
which is geometrically equivalent to a sample combinatorial disc.
\end{definition}

Let a \emph{sample combinatorial sphere\/} be a $2$-complex $C$
consisting of $1$ vertex, $1$ edge, and $2$ faces whose
attaching morphisms are isomorphisms of their characteristic boundaries
with the $1$-skeleton of~$C$.

\begin{definition}
A \emph{combinatorial sphere\/} is an arbitrary $2$-complex
which is geometrically equivalent to a sample combinatorial sphere.
\end{definition}

(Combinatorial discs and spheres are exactly those
$2$-com\-plexes whose geometric realisations are $2$-discs and $2$-spheres,
respectively.)

In a similar fashion, other combinatorial surfaces,
\eg, combinatorial tori, may be defined.
Such terms shall be used without further definitions.

\begin{definition}
A \emph{nontrivial singular combinatorial disc\/} is an arbitrary
$2$-com\-plex that can be obtained from a combinatorial sphere
by removing $1$ face
(or which can be turned into a combinatorial sphere by attaching $1$ face).
A \emph{trivial singular combinatorial disc\/} is a combinatorial complex
that consists of a single vertex.
\end{definition}


\begin{lemma}
\label{lemma:subcomplex__surface}
If\/ $C$ is a proper connected subcomplex
of a combinatorial surface\textup,
and\/ $C$ has at least\/ $1$ edge\textup,
then\/ $C$ can be turned into a connected combinatorial surface
by operations of attaching faces\textup.
If\/ $C$ is a proper subcomplex
of a connected combinatorial surface\textup,
then\/ $C$ is not a closed combinatorial surface\textup.
\end{lemma}


This lemma is not proved here because it is intuitively obvious, 
while its proof would probably be rather technical but hardly interesting.

\begin{definition}
The \emph{Euler characteristic\/} of a $2$-complex $C$
is denoted by $\chi_C$ and is defined by
$\chi_C=\norm{C(0)}-\norm{C(1)}+\norm{C(2)}$
where $\norm{C(i)}$ is the number of $i$-cells of $C$, $i=0,1,2$.
\end{definition}


\begin{lemma}
\label{lemma:positive_Euler_characteristic}
The maximal possible Euler characteristic of a closed connected
combinatorial surface is\/ $2$\textup, and among all
closed connected surfaces\textup,
only spheres have Euler characteristic\/ $2$\textup, and
only projective planes have Euler characteristic~$1$.
The maximal possible Euler characteristic of a
proper connected subcomplex of a combinatorial
surface is\/ $1$\textup, and
every such complex is a singular combinatorial disc\textup.
\end{lemma}


\begin{proof}
The first part of this lemma follows from the classification of compact
(or finite combinatorial) surfaces.
The second part follows from the first part together with
Lemma~\ref{lemma:subcomplex__surface}.
\end{proof}


Thus a combinatorial disc could be defined as
a non-closed connected finite combinatorial surface
of Euler characteristic $1$,
and a combinatorial sphere could be defined as
a connected finite combinatorial surface
of Euler characteristic~$2$.


\subsection{Maps}
\label{subsection:maps.ccvkd}

\begin{definition}
A \emph{nontrivial connected map\/} $\Delta$ consists of:
\begin{enumerate}
\item
    a finite connected combinatorial complex $C$;
\item
    a function that for every face $\Pi$ of $C$ chooses
    one of its simple cyclic c-paths,
    called the \emph{contour c-path}, or \emph{c-contour}, of $\Pi$;
\item
    a (possibly empty) indexed family of cyclic paths in $C$,
    called the \emph{contour paths}, or \emph{contours}, of~$\Delta$.
\end{enumerate}
The only requirement to this structure is that
a complex obtained from $C$ by attaching one new face
along each of the \emph{indexed contours\/} of $\Delta$ must be
a closed combinatorial surface.%
\footnote{It seems that a natural way to extend the notion of a map
would be to weaken this condition and require instead that in
the complex  obtained from $C$ by attaching faces along the
contours of $\Delta$, the links of all vertices are connected.
}
\end{definition}

An \emph{indexed contour\/} of a map is an element of the
indexed family of contours of the map together with its index.
One contour of a map may correspond to two distinct indexed contours.

The c-contour of a face $\Pi$ shall be denoted by $\ccntr\Pi$.
The image of the c-contour of a face $\Pi$ in $C^1$ is called
the \emph{contour path}, or \emph{contour}, of $\Pi$ and
shall be denoted by $\cntr\Pi$.
The cycle represented by the contour of a face $\Pi$ is called
the \emph{contour cycle\/} of $\Pi$.
Similarly, the cycle represented by a contour of a map $\Delta$
is called a \emph{contour cycle\/} of $\Delta$.
The contours of $\Delta$ shall be denoted by
$\cntr_1\Delta$, $\cntr_2\Delta$, $\cntr_3\Delta$,
et cetera.
If $\Delta$ has only one contour, then it can be alternatively
denoted by~$\cntr{}\Delta$.

\begin{definition}
A \emph{trivial map\/} is a combinatorial complex consisting
of a single vertex together with the trivial cyclic path in it
called its \emph{contour}.
\end{definition}

A \emph{map\/} in general consists of
a finite non-zero number of \emph{connected components},
each of which is either a nontrivial connected map, or a trivial map.

\begin{definition}
A map without contours is called \emph{closed}.
If $\Delta$ is a map without trivial connected components,
then a closed map obtained
from $\Delta$ by attaching new faces along its contours,
and choosing the c-contours of the new faces so that
the contours of $\Delta$ become the contours of the new faces,
is called a \emph{closure\/} of~$\Delta$.
\end{definition}

\begin{remark}
It is not possible to define a closure of a trivial map
in a similar way because no face in a combinatorial complex
can have boundary consisting of a single vertex.
\end{remark}

A map is closed if and only if its underlying complex
is a closed surface (see Lemma~\ref{lemma:subcomplex__surface}).

A closure of a map $\Delta$ is unique up to isomorphism.

The important notion of a \emph{submap\/} of a map shall be defined in
the subsequent subsection.

\begin{definition}
Two maps are called \emph{essentially isomorphic\/} if
there exists an isomorphism between their underlying
complexes which preserves the contours of the maps up to re-indexing
and/or replacing with cyclic shifts or cyclic shifts of the inverses.
\end{definition}

Two maps without trivial connected components are essentially isomorphic
if and only if there exists an isomorphism between their underlying
combinatorial complexes that extends to an isomorphism between
the underlying complexes of their closures.

\begin{definition}
A map is \emph{simple\/} if
all its contours are simple (cyclic) paths,
and distinct contours have no common vertices.
A map is \emph{semi-simple\/} if
in it every edge is incident to a face.
A map without faces is called \emph{degenerate}.
\end{definition}

\begin{definition}
A \emph{disc\/} map is either a trivial map, or any map with exactly
one contour which has a spherical closure.
An \emph{annular\/} map is an arbitrary map with exactly two contours
which has a spherical closure.
An \emph{elementary\/} map is
a spherical map with exactly $2$ faces,
whose $1$-skeleton is a combinatorial circle.
\end{definition}

The underlying $2$-complex of a disc map is a singular
combinatorial disc, and the underlying $2$-complex of simple
disc map is a combinatorial disc.


\begin{lemma}
\label{lemma:maps__positive_Euler_characteristic}
A connected map of Euler characteristic\/ $2$ is closed spherical\textup.
A connected map of Euler characteristic\/ $1$ is either disc\textup,
or\/ \textup(closed\/\textup) projective-planar\textup.
The maximal possible Euler characteristic of a connected map with\/
$n$ contours is\/~$2-n$\textup.
\end{lemma}


This lemma follows from Lemma~\ref{lemma:positive_Euler_characteristic}.

Every non-free arc of a map is either \emph{internal\/}
(is the image of two distinct c-arcs, and does not lie on any
contour cycle of the map),
or \emph{external\/}
(is the image of only one c-arc,
and lies on some contour cycle of the map).

\begin{definition}
A map is called \emph{contour-oriented\/} if
every oriented edge occurs in the contour of some face or
in some contour of the map.
\end{definition}

If $\Delta$ is a contour-oriented map and $\bar\Delta$ is its closure,
then the c-contours of the faces of $\bar\Delta$
induce an orientation of the underlying complex of~$\bar\Delta$.

It is convenient to have terms to express the idea that two
paths in the characteristic boundary of a face
``go in the same direction,'' and also to have a ``preferred direction''
in the boundary.
For that purpose let all nontrivial reduced c-path of each face of a map
be divided into \emph{positive\/} and \emph{negative\/}:

\begin{definition}
A c-path $p$ of a face $\Pi$ of a map is called
\emph{positive\/} if
it is a nontrivial subpath of $\ccntrcclof{\Pi}$.
A c-path $p$ is \emph{negative\/} if
$p^{-1}$ is positive.
\end{definition}


\subsection{Transformations of maps}
\label{subsection:transformations__maps.ccvkd}

\subsubsection{Removing a face}

If $\Pi$ is a face of a map $\Delta$, then
the \emph{submap\/} of $\Delta$ obtained by \emph{removing\/} $\Pi$
is the map $\Psi$ such that:
\begin{enumerate}
\item
    the underlying complex of $\Psi$ is obtained from
    the underlying complex of $\Delta$ by removing
    the face~$\Pi$,
\item
    the c-contours of faces of $\Psi$ are those inherited
    from $\Delta$, and
\item
    the indexed family of contours of $\Psi$ is obtained from the
    indexed family of contours of $\Delta$ by adding the contour of $\Pi$
    as a new indexed member, and possibly re-indexing the family.
\end{enumerate}

(\emph{Re-indexing\/} an indexed family means replacing the index set
with a new set of the same cardinality, and pre-composing the indexing
function with a bijection from the new index set onto the original
index set of the family.)

\subsubsection{Removing an arc}

If $u$ is a free arc of a map $\Delta$, then
a \emph{submap\/} of $\Delta$ obtained by \emph{removing\/} $u$
is a map $\Psi$
such that:
\begin{enumerate}
\item
    the underlying complex of $\Psi$ is obtained from
    the underlying complex of $\Delta$ by removing
    the arc $u$;
\item
    the c-contours of faces of $\Psi$ are those inherited
    from~$\Delta$;
\item
    the family of contours of $\Psi$ consists of, up to re-indexing,
    all the indexed contours of $\Delta$ that do not have
    common edges with $u$,
    together with additional $1$ or $2$ chosen as follows:
    \begin{enumerate}
    \item
        if removing $u$ from the underlying complex increases the number
        of connected components,
        then, first, take paths $p_1$, $p_2$, and $v$ such that:
        \begin{enumerate}
        \item
            $v$ and $v^{-1}$ are the oriented arcs
            associated with $u$, and
        \item
            $\langle vp_1v^{-1}p_2\rangle$
            is a contour cycle of $\Delta$,
        \end{enumerate}
        and second, assign new indices to
        some cyclic shift of $p_1^{\pm1}$ and some cyclic shift of $p_2^{\pm1}$,
        and take them as $2$ additional indexed contours of~$\Psi$;
    \item
        if removing $u$ from the underlying complex does not increase
        the number of connected components, and
        $u$ lies on only one contour cycle of $\Delta$
        (which implies that a closure of $\Delta$ is
        non-orientable),
        then, first, take paths $p_1$, $p_2$, and $v$ such that:
        \begin{enumerate}
        \item
            $v$ and $v^{-1}$ are the oriented arcs associated
            with $u$, and
        \item
            $\langle vp_1vp_2\rangle$ is a contour cycle
            of $\Delta$,
        \end{enumerate}
        and second, assign a new index to
        a cyclic shift of $(p_1p_2^{-1})^{\pm1}$,
        and take it as an additional indexed contour of~$\Psi$;
   \item
        if $u$ lies on two distinct contour cycles of $\Delta$,
        then, first, take paths $p_1$, $p_2$, and $v$,
        and indices $i$ and $j$ such that:
        \begin{enumerate}
        \item
            $v$ and $v^{-1}$ are the oriented arcs associated
            with $u$,
        \item
            $\langle vp_1\rangle=\langle\cntr_i\Delta\rangle$,
        \item
            either $\langle vp_2\rangle=\langle\cntr_j\Delta\rangle$,
            or $\langle p_2^{-1}v^{-1}\rangle
            =\langle\cntr_j\Delta\rangle$,
            and
        \item
            $i\ne j$,
        \end{enumerate}
        and second, assign a new index to
        a cyclic shift of $(p_1p_2^{-1})^{\pm1}$,
        and take it as an additional indexed contour of~$\Psi$.
    \end{enumerate}
\end{enumerate}

\begin{definition}
A map $\Psi$ is a \emph{submap\/} of a map $\Delta$ if
it can be obtained from $\Delta$ by operations of 
removing faces, removing free arcs, and removing connected components
(the last operation is self-explanatory).
\end{definition}



\begin{lemma}
\label{lemma:submap}
If\/ $\Gamma$ is a submap of a map\/ $\Delta$\textup,
and\/ $\bar\Delta$ is a closure of\/ $\Delta$\textup,
then\/\textup:
\begin{enumerate}
\item
    the underlying complex of\/ $\Gamma$ is a subcomplex of
    the underlying complex of\/~$\Delta$\textup;
\item
    the contours and c-contours of faces of\/ $\Gamma$ 
    are those inherited from\/~$\Delta$\textup;
\item
    for every contour\/ $q$ of\/ $\Gamma$ that is neither a contour
    of\/ $\Delta$\textup, nor the contour of a face of\/ $\Delta$\textup,
    there is\/ $n\in\mathbb N$ and
    there are subpaths\/ $p_0$\textup, $p_1$\textup{, \dots,} $p_n=p_0$
    of\/ $q$\textup,
    oriented arcs\/ $u_0$\textup, $u_1$\textup{, \dots,} $u_n=u_0$
    of\/ $\Delta$\textup,
    and c-paths\/
    $p_0'$\textup, $p_1'$\textup{, \dots,} $p_n'=p_0'$\textup,
    $u_0'$\textup, $u_1'$\textup{, \dots,} $u_n'=u_0'$\textup,
    $v_0'$\textup, $v_1'$\textup{, \dots,} $v_n'=v_0'$
    of faces of\/ $\bar\Delta$ that are not in\/ $\Gamma$\textup,
    such that\/\textup:
    \begin{enumerate}
    \item
        $\langle q\rangle=\langle p_0p_1\dots p_{n-1}\rangle$\textup,
    \item
        $u_0$\textup, $u_1$\textup{, \dots,} $u_n$
        are pairwise non-overlapping and maximal among oriented arcs of\/
        $\Delta$ that do not have common edges with\/~$\Gamma$\textup,
    \item
        the terminal vertices of\/
        $u_0$\textup, $u_1$\textup{, \dots,} $u_n$
        are vertices of\/~$q$\textup,
    \item
        for every\/ $i=0,\dots,n$\textup,
        the images of\/ $p_i'$\textup, $u_i'$\textup, and\/ $v_i'$ are\/
        $p_i$\textup, $u_i$\textup, and\/ $u_i^{-1}$\textup,
        respectively\textup,
        but\/ $v_i'\ne u_i^{\prime-1}$\textup,
    \item
        for every\/ $i=0,\dots,n-1$\textup,
        the product\/ $u_i'p_i'v_{i+1}'$ is a reduced
        c-path of a face of\/~$\bar\Delta$\textup.
    \end{enumerate}
\end{enumerate}
\end{lemma}


This lemma can be proved by induction on the number of operations
of removing (free) arcs used in obtaining $\Gamma$ from~$\Delta$.


\begin{lemma}
\label{lemma:subcomplex__map}
Every subcomplex of the underlying complex of
every map\/ $\Delta$ has a structure of a submap of\/ $\Delta$\textup,
which is unique up to essential isomorphism\textup.
\end{lemma}


The non-obvious part of this lemma is the ``uniqueness.''
It can be proved by showing that the set of essential isomorphism
classes of submaps of a given map containing a given subcomplex,
together with operations of removing faces, arcs,
and connected components, form a
\emph{confluent\/} and \emph{terminating rewriting system}.

For brevity, subcomplexes of the underlying complexes of maps shall be
called subcomplexes of the maps.


\subsubsection{Diamond move}

\begin{definition}
Two c-edges, or two oriented c-edges, of faces of a complex $C$ are called
\emph{contiguous\/} if
their images (in $C^1$) under the respective attaching morphisms coincide.
\end{definition}

\emph{Contiguity\/} of oriented c-edges is an equivalence relation.
Observe that any closed map is determined up to isomorphism by
c-contours of
its faces and the contiguity relation on oriented c-edges.
In other words, the contiguity relation tells how to glue faces together,
which together with a choice of c-contours of faces
``completely'' determines a closed map.
This observation allows one to define \emph{diamond moves\/}
on maps in terms of
changing the contiguity relation on oriented c-edges.

Consider an arbitrary \emph{closed\/} map $\Delta$.
Let $e_1$ and $e_2$ be two distinct oriented edges in $\Delta$
with a common terminal vertex $v$.

Choose a ``local orientation around $v$,'' \ie,
choose an orientation of the link of $v$ in $\Delta$.
Consider an arbitrary oriented c-edge $x$ of a face of $\Delta$
such that $v$ is the head-vertex of the image of $x$
(for example, the image of $x$ may be $e_1$ or $e_2$).
The ``positive'' end of $x$ naturally corresponds to
an end of some edge of $\vtxlnk v$,
which in turn is either ``positive'' or ``negative''
with respect to the chosen orientation of $\vtxlnk v$.
Call $x$ ``positive'' or ``negative'' accordingly.
(This terminology shall only be used in the subsequent definition.)
Thus every oriented edge entering $v$ has one ``positive''
and one ``negative'' pre-image under attaching morphisms.
Let $a$, $b$, $c$, and $d$ be, respectively, 
the ``positive'' and the ``negative'' pre-images of $e_1$,
and the ``positive'' and the ``negative'' pre-images of~$e_2$.

\begin{figure}\centering
\includegraphics{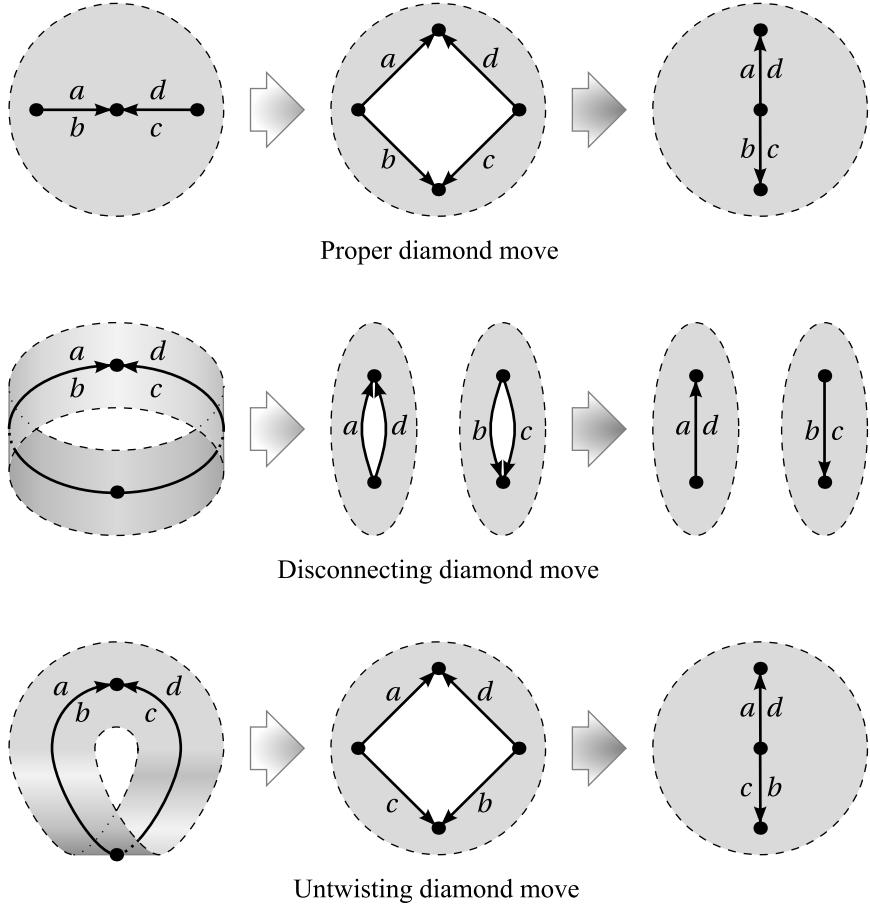}
\caption{Diamond Move.}
\label{figure:1}
\end{figure}

\begin{definition}
A map obtained from $\Delta$ by the \emph{diamond move\/}
along $e_1$ and $e_2$ is a (unique up to isomorphism)
closed map $\Gamma$ which has the same characteristic boundaries and
c-contours of faces, but in which $a$ is contiguous to $d$,
$b$ is contiguous to $c$, and the contiguity relation
on the oriented c-edges distinct from 
$a^{\pm1}$, $b^{\pm1}$, $c^{\pm1}$, $d^{\pm1}$ is the same as in $\Delta$
(see Fig.~\ref{figure:1}).
\end{definition}

Let $e_3$ and $e_4$ be the images in $\Gamma^1$ of 
$a$ and $d$, and $b$ and $c$, respectively.
There are natural bijections:
\begin{enumerate}
\item
    between the vertices of $\Delta$ not incident to $e_1$ and $e_2$,
    and the vertices of $\Gamma$ not incident to $e_3$ and~$e_4$;
\item
    between the edges of $\Delta$ distinct from $e_1$ and $e_2$,
    and the edges of $\Gamma$ distinct from $e_3$ and~$e_4$;
\item
    between all faces of $\Delta$ and all faces of~$\Gamma$.
\end{enumerate}

Informally and imprecisely speaking, the diamond move consists in
cutting the map $\Delta$ along the path $e_1e_2^{-1}$ and then gluing
the sides of the obtained diamond-shaped hole in a different way
than how it was before.

Consider now an arbitrary map $\Delta$ and
two distinct oriented edges $e_1$ and $e_2$ in $\Delta$
with a common terminal vertex, none of which is a loop.
The \emph{diamond move\/} along $e_1$ and $e_2$ consists
of, first, closing the connected component of $\Delta$
that contains $e_1$ and $e_2$;
second, applying the diamond move along $e_1$ and $e_2$ to
the closure; and finally, removing the faces
that were added when closing the component.

Suppose that neither $e_1$ nor $e_2$ is a loop.
If the initial vertices of $e_1$ and $e_2$ are distinct,
then the diamond move is called \emph{proper},
otherwise it is called \emph{improper}.
If the initial vertices of $e_1$ and $e_2$ coincide,
and the cyclic path $e_1e_2^{-1}$ does not switch orientation
in $\Delta$ (informally speaking, this means that some neighbourhood
of $e_1e_2^{-1}$ is orientable),
then the (improper) diamond move is called \emph{disconnecting};
otherwise, if the cyclic path $e_1e_2^{-1}$ does switch orientation,
the (improper) diamond move is called \emph{untwisting}.

Suppose that $e_1$ is a loop, while $e_2$ is not.
In this case the dimond move is called \emph{proper\/}
for the reason that will be clear from
Lemma \ref{lemma:diamond_move}.

If both $e_1$ and $e_2$ are loops, then the task to
suitably classify such a diamond move as either \emph{proper},
or \emph{improper disconnecting}, or \emph{improper untwisting\/}
is left to the reader.
(This case is admittedly more difficult, but analogous to the
previous two.)


\begin{lemma}
\label{lemma:diamond_move}
Consider an arbitrary map\/ $\Delta$ and a map\/ $\Gamma$
obtained from\/ $\Delta$ by a diamond move\textup.
Then\/
\begin{enumerate}
\item
    $\norm{\Gamma(1)}=\norm{\Delta(1)}$ and\/
    $\norm{\Gamma(2)}=\norm{\Delta(2)}$\textup;
\item
    if the diamond move is proper\textup, then\/
    $\norm{\Gamma(0)}=\norm{\Delta(0)}$\textup,
    $\chi_\Gamma=\chi_\Delta$\textup,
    and\/ $\bar\Gamma$ is geometrically equivalent to\/ $\bar\Delta$
    if\/ $\bar\Gamma$ and\/ $\bar\Delta$ are closures of\/
    $\Gamma$ and\/~$\Delta$\textup;
\item
    if the diamond move is untwisting\textup, then\/
    $\norm{\Gamma(0)}=\norm{\Delta(0)}+1$\textup,
    and\/ $\Gamma$ has the same number of connected components
    as\/~$\Delta$\textup;
\item
    if the diamond move is disconnecting\textup, then\/
    $\norm{\Gamma(0)}=\norm{\Delta(0)}+2$\textup,
    and\/ $\Gamma$ has either the same number of connected components
    as\/ $\Delta$\textup, or\/ $1$ more\textup;
\item
    $\chi_\Delta\le\chi_\Gamma\le\chi_\Delta+2$\textup.
\end{enumerate}
\end{lemma}


Proof of this lemma is left to the reader.

Properties of diamond moves can be described even in greater detail,
but this lemma is sufficient for many applications.






\subsection{Diagrams}
\label{subsection:diagrams.ccvkd}

\begin{definition}
If \GP{\mathfrak A}{\mathcal R}
is a group presentation, then a \emph{van Kampen diagram\/}
(or simply a \emph{diagram\/})
over \GP{\mathfrak A}{\mathcal R}
is a map together with a labelling of its oriented edges
such that every two mutually inverse oriented edges are
labelled with mutually inverse group letters from $\mathfrak A^{\pm1}$,
and each group word that ``reads'' on the contour of some face
belongs to~$\mathcal R^{\pm1}$.
\end{definition}

In every van Kampen diagram,
let the \emph{label\/} of an oriented edge $e$ be denoted by $\pathlbl(e)$,
and the label of a path $p$ be denoted by $\pathlbl(p)$.

\begin{definition}
Two faces of a diagram or of two distinct diagrams are said to be
\emph{congruent\/} if
their contour labels are either cyclic shifts of each other,
or cyclic shifts of the inverses of each other.
\end{definition}

\begin{definition}
Two diagrams are called \emph{essentially isomorphic\/} if
there exists a label-preserving isomorphism between their underlying
complexes which preserves the contours of the diagrams up to re-indexing
and/or cyclically shifting and/or inverting.
\end{definition}

\begin{definition}
A pair of distinct faces $\{\Pi_1,\Pi_2\}$ in a diagram
$\Delta$ is called \emph{immediately cancellable\/} if
there are paths $p_1$ and $p_2$ in $\Delta$ such that
\begin{enumerate}
\item
    $\cntrcclof{\Pi_1}\in\{\langle p_1\rangle,\langle p_1^{-1}\rangle\}$
    and
    $\cntrcclof{\Pi_2}\in\{\langle p_2\rangle,\langle p_2^{-1}\rangle\}$,
\item
    $p_1$ and $p_2$ have a common nontrivial initial subpath, and
\item
    $\pathlbl(p_1)=\pathlbl(p_2)$.
\end{enumerate}
A diagram $\Delta$ is called \emph{weakly reduced\/} if
it does not have immediately cancellable pairs of faces.
\end{definition}

Weakly reduced diagrams are exactly the diagrams
\emph{reduced\/} in the sense of \cite{LyndonSchupp:2001:cgt}.

\begin{definition}
A \emph{diamond move\/} in a diagram $\Delta$
(a \emph{diagrammatic diamond move\/})
is a diamond move in the underlying map of $\Delta$
along two oriented edges with identical labels,
followed by the natural labelling of the obtained map 
so as to obtain a diagram.%
\footnote{Diamond moves in diagrams correspond to \emph{bridge moves\/}
in \emph{pictures}, see \cite{Rourke:1979:ptg,Huebschmann:1981:a2cupW}.}
\end{definition}

\begin{definition}
A pair of distinct faces $\{\Pi_1,\Pi_2\}$ in a diagram
$\Delta$ is called \emph{cancellable\/} if
there exists a sequence of diamond moves that separates
these two faces into an elementary spherical subdiagram
(\ie, leads to a diagram in which the faces corresponding to
$\Pi_1$ and $\Pi_2$ form an elementary spherical connected component).
A diagram $\Delta$ is called \emph{reduced\/} if
it does not have cancellable pairs of faces.
\end{definition}


\begin{lemma}
\label{lemma:cancellable_pairs}
Immediately cancellable pairs are cancellable\textup.
Reduced diagrams are weakly reduced\textup.
\end{lemma}


Proof of this lemma is left to the reader.

It should be noted that there is a substantial distinction between
the diagrams (and their transformations) defined in this paper
and such classical objects as \emph{pictures\/}
(called \emph{standard diagrams\/} in \cite{Rourke:1979:ptg})
and \emph{$0$-refined diagrams\/}
(in \cite{Olshanskii:1989:gosg-rus,Olshanskii:1991:gdrg-eng}).
If $D^2$ is a $2$-dimensional disc, $K(\mathfrak A;\mathcal R)$
is the geometric realisation of a group presentation
\GP{\mathfrak A}{\mathcal R}, and
$K^1(\mathfrak A;\mathcal R)=K(\mathfrak A;\varnothing)$ is its
$1$-skeleton, which is a wedge of circles, then
both pictures and $0$-refined diagrams
over \GP{\mathfrak A}{\mathcal R}
can be used to represent arbitrary \emph{transverse\/}
(in the sense of \cite{BuoncristianoRourkeSanderson:1976:gaht})
continuous maps
$(D^2,\partial D^2)\to
(K(\mathfrak A;\mathcal R),K^1(\mathfrak A;\mathcal R))$
up to isotopy of the domain $D^2$.
Transverse maps in turn represent arbitrary continuous maps up
to homotopy.
Certain combinatorially defined transformations of pictures,
as well as of $0$-refined diagrams, represent homotopies between
transverse continuous maps.
The diagrams defined in this paper without introducing $0$-cells
are not suitable for representing arbitrary homotopy classes of maps
$(D^2,\partial D^2)\to
(K(\mathfrak A;\mathcal R),K^1(\mathfrak A;\mathcal R))$.
Nevertheless, they are an appropriate tool
for studying relations in groups and formulating useful results
(see Lemma~\ref{lemma:diagrams}).

If \GP{\mathfrak A}{\mathcal R} is a group presentation,
$G$ is the group presented by \GP{\mathfrak A}{\mathcal R},
and $w$ is a group word over $\mathfrak A$, then let $[w]_G$,
or $[w]_{\mathcal R}$, or simply $[w]$, denote the element of $G$
represented by~$w$.

The results of the following lemma are assumed to be well-known.


\begin{lemma}
\label{lemma:diagrams}
Let\/ $G$ be the group presented by\/ \GP{\mathfrak A}{\mathcal R}\textup.
Let\/ $w$\textup, $w_1$\textup, and\/ $w_2$ be arbitrary group words
over\/ $\mathfrak A$\textup, and\/ $n$ be a natural number\textup.
Then
\begin{enumerate}
\item
\label{item:lemma.diagrams.1}
    if there exists a disc diagram\/ $\Delta$ over\/
    \GP{\mathfrak A}{\mathcal R} such that\/
    $\pathlbl(\cntr{}\Delta)=w$\textup, then\/ $[w]=1$\textup;
\item
\label{item:lemma.diagrams.2}
    if\/ $[w]=1$\textup, then there exists a reduced
    disc diagram\/ $\Delta$ over\/ \GP{\mathfrak A}{\mathcal R}
    such that\/ $\pathlbl(\cntr{}\Delta)=w$\textup;
\item
\label{item:lemma.diagrams.3}
    if there exists a contour-oriented annular diagram\/ $\Delta$
    over\/ \GP{\mathfrak A}{\mathcal R} such that\/
    $\pathlbl(\cntr_1\Delta)=w_1$ and\/
    $\pathlbl(\cntr_2\Delta)^{-1}=w_2$\textup, then\/
    $[w_1]$ and\/ $[w_2]$ are conjugate in\/ $G$\textup;
\item
\label{item:lemma.diagrams.4}
    if\/ $[w_1]$ and\/ $[w_2]$ are conjugate in\/ $G$\textup, then
    either\/ $[w_1]=[w_2]=1$\textup, or there exists a contour-oriented
    reduced annular diagram\/ $\Delta$
    over\/ \GP{\mathfrak A}{\mathcal R} such that\/
    $\pathlbl(\cntr_1\Delta)=w_1$ and\/
    $\pathlbl(\cntr_2\Delta)^{-1}=w_2$\textup;
\item
\label{item:lemma.diagrams.5}
    if there exists a one-contour diagram\/ $\Delta$
    over\/ \GP{\mathfrak A}{\mathcal R}
    such that\/ $\pathlbl(\cntr{}\Delta)=w$ and
    the underlying complex of a closure of\/ $\Delta$
    is a combinatorial sphere with\/ $n$ handles\textup, then\/
    $[w]\in[G,G]$ and\/ $\cl_G([w])\le n$\textup;
\item
\label{item:lemma.diagrams.6}
    if\/ $\cl_G([w])=n$\textup, then there exists a one-contour
    reduced diagram\/ $\Delta$ over\/ \GP{\mathfrak A}{\mathcal R}
    such that the underlying complex of a closure of\/ $\Delta$
    is a combinatorial sphere with\/ $n$ handles\textup,
    and\/ $\pathlbl(\cntr{}\Delta)=w$\textup.
\end{enumerate}
\end{lemma}


\begin{proof}[Outline of a proof]
Parts \thetag{\ref{item:lemma.diagrams.1}},
\thetag{\ref{item:lemma.diagrams.2}},
\thetag{\ref{item:lemma.diagrams.3}}, and
\thetag{\ref{item:lemma.diagrams.4}} follow,
for example, from Theorem V.1.1 and Lemmas V.1.2, V.5.1, and V.5.2 of
\cite{LyndonSchupp:2001:cgt}.
See also Lemmas 11.1 (van Kampen Lemma) and 11.2 in
\cite{Olshanskii:1989:gosg-rus,Olshanskii:1991:gdrg-eng}
(all results there are formulated in terms of $0$-refined diagrams).

Here is an outline of a proof of parts \thetag{\ref{item:lemma.diagrams.5}}
and \thetag{\ref{item:lemma.diagrams.6}}.

Suppose that $\Delta$ is a one-contour diagram
over \GP{\mathfrak A}{\mathcal R}
such that the underlying complex of a closure of $\Delta$
is a combinatorial sphere with $n$ handles,
and $\pathlbl(\cntr{}\Delta)=w$.
Let $o$ be the initial vertex of $\cntr{}\Delta$.
Consider the (combinatorial) fundamental group $\pi_1(\Delta,o)$
of $\Delta$ with base-point $o$.
It can be shown from the definition of a combinatorial handled sphere
(which is easy to formulate) that in $\pi_1(\Delta,o)$,
the homotopy class of $\cntr{}\Delta$ is the product of $n$ commutators.
Therefore $\cl_G([w])\le n$,
since there is a homomorphism $\pi_1(\Delta,o)\to G$ which maps
the homotopy class of $\cntr{}\Delta$ to~$[w]$.

Now suppose $\cl_G([w])=n$.
Let $x_1$, \dots, $x_n$, $y_1$, \dots, $y_n$ be group words over
$\mathfrak A$ such that
$[w]=\bigl[[x_1],[y_1]\bigr]\dots\bigl[[x_n],[y_n]\bigr]$ in $G$.
Let $\Psi$ be a disc diagram over \GP{\mathfrak A}{\mathcal R}
such that $\pathlbl(\cntr{}\Delta)=[x_1,y_1]\dots[x_n,y_n]w^{-1}$
(here part \thetag{\ref{item:lemma.diagrams.2}} of this lemma is used).
At this point $0$-refinement of $\Psi$ is needed.

Definition and explanation of $0$-refinement are given
in \cite{Olshanskii:1989:gosg-rus,Olshanskii:1991:gdrg-eng}.
In a $0$-refined diagram, faces and edges are usually divided into $2$
classes: $0$-edges and $0$-faces, and all the other,
``regular,'' edges and faces.
Here the terminology shall be slightly different.
The class of $0$-faces shall be subdivided into a class of
\emph{$0$-faces\/} and a class of \emph{$1$-faces\/};
``regular'' edges shall be called \emph{$1$-edges}, and ``regular'' faces
shall be called \emph{$2$-faces}.
Thus, all the edges of a $0$-refined diagram are divided
into $0$-edges and $1$-edges, and all the faces are divided into
$0$-faces, $1$-faces, and $2$-faces.
The requirements on the labelling of a $0$-refined diagram over
\GP{\mathfrak A}{\mathcal R} are the following:
\begin{enumerate}
\item
    the label of every oriented $0$-edge is $1$
    (the symbol `$1$' here is regarded as a new group letter such that
    $1^{-1}=1^{+1}=1$);
\item
    the label of every oriented $1$-edge is an element of
    $\mathfrak A^{\pm1}$, and, as usual, mutually inverse oriented
    edges are labelled with mutually inverse group letters;
\item
    the label of the contour of every $0$-face is of the form $1^k$;
\item
    the label of the contour of every $1$-face is of the form
    $1^kx1^lx^{-1}1^m$ where $x\in\mathfrak A^{\pm1}$; and
\item
    the label of the contour of every $2$-face is 
    an element of $\mathcal R^{\pm1}$.
\end{enumerate}

Let $\tilde\Psi$ be a $0$-refinement of $\Psi$ such that
$\cntr{}\tilde\Psi$ is a simple cyclic path,
and $\pathlbl(\cntr{}\tilde\Psi)=\pathlbl(\cntr{}\Psi)$.
Let $p_1$, \dots, $p_{2n}$, $q_1$, \dots, $q_{2n}$, and $t$
be the paths such that
\begin{enumerate}
\item
    $\cntr{}\tilde\Delta=p_1p_2q_1^{-1}q_2^{-1}\dots
    p_{2n-1}p_{2n}q_{2n-1}^{-1}q_{2n}^{-1}t^{-1}$,
\item
    $\pathlbl(p_{2i-1})=\pathlbl(q_{2i-1})=x_i$ and
    $\pathlbl(p_{2i})=\pathlbl(q_{2i})=y_i$ for $i=1,\dots,n$, and
\item
    $\pathlbl(t)=w$.
\end{enumerate}

Let $\Delta_0$ be the ($0$-refined) diagram obtained from $\tilde\Psi$
by ``gluing'' together each pair of paths $p_i$ and $q_i$, $i=1,\dots,n$,
and choosing $t$ (or rather the copy of $t$ in $\Delta_0$)
as the contour of $\Delta_0$.
Let $\bar\Delta_0$ be a closure of $\Delta_0$.
Then the underlying complex of $\bar\Delta_0$ is a combinatorial sphere
with $n$ handles.
Let $\Theta$ be the ``improper'' face of $\bar\Delta_0$
(which is not a face of $\Delta_0$); this face is to be regarded as
a $2$-face in the sense of $0$-refinement.

Eliminate all $0$-edges and $0$-faces of $\bar\Delta_0$ by
\emph{collapsing\/} $0$-edges.
If $e$ is a $0$-edge which is not a loop and not the only edge of some
connected component of the diagram, then
the meaning of collapsing $e$ is clear.
If $e$ is the only edge of some connected component of the diagram,
then collapsing $e$ means removing this component all together.
Consider a $0$-edge $e$ which is a loop.
If $e$ is the only edge incident to some $0$-face $\Pi_1$, and
$e$ is incident to another face $\Pi_2$ which is incident to
some edge distinct from $e$, then collapsing $e$ results in removing
$e$ and $\Pi_1$, and shortening the contour of $\Pi_2$ by $1$.
If $e$ is incident to two distinct faces $\Pi_1$ and $\Pi_2$,
both of which are also incident to some other edges,
then collapsing $e$ results in removing $e$,
possibly doubling the end-vertex of $e$
(unless $e$ switches orientation in the diagram),
and shortening the contours of $\Pi_1$
and $\Pi_2$ by~$1$.
Let $\bar\Delta_1$ be the closed map obtained from $\bar\Delta_0$
by collapsing one-by-one all $0$-edges.
Then $\bar\Delta_1$ does not have any $0$-edges or $0$-faces.
Clearly, $\bar\Delta_1$ is orientable, since so is~$\bar\Delta_0$.

Observe that if the operation of collapsing an edge increases the number
of connected components, then it increases it only by $1$,
and simultaneously increases the Euler characteristic by $2$,
and if it decreases the
Euler characteristic, then it decreases it at most by $2$ and
simultaneously decreases the number of connected components
(recall Lemma~\ref{lemma:positive_Euler_characteristic}).
Therefore, if $k$ is the number of connected components of $\bar\Delta_1$,
then $\chi_{\bar\Delta_1}\ge\chi_{\bar\Delta_0}+2(k-1)$.
Let $\bar\Delta_2$ be the connected component of $\bar\Delta_1$
that contains the face $\Theta$.
By Lemma~\ref{lemma:positive_Euler_characteristic},
$\chi_{\bar\Delta_2}\ge\chi_{\bar\Delta_0}$.

The number of connected components and the Euler
characteristic of any map that can be obtained from a given closed map
by diamond moves are both bounded from above.
Indeed, the number of connected components is bounded by the number of
faces, and hence, by Lemma~\ref{lemma:positive_Euler_characteristic},
the Euler characteristic is bounded by $2$ times the number of faces.
Let $\bar\Delta_3$ be a map of maximal Euler characteristic that
can be obtained from $\bar\Delta_2$ by diamond moves.
Since diamond moves do not decrease the Euler characteristic,
and improper diamond moves increase it
(see Lemma~\ref{lemma:diamond_move}),
no improper diamond move is applicable to $\bar\Delta_3$,
nor to any diagram obtained from $\bar\Delta_3$
by any sequence of diamond moves.

Let $\bar\Delta$ be the connected component of $\bar\Delta_3$
that contains the face $\Theta$.
Then $\bar\Delta$ does not contain any $1$-faces
(otherwise an improper diamond move would be applicable
to $\bar\Delta$).
Every diamond move that increases the number of connected components,
increases it by $1$ and simultaneously
increases the Euler characteristic by $2$.
Therefore $\chi_{\bar\Delta}\ge\chi_{\bar\Delta_2}$.
Since $\bar\Delta$ is oriented, its underlying complex is a combinatorial
sphere with at most $n$ handles,
but the number of handles cannot be less than $n$,
as follows from part~\thetag{\ref{item:lemma.diagrams.5}}.

Let $\Delta$ be the subdiagram of $\bar\Delta$ obtained
by removing $\Theta$.
The diagram $\Delta$ is reduced, because otherwise some
improper diamond move would be applicable to some
diagram obtained from $\bar\Delta$ by proper diamond moves.
The diagram $\Delta$ is such as desired.
\end{proof}



\section{Estimating Lemmas}
\label{section:estimating_lemmas}

Lemmas of this and the subsequent sections are rather technical.
It is advisable that the reader first takes a look at the proofs of
Propositions \ref{proposition:Gn_lcw},
\ref{proposition:G_icw}, and \ref{proposition:Gn_G_swcp} in
Section~\ref{section:proof_theorems}.

If $X$ is a set, then $\norm{X}$ shall denote the cardinality of $X$.
Assume the usual definitions and notation concerning binary relations
(subsets of Cartesian products).
In particular, if $R$ is a relation and $X$ is a set, then 
$$
R(X)=\{\,y\mid(\exists x\in X)(x\mathrel Ry)\,\}.
$$


\begin{named_lemma:p.hall}[Philip~Hall, 1935]
\label{lemma:Halls_lemma}
Let\/ $A$ and\/ $B$ be two finite sets\textup,
and\/ $R$ be a relation from\/ $A$ to\/ $B$ \textup(\ie\textup,
$R\subset A\times B$\textup{).}
Then the following are equivalent\/\textup:
\begin{itemize}
\item[\textup{(I)}]
    There exists an injection\/ $h\!:A\to B$ such that
    for each\/ $x\in A$\textup, $x\mathrel Rh(x)$\textup.
\item[\textup{(II)}]
    For each subset\/ $X$ of\/ $A$\textup,
    $\norm{R(X)}\ge\norm{X}$\textup.
\end{itemize}
\end{named_lemma:p.hall}


\begin{corollary.named_lemma:p.hall}
\label{corollary.lemma:1.Halls_lemma}
Let\/ $A$ and\/ $B$ be two finite sets\textup,
and\/ $R$ be a relation from\/ $A$ to\/ $B$\textup.
Let\/ $w$ be a function from\/ $B$ to\/ $\mathbb N\cup\{0\}$\textup.
Then the following are equivalent\/\textup:
\begin{itemize}
\item[\textup{(I)}]
    There exists a function\/ $h\!:A\to B$ such that\/\textup:
    \begin{enumerate}
    \item
        for each\/ $x\in A,$ $x\mathrel Rh(x)$\textup, and
    \item
        for each\/ $y\in B$\textup, the full pre-image of\/ $y$
        under\/ $h$ consists of at most\/ $w(y)$ elements\textup.
    \end{enumerate}
\item[\textup{(II)}]
    For each subset\/ $X$ of\/ $A$\textup,
    $\displaystyle\sum_{y\in R(X)}w(y)\ge\norm{X}$\textup.
\item[\textup{(III)}]
    For each subset\/ $Y$ of\/ $B$\textup,
    $\displaystyle\sum_{y\in Y}w(y)
    \ge\norm{\{\,x\mid R(\{x\})\subset Y\,\}}$\textup.
\end{itemize}
\end{corollary.named_lemma:p.hall}


Proofs of the lemma and the corollary may be found, for example,
in \cite{Hall:1935:ors,Muranov:2005:dsmcbgbsg}.
(The equivalence of items
(II) and (III)
of the corollary is not proved
in those papers, but is easy to verify.)

\begin{definition}
A c-path is called \emph{regular\/} if
its image in the $1$-skeleton of the complex is reduced and nontrivial.
A c-pseudo-arc is \emph{regular\/} if
the associated oriented c-pseudo-arcs are such.
\end{definition}

\begin{definition}
An \emph{$S_1$-map\/} is a map together with a system of
\emph{selected\/} c-paths of its faces satisfying
the following conditions:
\begin{enumerate}
\item
    all selected c-paths are regular (in particular, they are
    oriented c-pseudo-arcs),
\item
    the inverse path of every selected c-path is selected, and
\item
    every nontrivial subpath of every selected c-path is selected.
\end{enumerate}
A c-pseudo-arc of a face in an $S_1$-map is \emph{selected\/} if
the associated oriented c-pseudo-arcs are selected.
An arc in an $S_1$-map is \emph{selected\/} if
this arc is internal and both c-arcs that map to it
(by attaching morphisms) are selected.
\end{definition}

This definition of $S_1$-maps is similar to
the definition of S-maps in \cite{Muranov:2005:dsmcbgbsg},
but it is adapted to the more general definition of maps
(one of the generalisations is that
maps now are allowed to be non-orientable).
In \cite{Muranov:2005:dsmcbgbsg}, ``S'' stood for ``selection,''
and here it stands for ``structure,''
because an $S_1$-map is a map with additional structure.
So are $S_2$-maps and $S$-maps,
which shall be defined and used below.

\begin{definition}
A set $X$ of c-pseudo-arcs \emph{encloses\/} a simple disc
submap $\Phi$ if
\begin{enumerate}
\item
    elements of $X$ are c-pseudo-arcs of faces which do not belong
    to $\Phi$ (are ``outside'' of $\Phi$),
\item
    for every c-pseudo-arc from $X$, one of the associated oriented
    c-pseudo-arcs maps to a subpath of $\cntrcclof{{}\Phi}$, and
\item
    $\cntr{}\Phi$ can be decomposed into a product of paths
    each of which is the image of a subpath of an oriented
    c-pseudo-arc associated with an element of~$X$.
\end{enumerate}
\end{definition}

\begin{definition}
Let $\Delta$ be an $S_1$-map, $\Phi$ its simple disc submap, 
and $n\in\mathbb N$.
The $S_1$-map $\Delta$ is said to satisfy
the \emph{condition\/ $\mathsf Z(n)$ relative to\/ $\Phi$} if
every set of selected c-pseudo-arcs enclosing $\Phi$ in $\Delta$
has at least $n+1$
element.
\end{definition}


\begin{lemma}
\label{lemma:simple_submap__disc_map__partitioned_contour}
Let\/ $\Psi$ be a non-degenerate disc map\textup.
Suppose\/ $\cntr{}\Psi=p_1\dots p_n$\textup, $n\in\mathbb N$\textup,
where\/ $p_1$\textup{, \dots,} $p_n$ are reduced paths\textup.
Then there exist a maximal simple disc submap\/ $\Phi$ of\/ $\Psi$
and simple paths\/ $q_1$\textup{, \dots,} $q_m$\textup,
$1\le m\le n$\textup, such that\/ $\cntr{}\Phi=q_1\dots q_m$
and there are\/ $i_1$\textup{, \dots,} $i_m$ such that\/
$1\le i_1<\dots<i_m\le n$ and for every $j=1,\dots,m$\textup,
$q_j$ is a subpath of\/~$p_{i_j}$\textup.
\end{lemma}


(This Lemma is similar to Proposition~3.1
in~\cite{Muranov:2005:dsmcbgbsg}.)


\begin{proof}
Without loss of generality, assume that all the paths
$p_1$, \dots, $p_n$ are nontrivial, and
that the contour of $\Psi$ is cyclically reduced.
For if it is not, then the terminal vertex of one of the paths
$p_1$, \dots, $p_n$ has degree $1$ in $\Psi$.
Remove this vertex together with the incident edge,
and ``shorten'' or remove each of the paths
from among $p_1$, \dots, $p_n$ that start or end at this vertex
(the contour of $\Psi$ is also ``shortened'').
If the lemma holds for the new disc map and the new set of
paths, it is clear that it holds for the initial ones.
Thus it can be assumed that $\cntr{}\Psi$ is cyclically reduced.

The conclusion is obvious if $\Psi$ is simple.
Assume $\Psi$ is not simple.
Then it has two maximal simple disc submaps
whose contours are subpaths of $\cntrcclof{{}\Psi}$, and
which are either disjoint, or have only one vertex in common.
Let $\Phi_1$ and $\Phi_2$ be such maximal simple disc submaps.

If $\cntr{}\Phi_2$ is a subpath of one of the paths $p_1$, \dots, $p_n$,
then take $\Phi=\Phi_2$, $m=1$, $q_1=\cntr{}\Phi$,
and see that the conclusion holds.

Suppose $\cntr{}\Phi_2$ is not a subpath of any one
of the paths $p_1$, \dots, $p_n$.
If the initial vertex of $\cntr{}\Psi$ is not
in $\Phi_1$, then let $\Phi=\Phi_1$.
If the initial vertex of $\cntr{}\Psi$ is
in $\Phi_1$, then let $\Phi$ be the map obtained from $\Phi_1$
by cyclically shifting its contour so that
$\cntr{}\Phi$ starts at the same vertex as $\cntr{}\Psi$.
The initial vertices of the $n+1$ paths
$p_1$, \dots, $p_n$, $\cntr{}\Phi_1$ divide the simple
path $\cntr{}\Phi$ into at most $n$ simple subpaths.
Denote these subpaths by $q_1$, \dots, $q_m$ so that
$\cntr{}\Phi=q_1\dots q_m$.
The submap $\Phi$ and the path $q_1$, \dots, $q_m$ are
the desired ones.
\end{proof}


The following notation is used in
Estimating Lemma~\ref{estimating_lemma:selected_arcs}
and throughout the rest of this paper:
if $\Pi$ is a face of an $S_1$-map $\Delta$, then let
$\kappa_{\Delta}(\Pi)$, or $\kappa(\Pi)$, denote the number of
maximal selected c-pseudo-arcs of $\Pi$, and
$\kappa_{\Delta}'(\Pi)$, or $\kappa'(\Pi)$, denote the number of maximal
piece-wise selected regular c-pseudo-arcs of~$\Pi$.
Note that $\kappa_{\Delta}'(\Pi)\le\kappa_{\Delta}(\Pi)$.
Note also that if all c-pseudo-arcs of $\Pi$ are selected,
as well as if no c-pseudo-arc of $\Pi$ is selected, then
$\kappa_{\Delta}(\Pi)=\kappa_{\Delta}'(\Pi)=0$.

Recall that an \emph{elementary map\/} is
a spherical map with exactly $2$ faces,
whose $1$-skeleton is a combinatorial circle.
Elementary maps are ``bad'' in the sense that the conclusion of
Estimating Lemma~\ref{estimating_lemma:selected_arcs}
may fail for them
(but only if all c-pseudo-arcs are selected,
and hence $\kappa=\kappa'=0$).
They are also ``inconvenient'' in the sense that their
distinct maximal selected arcs can overlap.


\begin{estimating_lemma}[First Estimating Lemma]
\label{estimating_lemma:selected_arcs}
Let\/ $\Delta$ be a non-elementary connected\/ $S_1$-map\textup,
or an elementary\/ $S_1$-map which has a maximal selected
c-pseudo-arc\textup.
Let\/ $A$ be a set of selected internal arcs of\/ $\Delta$ such that
no two distinct elements of\/ $A$ are subarcs of the same
selected arc\textup.
Let\/ $C$ and\/ $D$ be sets of faces of\/ $\Delta$
such that\/\textup:
\begin{enumerate}
\item
    $C$ contains all faces incident to arcs from\/ $A$\textup, and
\item
    $\Delta$ satisfies the condition\/ $\mathsf Z(2)$
    relative to every simple disc submap that
    does not contain any faces from\/ $D$ and does not contain
    at least one arc from\/~$A$\textup.
\end{enumerate}
Let\/ $c_\Delta$ be the number of contours of\/ $\Delta$\textup.
Then either\/ $A$ is empty\textup, or
$$
\norm{A}\le\sum_{y\in C}\!
\bigl(3+\kappa_\Delta(y)+\kappa_\Delta'(y)\bigr)
+2\norm{D\setminus C}-c_\Delta-3\chi_\Delta.
$$
Furthermore\textup, if\/ $B$ is a subset of\/ $C$\textup,
there exist a subset\/ $E\subset A$ and a function\/
$f\!:A\setminus E\to B$ such that\/\textup:
\begin{enumerate}
\item
    either\/ $E$ is empty\textup, or
    \begin{align*}
    \norm{E}&\le\sum_{y\in (C\setminus B)\setminus D}\!
    \bigl(3+\kappa_\Delta(y)+\kappa_\Delta'(y)\bigr)\\
    &\qquad+\sum_{y\in (C\setminus B)\cap D}\!
    \bigl(1+\kappa_\Delta(y)+\kappa_\Delta'(y)\bigr)\\
    &\qquad+2\norm{D}-c_\Delta-3\chi_\Delta;
    \end{align*}
\item
    for every\/ $x\in A\setminus E$\textup,
    $f(x)$ is incident to\/~$x$\textup;
\item
    for every\/ $y\in B$\textup,
    the full pre-image of\/ $y$ under\/ $f$
    consists of at most\/
    $3+\kappa_\Delta(y)+\kappa_\Delta'(y)$ elements\/\textup;
\item
    for every\/ $y\in D$\textup,
    the full pre-image of\/ $y$ under\/ $f$
    consists of at most\/
    $1+\kappa_\Delta(y)+\kappa_\Delta'(y)$ elements\textup.
\end{enumerate}
\end{estimating_lemma}


\begin{proof}
If $A$ is empty, then there is nothing to prove
(meaning the proof is easy).
Assume it is non-empty.

It suffices to prove this lemma in the case $\Delta$ is closed.
(To prove the statement in the case $\Delta$ is not closed,
apply this lemma to a closure of $\Delta$, the same sets $A$, $B$, $C$,
and the set $D$ extended by including the attached ``improper'' faces).
Hence assume without loss of generality that $\Delta$ is closed.

Let $K$ be the set of all connected components of a submap
obtained from $\Delta$ by removing
all the faces that are in $C$ and all the arcs that are in $A$.

For every element $\Psi$ of $K$, let $d(\Psi)$ denote
the number of arcs in $A$ that have exactly one end-vertex in $\Psi$,
plus twice the number of arcs in $A$ that have both end-vertices 
in $\Psi$.
(Thus $d$ is analogous to vertex degree.)

Clearly,
$$
\sum_{\Psi\in K}\chi_\Psi-\norm{A}+\norm{C}=\chi_\Delta,
$$
and
$$
\sum_{\Psi\in K}d(\Psi)=2\norm{A}.
$$
Using these two equations, one has
\begin{align*}
\norm{A}&=3\norm{C}+3\sum_{\Psi\in K}\chi_\Psi-2\norm{A}-3\chi_\Delta\\
&=3\norm{C}+\sum_{\Psi\in K}\bigl(3\chi_\Psi-d(\Psi)\bigr)-3\chi_\Delta.
\end{align*}

By Lemmas \ref{lemma:positive_Euler_characteristic} and
\ref{lemma:maps__positive_Euler_characteristic},
the Euler characteristic of each element of $K$ is at most $1$,
and if the Euler characteristic of $\Psi\in K$ is $1$,
then $\Psi$ is a disc map.
Let
$$
K_i'=\{\,\Psi\in K\,|\,d(\Psi)=i\ \text{and}\ \chi_\Psi=1\,\}
\quad\text{for}\quad i=0,1,2,\dots.
$$
Each element of each $K_i'$ is a disc map.
Observe that $K_0'=\varnothing$.
Therefore,
$$
\norm{A}\le3\norm{C}+2\norm{K_1'}+\norm{K_2'}-3\chi_\Delta.
$$
To complete the proof, essentially it is only left to demonstrate that
$$
2\norm{K_1'}+\norm{K_2'}
\le\sum_{\Pi\in C}\!
\bigl(\kappa_\Delta(\Pi)+\kappa_\Delta'(\Pi)\bigr)+2\norm{D\setminus C},
$$
and then to apply the corollary of Hall's Lemma.

For $i=1,2$, let $K_i''$ be the set of those elements of $K_i'$
whose face sets are disjoint with $D$
(\ie, such $\Psi\in K_i'$ that $\Psi(2)\cap D=\varnothing$).
Clearly,
$$
\norm{K_1'\setminus K_1''}+\norm{K_2'\setminus K_2''}\le\norm{D\setminus C}.
$$
Now it is to be shown that
$$
2\norm{K_1''}+\norm{K_2''}
\le\sum_{\Pi\in C}\!\bigl(\kappa_\Delta(\Pi)+\kappa_\Delta'(\Pi)\bigr).
$$

Let $L$ be the set of all positive (for definiteness)
maximal selected c-paths of all face from $C$.
Let $L'$ be the set of all elements of $L$
that are terminal subpaths of maximal piece-wise selected
regular c-path.
Clearly,
$$
\norm{L}=\sum_{\Pi\in C}\!\kappa_\Delta(\Pi)
\quad\text{and}\quad
\norm{L'}=\sum_{\Pi\in C}\!\kappa_\Delta'(\Pi).
$$
Note that the image in $\Delta^1$ of every element of $L$
has a common vertex with at least one element of $K$
(because of the maximality of elements of~$L$).

Let a function $h\!:L\to K$ be defined as follows:
$h(x)$ is the element of $K$ such that
the image of some terminal subpath of $x$ has a common
vertex with $h(x)$ and no common vertices with any other
element of~$K$.

Assign weights to all elements of $L$ so that
the weight of every element of $L'$ is $2$, and the weight of every
element of $L\setminus L'$ is~$1$.
Let the weight of every element of $K$ be the sum
of the weights of all elements of its full pre-image under~$h$.

Consider an arbitrary $\Psi\in K_1''$.
Let $v$ be the oriented arc that represents an element of $A$
and whose terminal vertex is in $\Psi$.
Let $\Pi$ be the face incident to $v$.
Let $v'$, $q'$, and $u'$ be c-paths of $\Pi$ such that
$v'q'u'$ is a c-path of $\Pi$ as well,
the images of $v'$ and $u'$ are $v$ and $v^{-1}$, respectively,
and the image of $q'$ represents $\cntrcclof{\Psi}$
(see Fig.~\ref{figure:2}).
Both $v'$ and $u'$ are selected c-paths.
Let $q$ be the image of~$q'$.

\begin{figure}\centering
\includegraphics{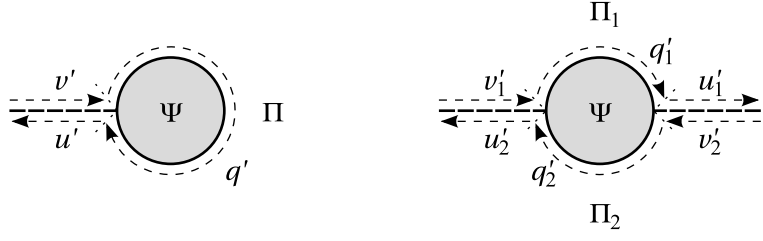}
\caption{Illustrations of $\Psi\in K_1''$ (left) and 
$\Psi\in K_2''$ (right).}
\label{figure:2}
\end{figure}

Suppose $\Psi$ is degenerate.
Then $vqv^{-1}$ is not reduced,
which means that $v'q'u'$ is not regular.
Consider the maximal positive piece-wise selected regular c-path $x$
containing either $v'$ or $u^{\prime-1}$ as a subpath.
Clearly, the image of the terminal vertex of $x$ is in $\Psi$,
and therefore $x$ has a (terminal) subpath which is an element of $L'$.
The image of this element of $L'$ under $h$ is $\Psi$.
Therefore, the weight of $\Psi$ is at least~$2$.

Suppose $\Psi$ is non-degenerate.
By Lemma~\ref{lemma:simple_submap__disc_map__partitioned_contour},
there is a simple disc submap $\Phi$ of $\Psi$
whose contour is a subpath of $q$.
Since $\Delta$ satisfies $\mathsf Z(2)$ relative to $\Phi$,
the (nontrivial) c-path $q'$ cannot be selected
and cannot be the product of two selected c-paths.
Therefore, the full pre-image of $\Psi$ under $h$ must contain either
an element of $L'$, or at least $2$ distinct elements of $L$.
(If the pre-image does not contain any element of $L'$, then
$v'q'u'$ is piece-wise selected;
if additionally the pre-image consisted of a single element,
then $q'$ would be selected or would be
the product of two selected c-paths, which is impossible.)
In either case the weight of $\Psi$ is at least~$2$.

Consider an arbitrary $\Psi\in K_2''$.
Let $v_1$ and $v_2$ be the two oriented arcs
that represent elements of $A$
and whose terminal vertices are in $\Psi$.
Let $v_1'$, $q_1'$, $u_1'$, $v_2'$, $q_2'$, and $u_2'$
be c-paths, and $q_1$ and $q_2$ be paths such that:
\begin{enumerate}
\item
    $v_1'q_1'u_1'$ and $v_2'q_2'u_2'$ are c-paths
    (\ie, the products are defined),
\item
    the images of $v_1'$, $q_1'$, $u_1'$, $v_2'$, $q_2'$, $u_2'$
    are $v_1$, $q_1$, $v_2^{-1}$, $v_2$, $q_2$, $u_1^{-1}$,
    respectively, and
\item
    $\langle q_1q_2\rangle=\cntrcclof{{}\Psi}$ (see Fig.~\ref{figure:2}).
\end{enumerate}
The c-paths $v_1'$, $u_1'$, $v_2'$, and $u_2'$ are selected.
Let $\Pi_1$ and $\Pi_2$ be the faces to which the c-paths
$v_1'q_1'u_1'$ and $v_2'q_2'u_2'$ respectively belong.

Suppose $\Psi$ is degenerate.
If $v_1q_1v_2^{-1}$ or $v_2q_2v_1^{-1}$ is not reduced,
which means that $v_1'q_1'u_1'$ or $v_2'q_2'u_2'$ is not regular,
then the pre-image of $\Psi$ under $h$ contains at least one
element of $L'$.
Suppose now that both $v_1q_1v_2^{-1}$ and $v_2q_2v_1^{-1}$
are reduced.
Then they are inverse to each other.
Therefore, they are oriented arcs in $\Delta$, unless $v_1=v_2^{-1}$.
If $v_1=v_2^{-1}$, then $\Delta$ is elementary,
$K=\{\Psi\}$, $L\ne\varnothing$, and hence the full pre-image of
$\Psi$ under $h$ is non-empty.
Hence, suppose that $v_1q_1v_2^{-1}$ and $v_2q_2v_1^{-1}$ are
mutually inverse oriented arcs of $\Delta$.
The associated non-oriented arc cannot be selected because otherwise
it would be a selected arc containing two distinct elements of $A$
as subarcs.
Therefore, at least one of the c-paths
$v_1'q_1'u_1'$ or $v_2'q_2'u_2'$ is not selected, which implies that
the maximal selected c-path containing one of the c-paths
$v_1'$, $u_1^{\prime-1}$, $v_2'$, or $u_2^{\prime-1}$
as a subpath is mapped by $h$ to $\Psi$.
Thus, the weight of every degenerate element of $K_2''$
is at least~$1$.

Suppose $\Psi$ is non-degenerate.
By Lemma~\ref{lemma:simple_submap__disc_map__partitioned_contour},
there is a simple disc submap $\Phi$ of $\Psi$
whose contour either is a subpath of one of the paths $q_1$ or $q_2$,
or is the product of a subpath of $q_1$ and a subpath of $q_2$.
Since $\Delta$ satisfies $\mathsf Z(2)$ relative to $\Phi$,
at least one of the c-paths $q_1'$ or $q_2'$ is nontrivial
but not selected.
Therefore, the full pre-image of $\Psi$ under $h$ is not empty,
and the weight of $\Psi$ is at least~$1$.

On one hand, the sum of the weights of all elements of
$K_1''\sqcup K_2''$ is at least $2\norm{K_1''}+\norm{K_2''}$.
On the other hand, it equals the sum of the weights of all elements
of $L$, which is $\norm{L}+\norm{L'}$.
Therefore,
\begin{align*}
2\norm{K_1'}+\norm{K_2'}
&=2\norm{K_1''}+\norm{K_2''}
+2\norm{K_1'\setminus K_1''}+\norm{K_2'\setminus K_2''}\\
&\le\norm{L}+\norm{L'}+2\norm{D\setminus C}\\
&\le\sum_{\Pi\in C}\!
\bigl(\kappa_\Delta(\Pi)+\kappa_\Delta'(\Pi)\bigr)+2\norm{D\setminus C}.
\end{align*}
This gives
$$
\norm{A}\le\sum_{\Pi\in C}\!
\bigl(3+\kappa_\Delta(\Pi)+\kappa_\Delta'(\Pi)\bigr)
+2\norm{D\setminus C}-3\chi_\Delta.
$$

Now the first part of the statement of the lemma is proved.
Apply it to all subset of $A$.
More precisely, take an arbitrary subset $X$ of $A$,
take the subset $Y$ of $C$ consisting of all the faces incident to
elements of $X$,
and apply the proved part of the lemma to conclude that
$$
\norm{X}\le\sum_{y\in Y}\!
\bigl(3+\kappa_\Delta(y)+\kappa_\Delta'(y)\bigr)
+2\norm{D\setminus Y}-3\chi_\Delta.
$$
Let $w$ be the function on $C$ defined as follows:
$$
w(\Pi)=
\begin{cases}
3+\kappa_\Delta(\Pi)+\kappa_\Delta'(\Pi)
\quad\text{for}\quad\Pi\in C\setminus D,\\
1+\kappa_\Delta(\Pi)+\kappa_\Delta'(\Pi)
\quad\text{for}\quad\Pi\in C\cap D.
\end{cases}
$$
In terms of $w$, have
\begin{align*}
\norm{X}&\le\sum_{y\in Y}\!w(y)+2\norm{D}-3\chi_\Delta\\
&\le\sum_{y\in (Y\cap B)}\!w(y)
+\sum_{y\in (C\setminus B)}\!w(y)+2\norm{D}-3\chi_\Delta.
\end{align*}

Now apply the corollary of Hall's Lemma to verify the remaining part
of the lemma.
Let $\omega$ be an arbitrary element not in $B$.
Define a binary relation $R\subset A\times(B\cup\{\omega\})$ as follows:
$x\mathrel Ry$ if and only if $x\in A$ and either $y=\omega$, or $y\in B$
and $y$ is incident to $x$.
Use the corollary of Hall's Lemma and the last inequality
to conclude that there is a function $h\!:A\to B\cup\{\omega\}$
such that:
\begin{enumerate}
\item
    $\norm{h^{-1}(\omega)}\le\max\{0,
    \sum_{y\in (C\setminus B)}\!w(y)+2\norm{D}-3\chi_\Delta\}$;
\item
    for every $x\in A$, either $h(x)$ is incident to $x$, or
    $h(x)=\omega$;
\item
    for every $y\in B$,
    $\norm{h^{-1}(y)}\le w(y)$.
\end{enumerate}
To complete the proof of the second part,
take $E=h^{-1}(\omega)$ and $f=h|_{A\setminus E}$.
\end{proof}


\begin{definition}
A \emph{graded\/} map is a map $\Delta$ together with a function
$\rank_\Delta\!:\Delta(2)\to J$ where $J$ is an arbitrary set or
algebraic structure.
The \emph{rank\/} of a face $\Pi$ of $\Delta$ is $\rank(\Pi)$.
Two faces are called \emph{rank-equivalent\/} if
their ranks are equal.
\end{definition}

\begin{definition}
An \emph{$S_2$-map\/} is a graded map together with
a system of \emph{exceptional\/} arcs such that:
\begin{enumerate}
\item
    distinct exceptional arcs do not overlap,
\item
    every exceptional arc is incident to a face, and
\item
    faces incident to the same exceptional arc are of the same rank.
\end{enumerate}
\end{definition}

Assign a \emph{rank\/} to every exceptional arc of an $S_2$-map
according to the rank of the incident faces.
Exceptional arcs of the same rank shall be called \emph{rank-equivalent}.

Consider an arbitrary \emph{connected\/} $S_2$-map $\Delta$.
For every $j$, 
let $\Gamma_j$ denote the subcomplex of $\Delta$ obtained by
removing all the faces of rank $j$ and all the internal exceptional arcs
of rank~$j$.

\begin{definition}
The $S_2$-map $\Delta$ is said to satisfy
the \emph{condition\/ $\mathsf Y$} if
for every $j$ such that $\Delta$ has an internal exceptional
arc of rank $j$,
the number of connected component of $\Gamma_j$ that either have
Euler characteristic $1$ or contain a rank-$j$ (external) exceptional
arc of $\Delta$
is less than or equal to the number of faces of $\Delta$ of rank~$j$.
\end{definition}

Note that every connected component of $\Gamma_j$ which contains an
(external) exceptional arc of $\Delta$ of rank $j$
is the underlying subcomplex of a map with at least $2$ contours,
and hence has non-positive Euler characteristic.



\begin{estimating_lemma}[Second Estimating Lemma]
\label{estimating_lemma:exceptional_arcs}
Let\/ $\Delta$ be a connected\/ $S_2$-map satisfying
the condition\/ $\mathsf Y$\textup.
For every\/ $j$\textup, let\/ $A_j$ denote the set of all the
internal exceptional arcs of\/ $\Delta$ of rank\/ $j$\textup, and\/
$B_j$ denote the set of all the faces of\/ $\Delta$ of rank\/ $j$\textup.
For every\/ $j$\textup, let\/ $\varepsilon(j)=1$
if\/ $\Delta$ has an external exceptional arc of rank\/ $j$\textup,
and let\/ $\varepsilon(j)=0$ otherwise\textup.
Then for every\/ $j$\textup, either\/ $A_j$ is empty\textup, or
$$
\norm{A_j}\le2\norm{B_j}-\varepsilon(j)-\chi_\Delta.
$$
Furthermore\textup, there exists a set\/ $E$ such that\/\textup:
\begin{enumerate}
\item
    either\/ $E$ is empty\textup, or\/
    $\norm{E}\le-\chi_\Delta$\textup, and
\item
    for every\/ $j$\textup,
    $\norm{A_j\setminus E}\le2\norm{B_j}-\varepsilon(j)$\textup.
\end{enumerate}
\end{estimating_lemma}


\begin{proof}
For every set $J$,
let $A_J=\bigcup_{j\in J}A_j$, $B_J=\bigcup_{j\in J}B_j$,
and let $\Gamma_J$ be the subcomplex obtained from the underlying
complex of $\Delta$ by removing all the faces that are in $B_J$
and all the arcs that are in~$A_J$.

It follows from the condition $\mathsf Y$ that for every $j$
such that $A_j\ne\varnothing$,
the number of connected components of $\Gamma_{\{j\}}$
of Euler characteristic $1$ is at most $\norm{B_j}-\varepsilon(j)$.
Observe also that for every $j$, $\norm{B_j}-\varepsilon(j)\ge0$.

Let $J$ be an arbitrary set such that $A_J$ is non-empty.
It is to be shown that
$$
\chi_{\Gamma_J}\le\norm{B_J}-\sum_{j\in J}\varepsilon(j).
$$

Let $K$ be the set of all connected components of $\Gamma_J$.
By Lemmas~\ref{lemma:positive_Euler_characteristic} and
\ref{lemma:maps__positive_Euler_characteristic},
the Euler characteristic of each elements of $K$ is at most $1$,
and every element of $K$ of Euler characteristic $1$ is
the underlying complex of a disc submap of $\Delta$.
Let $K'$ be the set of all the elements of $K$
of Euler characteristic~$1$.

Define a function $f\!:K'\to J$ as follows.
Consider an arbitrary $\Psi\in K'$.
Let $A_J^{(\Psi)}$ be the set of all the elements of $A_J$
that have an end-vertex in $\Psi$.
Since $\Delta$ is connected and $A_J\ne\varnothing$
(and hence at least one arc has been removed in the process
of obtaining $\Gamma_J$),
the set $A_J^{(\Psi)}$ is non-empty.
Since $\Psi$ is the underlying complex of a disc submap of $\Delta$
(and a disc map has only $1$ contour),
all elements of $A_J^{(\Psi)}$
are of the same rank (see Lemma~\ref{lemma:submap}).
Let $f(\Psi)$ be the rank of the elements of~$A_J^{(\Psi)}$.
Then $A_{f(\Psi)}\supset A_J^{(\Psi)}\ne\varnothing$.

Because $\Delta$ satisfies the condition $\mathsf Y$, and 
every $\Psi\in K'$ is a connected component of $\Gamma_{\{f(\Psi)\}}$,
it follows that 
$$
\norm{\{\,\Psi\in K'\mid f(\Psi)=j\,\}}\le\norm{B_j}-\varepsilon(j)
\qquad\text{for every $j$}.
$$
Therefore,
$$
\chi_{\Gamma_J}=\sum_{\Psi\in K}\chi_\Psi\le\norm{K'}
\le\sum_{j\in J}\bigl(\norm{B_j}-\varepsilon(j)\bigr).
$$

Let $J$ be an arbitrary set.
Then $\chi_\Delta=\chi_{\Gamma_J}-\norm{A_J}+\norm{B_J}$,
and therefore
$$
\norm{A_J}=\norm{B_J}+\chi_{\Gamma_J}-\chi_\Delta.
$$
In the case $A_J$ is non-empty, obtain:
$$
\norm{A_J}\le\sum_{j\in J}\bigl(2\norm{B_j}-\varepsilon(j)\bigr)
-\chi_\Delta
\le\sum_{j}\bigl(2\norm{B_j}-\varepsilon(j)\bigr)-\chi_\Delta.
$$
In particular this proves the first part of the statement
(take $J$ to be the one-element set $\{j\}$).
To prove the second part, take $J$ to be the set of all $j$ such that
$\norm{A_j}>2\norm{B_j}-\varepsilon(j)$,
and observe from the last inequality that a desired set
$E\subset A_J$ exists.
\end{proof}



\section{$S$-maps}
\label{section:S-maps}

\begin{definition}
An $S$-map is a map together with structures of an $S_1$-map 
and an $S_2$-map such that every internal exceptional arc is selected,
and every external exceptional arc lies on the image of a selected
c-path.
\end{definition}

Every submap of an $S$-map has a natural structure of an $S$-map.
If $\Gamma$ is an $S$-submap of an $S$-map $\Delta$, then an arc
of $\Gamma$ is exceptional in $\Gamma$ if and only if it is exceptional
in $\Delta$ and is incident to a face of~$\Gamma$.

\begin{definition}
An $S$-map $\Delta$ is said to satisfy the
\emph{condition\/ $\mathsf D(\lambda,\mu,\nu)$} relative to
a submap $\Gamma$ if
$\lambda$, $\mu$, and $\nu$ are functions defined on $\Gamma(2)$
(and possibly elsewhere)
with values in $[0,1]$ 
such that the following three conditions hold:
\begin{description}
\item[$\mathsf D_1(\lambda)$]
    if $\Pi$ is a face of $\Gamma$, and
    $L$ is the number 
    of non-selected c-edges of $\Pi$,
    then
    $$
    L\le\lambda(\Pi)\abs{\cntr\Pi};
    $$
\item[$\mathsf D_2(\mu)$]
    if $\Pi$ is a face of $\Gamma$, 
    $u$ is a selected internal arc of $\Delta$ incident to $\Pi$,
    and $M$ is the number of the edges of $u$
    that do not lie on any exceptional arc, then
    $$
    M\le\mu(\Pi)\abs{\cntr\Pi};
    $$
\item[$\mathsf D_3(\nu)$]
    if $\Pi$ is a face of $\Delta$,
    $p$ is a simple path in $\Delta$ which is the image of a selected
    c-path of $\Pi$, and
    $N$ is the sum of the lengths of all the 
    exceptional arcs of $\Gamma$ that lie on $p$,
    then
    $$
    N\le\nu(\Theta)\abs{\cntr\Theta}
    $$
    for every face $\Theta$ of $\Gamma$ such that
    $\rank(\Theta)=\rank(\Pi)$.
\end{description}
The $S$-map $\Delta$ is said to satisfy the
\emph{condition $\mathsf D(\lambda,\mu,\nu)$} 
(absolutely) if
it satisfies it relative to itself.
\end{definition}

Let $\mathsf D_2'(\mu)$ denote the condition obtained from
$\mathsf D_2(\mu)$ by replacing ``\dots\ $M\le\mu(\Pi)\abs{\cntr\Pi}$''
with
``\dots\ $M\le\mu(\Theta)\abs{\cntr\Theta}$
for every face $\Theta$ of $\Gamma$ such that
$\rank(\Theta)\ge\rank(\Pi)$.''

Let $\mathsf D_3'(\nu)$ denote the condition obtained from
$\mathsf D_3(\nu)$ by replacing
``\dots\ such that $\rank(\Theta)=\rank(\Pi)$'' with
``\dots\ such that $\rank(\Theta)\ge\rank(\Pi)$.''

\begin{definition}
The \emph{condition\/ $\mathsf D'(\lambda,\mu,\nu)$} is the conjunction
of the conditions
$\mathsf D_1(\lambda)$, $\mathsf D_2'(\mu)$, and $\mathsf D_3'(\nu)$.
An $S$-map $\Delta$ is said to satisfy the
\emph{condition $\mathsf D'(\lambda,\mu,\nu)$} 
absolutely if
it satisfies it relative to itself.
\end{definition}

Note that
if an $S$-map $\Delta$ satisfies $\mathsf D(\lambda,\mu,\nu)$
or $\mathsf D'(\lambda,\mu,\nu)$
relative to a submap $\Gamma$, then 
$\Delta$ satisfies the same condition relative to
every submap of $\Gamma$ as well.

The condition $\mathsf D$ will be used in the proof of
Theorem~\ref{theorem:bsglcw}, and the somewhat stronger
condition $\mathsf D'$ will be used in the proof of 
Theorem~\ref{theorem:sgicw}.


\begin{inductive_lemma}[Inductive Lemma]
\label{lemma:inductive_lemma}
Let\/ $\Delta$ be an\/ $S$-map\textup,
and\/ $\Phi$ be a simple disc $S$-submap of\/ $\Delta$\textup.
Assume\/ $\Delta$ satisfies the condition\/ $\mathsf Z(2)$
relative to every proper simple disc submap of\/ $\Phi$\textup,
$\Phi$ satisfies the condition\/ $\mathsf Y$\textup,
and\/ $\Delta$ satisfies\/ $\mathsf D(\lambda,\mu,\nu)$ 
relative to\/ $\Phi$\textup.
Suppose
$$
\lambda+(3+\kappa+\kappa')\mu+2\nu\le\frac{1}{2}
$$
point-wise on\/ $\Phi$ 
\textup(\ie\textup, for every face of\/ $\Phi$\textup{).}
Then\/ $\Delta$ satisfies\/ $\mathsf Z(2)$
relative to\/~$\Phi$\textup.
\end{inductive_lemma}


\begin{proof}
Suppose $\Delta$ does not satisfy $\mathsf Z(2)$ relative to~$\Phi$.

Let $\bar\Phi$ be a (spherical) closure of $\Phi$.
Note that the $1$-skeleton of $\bar\Phi$ is a subcomplex of 
the $1$-skeleton of $\Delta$.
Let $\Theta$ be the face of $\bar\Phi$ that is not in $\Phi$
(the \emph{improper\/} face).
Endow $\bar\Phi$ with a structure of an $S_1$-map by selecting all the
c-paths of faces of $\Psi$ that are selected in $\Psi$,
and selecting those c-paths of $\Theta$ whose images in $\Phi$
coincide with images of selected c-paths of faces that
are in $\Delta(2)\setminus\Phi(2)$.

Since $\Delta$ satisfies $\mathsf Z(2)$
relative to every proper simple disc submap of $\Phi$,
so does~$\bar\Phi$.

Since $\Delta$ does not satisfy $\mathsf Z(2)$ relative to $\Phi$,
it follows that $\kappa_{\bar\Phi}'(\Theta)=0$ and
$\kappa_{\bar\Phi}(\Theta)\le2$.

Let $A'$ be the set of all the exceptional arcs of $\Delta$ 
that are internal in $\Phi$,
and $A''$ be the set of all the exceptional arcs of $\Delta$ 
that are external in~$\Phi$.

Let $A$ be a set of pair-wise non-overlapping
selected (internal) arcs of $\bar\Phi$ such that
every selected edge of $\bar\Phi$ lies on an element of $A$,
every element of $A'\sqcup A''$ lies on an element of $A$,
and the cardinality of $A$ is the minimal possible.
Then it is easy to see that no two distinct elements of $A$
are subarcs of the same selected arc.

Consider a special case: suppose that $\bar\Phi$ is an elementary map
in which all c-paths are selected.
This implies that $\Delta$ itself is an elementary map
in which all c-paths are selected.
Then $\Phi$ has a single face $\Pi$, 
the set $A$ consists of a single element $u$, 
and one of the oriented arcs representing $u$
is a cyclic shift of $\cntr\Pi$.
Hence, as follows from $\mathsf D_2(\mu)$ and $\mathsf D_3(\nu)$,
$$
\abs{\cntr\Pi}=\abs{u}\le(\mu(\Pi)+\nu(\Pi))\abs{\cntr\Pi}
<\frac{1}{2}\abs{\cntr\Pi}.
$$
This gives a contradiction, and hence either 
$\bar\Phi$ is non-elementary, 
or at least it has a maximal selected c-path.
Therefore, Estimating Lemma~\ref{estimating_lemma:selected_arcs}
can be applied.

Apply Estimating Lemma~\ref{estimating_lemma:selected_arcs}
to $\bar\Phi$, $A$, $\Phi(2)$ (in the role of the set $B$), 
$\bar\Phi(2)$ (in the role of the set $C$),
and $\{\Theta\}$ (in the role of the set~$D$).
Let $f$ be a function $A\to\Phi(2)$ such that:
\begin{enumerate}
\item
    for every $x\in A$, $f(x)$ is incident to $x$, and
\item
    for every $y\in\Phi(2)$,
    the full pre-image of $y$ under $f$ consists of at most
    $3+\kappa_{\bar\Phi}(y)+\kappa_{\bar\Phi}'(y)$ elements.
\end{enumerate}
(Since
$1+\kappa_{\bar\Phi}(\Theta)+\kappa_{\bar\Phi}'(\Theta)+2-3\chi_{\bar\Phi}
\le-1\le0$,
the ``set $E$'' is empty.)

For every $j$, let $B_j$ be the the set of all rank-$j$ faces of $\Phi$,
and $A_j'$ be the the set of all rank-$j$ elements of $A'$.
As in Estimating Lemma~\ref{estimating_lemma:exceptional_arcs},
for every $j$, let $\varepsilon(j)=1$ if $A''$ has an element of rank $j$,
and $\varepsilon(j)=0$ otherwise.

By Estimating Lemma~\ref{estimating_lemma:exceptional_arcs} 
applied to $\Phi$,
$$
\norm{A_j'}\le\max\{0,2\norm{B_j}-\varepsilon(j)-1\}
\qquad\text{for every $j$}.
$$

Let $p_1$ and $p_2$ be paths such that
$\langle p_1p_2\rangle=\cntrcclof{{}\Phi}$,
$p_1$ is the image of a selected c-path of some face
$\Pi_1\in\Delta(2)\setminus\Phi(2)$,
and $p_2$ either is trivial, 
or is the image of a selected c-path of some
$\Pi_2\in\Delta(2)\setminus\Phi(2)$
($\Pi_1$ and $\Pi_2$ are not assumed to be distinct).
Moreover, choose such paths $p_1$ and $p_2$ so that 
every element of $A''$ lie on one of them.
Such paths $p_1$ and $p_2$ exist because $\Delta$ does not satisfy
$\mathsf Z(2)$ relative to~$\Phi$.

For $i=1,2$, let $A^{\prime\prime(i)}$ be the set of those 
elements of $A''$ that lie on~$p_i$.
Clearly, for each $i$, all elements of $A^{\prime\prime(i)}$
have the same rank.
If $A^{\prime\prime(i)}\ne\varnothing$ and $j$
is the rank of elements of $A^{\prime\prime(i)}$, then
$\varepsilon(j)=1$.

Let $J$ be the set of ranks of all elements of $A'\sqcup A''$.
For every $j\in J$, let 
$$
n(j)=\min_{\Pi\in B_j}\nu(\Pi)\abs{\cntr\Pi}.
$$
Then, as follows from $\mathsf D_3(\nu)$,
$$
\sum_{x\in A^{\prime\prime(i)}}\!\abs{x}
\le\sum_{j\in J}\varepsilon(j)n(j)\le\sum_{j\in J}n(j)
\qquad\text{for $i=1,2$}.
$$

Estimate the total number of edges of all the elements 
of $A'\sqcup A''$.
Denote this number by $N$.
By $\mathsf D_3(\nu)$, obtain: 
\begin{align*}
N&=\sum_{j\in J}\sum_{x\in A_j'}\!\abs{x}
+\sum_{x\in A^{\prime\prime(1)}}\!\abs{x}
+\sum_{x\in A^{\prime\prime(2)}}\!\abs{x}\\
&\le\sum_{j\in J}(2\norm{B_j}-\varepsilon(j)-1)n(j)
+\sum_{j\in J}\varepsilon(j)n(j)+\sum_{j\in J}n(j)\\
&=\sum_{j\in J}2\norm{B_j}n(j)
\le\sum_{\Pi\in\Phi(2)}\!2\nu(\Pi)\abs{\cntr\Pi}.
\end{align*}

Estimate the total number of the edges of elements of $A$
that are not edges of elements of $A'\sqcup A''$.
Denote this number by $M$.
By $\mathsf D_2(\mu)$, obtain:
$$
M=\sum_{\Pi\in\Phi(2)}\sum_{x:f(x)=\Pi}\!\abs{x}
\le\sum_{\Pi\in\Phi(2)}\!
(3+\kappa(\Pi)+\kappa'(\Pi))\mu(\Pi)\abs{\cntr\Pi}.
$$

Estimate the total number of the edges of $\Phi$
that are not edges of elements of $A$.
Denote this number by $L$.
By $\mathsf D_1(\lambda)$, obtain:
$$
L\le\sum_{\Pi\in\Phi(2)}\!\lambda(\Pi)\abs{\cntr\Pi}.
$$

Thus, on one hand,
\begin{align*}
\norm{\Phi(1)}&=L+M+N\\
&\le\sum_{\Pi\in\Phi(2)}\!
\bigl(\lambda(\Pi)
+(3+\kappa(\Pi)+\kappa'(\Pi))\mu(\Pi)
+2\nu(\Pi)\bigr)\abs{\cntr\Pi}\\
&\le\sum_{\Pi\in\Phi(2)}\frac{1}{2}\abs{\cntr\Pi};
\end{align*}
on the other hand,
$$
\norm{\Phi(1)}=\frac{1}{2}\sum_{\Pi\in\Phi(2)}\!\abs{\cntr\Pi}
+\frac{1}{2}\abs{\cntr{}\Phi}.
$$
This gives a contradiction.
\end{proof}


\begin{lemma}
\label{lemma:main_1}
Let\/ $\Delta$ be an\/ $S$-map\textup, $c_\Delta$ be the number
of contours of\/ $\Delta$\textup.
Suppose\/ $c_\Delta+3\chi_\Delta\ge0$\textup.
Assume\/ $\Delta$ satisfies the conditions\/ $\mathsf Y$
and\/ $\mathsf D(\lambda,\mu,\nu)$ \textup(absolutely\/\textup{).}
Let\/ $\gamma=\lambda+(3+\kappa_\Delta+\kappa_\Delta')\mu+2\nu$\textup.
Suppose\/ $\gamma(\Pi)\le1/2$ for every face\/ $\Pi$ of\/ $\Delta$\textup.
Let\/ $T$ be the set of all the edges of\/ $\Delta$ that are
incident to faces\textup. 
Let\/ $S$ be the set of all those elements of\/ $T$ that are external
edges of $\Delta$ and are the images of selected c-edges\textup,
Then
$$
\sum_i\abs{\cntr_i\Delta}\ge\norm{S}
\ge\norm{T}-\sum_{\Pi\in\Delta(2)}\gamma(\Pi)\abs{\cntr\Pi}
\ge\sum_{\Pi\in\Delta(2)}(1-2\gamma(\Pi))\abs{\cntr\Pi}.
$$
\end{lemma}


\begin{proof}
It suffices to prove that
$$
\norm{S}\ge\norm{T}-\sum_{\Pi\in\Delta(2)}\gamma(\Pi)\abs{\cntr\Pi}.
$$
One of the other two inequalities is obvious, and the other follows
form a simple computation similar to that in 
Remark 6.1 of \cite{Muranov:2005:dsmcbgbsg} or
in Proposition 4.1 of \cite{Muranov:2007:otfgfrfb}.

Using induction and Inductive Lemma, obtain that $\Delta$ satisfies
the condition $\mathsf Z(2)$ relative to every simple disc submap.

Let $N$ be the sum of the lengths of all the internal exceptional arcs
of $\Delta$,
$M$ be the number of the selected internal edges of $\Delta$
that do not belong to any exceptional arc, and
$L$ be the number of (non-selected) edges of $\Delta$ that
are the images of non-selected c-edges.
Then
$$
\norm{S}=\norm{T}-L-M-N.
$$

Using Estimating Lemma~\ref{estimating_lemma:exceptional_arcs} and 
the condition $\mathsf D_3(\nu)$,
obtain:
$$
N\le\sum_{\Pi\in\Delta(2)}\!2\nu(\Pi)\abs{\cntr\Pi}.
$$
Similarly to the proof of Inductive Lemma 
(but with no need for using a closure of $\Delta$) obtain that
$$
M\le\sum_{\Pi\in\Delta(2)}\!
(3+\kappa_\Delta(\Pi)+\kappa_\Delta'(\Pi))\mu(\Pi)\abs{\cntr\Pi}
$$
(using Estimating Lemma~\ref{estimating_lemma:selected_arcs} and 
the condition $\mathsf D_2(\mu)$),
and
$$
L\le\sum_{\Pi\in\Delta(2)}\!\lambda(\Pi)\abs{\cntr\Pi}
$$
(using the condition $\mathsf D_1(\lambda)$).
Therefore, 
$$
L+M+N\le\sum_{\Pi\in\Delta(2)}\gamma(\Pi)\abs{\cntr\Pi},
$$
which completes the proof.
\end{proof}



\section{Asphericity and torsion}
\label{section:asphericity_torsion}

Definitions of
\emph{aspherical\/} (\emph{A\/}),
\emph{combinatorially aspherical\/} (\emph{CA\/}),
\emph{diagrammatically aspherical\/} (\emph{DA\/}),
\emph{singularly aspherical\/} (\emph{SA\/}), and
\emph{Cohen-Lyndon aspherical\/} (\emph{CLA\/})
presentations may be found in \cite{ChiColHue:1981:agp}.
It should be noted that none of these definitions
requires the set of relators to consist of only reduced elements.
Moreover, \emph{group presentations\/} are regarded in a way
that a priori allows for repetition of relators
(instead of sets of relators,
presentations have \emph{indexed families\/} of relators).
Only diagrammatic and singular asphericities shall be used in this paper.

The version of asphericity defined in
\cite{Olshanskii:1989:gosg-rus,Olshanskii:1991:gdrg-eng}
is equivalent to diagrammatic asphericity by Theorem~32.2 therein.

The following is another equivalent definition of diagrammatic
asphericity:

\begin{definition}
A group presentation is \emph{diagrammatically aspherical\/} if
every spherical diagram over this presentation can be transformed
by a sequence of diamond moves into a diagram whose
all connected components are elementary spherical diagrams. 
\end{definition}

Proof of equivalence is left to the reader.

\begin{definition}
A group presentation \GP{\mathfrak A}{\mathcal R} is 
\emph{singularly aspherical\/} if
it is diagrammatically aspherical, no element of $\mathcal R$
represents a proper power in the free group \GP{\mathfrak A}{\varnothing},
and no two distinct elements of $\mathcal R$
are conjugate or conjugate to each other's inverses
in~\GP{\mathfrak A}{\varnothing}. 
\end{definition}

\begin{definition}
Call a group (\emph{A\/}), (\emph{CA\/}), (\emph{DA\/}), 
(\emph{SA\/}), or (\emph{CLA\/}), accordingly,
if
it has a presentation which is such.
\end{definition}

Interesting results on relations between different concepts
of asphericity
(of which, by the way, combinatorial asphericity is the weakest,
and singular asphericity is in a sense one of the strongest),
and classification of torsion elements
in combinatorially aspherical groups may be found in
\cite{Huebschmann:1979:ctagscg,ChiColHue:1981:agp}.
Combinatorial asphericity is also discussed in
\cite{Huebschmann:1980:htcap} in great detail.


\begin{remark}
It follows directly from the definition of singular asphericity
that a group is singularly aspherical if and only if
it has a diagrammatically aspherical presentation without proper
powers among relators,
relators being viewed as elements of the free group on the set
of generators.
\end{remark}

\begin{lemma}
\label{lemma:singularly_aspherical_groups}
Singularly aspherical groups are torsion-free\textup.
\end{lemma}

\begin{proof}
Let $G$ be an arbitrary singularly aspherical group.
Let \GP{\mathfrak A}{\mathcal R}
be a singularly aspherical presentation of $G$.
Then the relation module $M$ of
\GP{\mathfrak A}{\mathcal R}
is a free $G$-module by Corollary~32.1
in \cite{Olshanskii:1989:gosg-rus,Olshanskii:1991:gdrg-eng}.
Therefore, there exists a finite-length free resolution of $\mathbb Z$
over $\mathbb ZG$:
$$
0\to M\to\bigoplus_{x\in\mathfrak A}\mathbb ZG
\to\mathbb ZG\to\mathbb Z\to0,
$$
where $\mathbb ZG$ and $\bigoplus_{x\in\mathfrak A}\mathbb ZG$
are identified with the (free) $G$-modules of, respectively,
$0$- and $1$-dimensional cellular chains of the Cayley complex of
\GP{\mathfrak A}{\mathcal R}.

Suppose now that $G$ has torsion.
Let $C$ be a nontrivial finite cyclic subgroup of~$G$.
Every free $G$-module may be naturally regarded as a free $C$-module.
Hence the above resolution may be viewed
a free resolution of $\mathbb Z$ over~$\mathbb ZC$.
This contradicts the fact that all odd-dimensional homology groups
of any nontrivial finite cyclic group are nontrivial
(see~\cite{Brown:1994:cg}).
\end{proof}



\section{Proof of the theorems}
\label{section:proof_theorems}

Theorems \ref{theorem:bsglcw} and \ref{theorem:sgicw}
are proved in this section by parallel series of similar arguments.
It is convenient in both cases to use the notation of
Section~\ref{section:group_presentations}.

There is a conflict in notation between Subsections
\ref{subsection:bsglcw.group_presentations} and
\ref{subsection:sgicw.group_presentations},
but it shall not cause confusion if
the notation of Subsection~\ref{subsection:bsglcw.group_presentations}
is used only in the context of proving Theorem~\ref{theorem:bsglcw},
while considering \GP{\mathfrak A}{\mathcal R_n}, $n\in\mathbb N$, and
the notation of Subsection~\ref{subsection:sgicw.group_presentations}
is used only in the context of proving Theorem~\ref{theorem:sgicw},
while considering \GP{\mathfrak A}{\mathcal R_\infty}.

It is convenient to assume in this section that to every diagram
under consideration there is assigned
a \emph{sort\/} which is either $(\mathrm{I}.n)$, $n\in\mathbb N$,
or $(\mathrm{II})$.
More precisely, every diagram or presentation considered in this section
is always equipped with a \emph{sort attribute}.
Diagrams of sorts $(\mathrm{I}.n)$, $n\in\mathbb N$,
will be used for proving Theorem~\ref{theorem:bsglcw},
and most diagrams of sort $(\mathrm{I}.n)$ under consideration
will be diagrams over \GP{\mathfrak A}{\mathcal R_n}.
Similarly, diagrams of sort $(\mathrm{II})$
will be used for proving Theorem~\ref{theorem:sgicw},
and most diagrams of this sort under consideration
will be diagrams over \GP{\mathfrak A}{\mathcal R_\infty}.
Whenever the sort is not assigned explicitly, it shall be
assumed in the most natural way,
but a priori the sort is not determined by the diagram itself.
The purpose of this convention is to unambiguously use the
same term in relation to a diagram in different senses depending
on the context (on the sort of the diagram).

Define sets of indices $I_n$, $n\in\mathbb N$, and $I_\infty$ as follows:
$$
I_n=\{\,j\,|\,\mathcal R_n^{(j)}\ne\mathcal R_n^{(j-1)}\,\},\qquad
I_\infty=\{\,j\,|\,\mathcal R_\infty^{(j)}\ne\mathcal R_\infty^{(j-1)}\,\}.
$$
Then
$$
\mathcal R_n=\{\,r_j\,|\,j\in I_n\,\},\qquad
\mathcal R_\infty=\{\,r_{j,1},r_{j,2}\,|\,j\in I_\infty\,\}.
$$

Let $W_n=\{\,w_j\,|\,j\in I_n\,\}$ where $w_j$ are the group words defined
in Subsection \ref{subsection:bsglcw.group_presentations}.
Let $W_\infty=\{\,w_j\,|\,j\in I_\infty\,\}$ where $w_j$ are
the group words defined
in Subsection~\ref{subsection:sgicw.group_presentations}.

\begin{definition}
A graded $S_1$-diagram $\Delta$ of sort $(\mathrm{I}.n)$
is called \emph{correct\/} if
\begin{enumerate}
\item
    the rank of every face of $\Delta$ is in $I_n\sqcup\{0\}$;
\item
    for every face $\Pi$ of rank $j$,
    $\pathlbl(\cntr\Pi)\in\{r_j^{\pm1}\}$
    if $j\in I_n$,
    and $\pathlbl(\cntr\Pi)$ is the concatenation of several copies of
    $z_1^{\pm1}$ and $z_2^{\pm1}$ if $j=0$;
\item
    for every face $\Pi$ of rank $j\ne0$,
    $\Pi$ has c-paths
    $s_1$, \dots, $s_{2n+2}$, $s_1'$, \dots, $s_{2n+2}'$,
    $t_1$, \dots, $t_{2n+2}$, and $t_0$ such that:
    \begin{enumerate}
    \item
        $s_1t_1s_1's_2t_2s_2'\dots s_{2n+2}t_{2n+2}s_{2n+2}'t_0
        \in\{\ccntr\Pi,(\ccntr\Pi)^{-1}\}$,
    \item
        \begin{enumerate}
        \item
            $\pathlbl(s_i)=u_{ji}$ and $\pathlbl(s_i')=u_{ji}^{-1}$
            for every $i=1,\dots,2n+2$,
        \item
            $\pathlbl(t_1)=\dots=\pathlbl(t_{2n+2})=w_j$,
        \item
            $\pathlbl(t_0)=v_j^{-1}$,
        \end{enumerate}
    \item
        a c-path of $\Pi$ is selected if and only if
        it is a nontrivial subpath of one of the following
        $8n+8$ paths: $s_1^{\pm1}$, \dots, $s_{2n+2}^{\pm1}$,
        $s_1^{\prime\pm1}$, \dots, $s_{2n+2}^{\prime\pm1}$
    \end{enumerate}
    (in particular, $\kappa_\Delta(\Pi)=4n+4$);
\item
    for every face $\Pi$ of rank $0$,
    all c-pseudo-arcs of $\Pi$ are selected
    (in particular, $\kappa_\Delta(\Pi)=0$);
\item
    if two faces are congruent
    (\ie, their contour labels are cyclic shifts of each other
    or cyclic shifts of the inverses of each other), then either
    these faces have the same rank, or
    the rank of at least one of these faces is~$0$.
\end{enumerate}
\end{definition}

\begin{definition}
A graded $S_1$-diagram $\Delta$ of sort $(\mathrm{II})$
is called \emph{correct\/} if
\begin{enumerate}
\item
    the rank of every face of $\Delta$ is in $I_\infty\sqcup\{0\}$;
\item
    for every face $\Pi$ of rank $j$,
    $\pathlbl(\cntr\Pi)
    \in\{r_{j,1}^{\pm1},r_{j,2}^{\pm1}\}$
    if $j\in I_\infty$,
    and $\pathlbl(\cntr\Pi)$ is the concatenation of several copies of
    $z_1^{\pm1}$ and $z_2^{\pm1}$ if $j=0$;
\item
    for every face $\Pi$ of rank $j\ne0$,
    $\Pi$ has c-paths $q$,
    $s_1$, \dots, $s_{2j+2}$, $s_1'$, \dots, $s_{2j+2}'$,
    $t_1$, \dots, $t_{2j+2}$, and $t_0$ such that:
    \begin{enumerate}
    \item
        $s_1t_1s_1's_2t_2s_2'\dots s_{2j+2}t_{2j+2}s_{2j+2}'t_0
        =q\in\{\ccntr\Pi,(\ccntr\Pi)^{-1}\}$,
    \item
        either
        \begin{enumerate}
        \item
            $\pathlbl(q)=r_{j,1}$,
        \item
            $\pathlbl(s_i)=u_{ji}$ and $\pathlbl(s_i')=u_{ji}^{-1}$
            for every $i=1,\dots,2j+2$,
        \item
            $\pathlbl(t_1)=\dots=\pathlbl(t_{2j+2})=w_j$,
        \item
            $\pathlbl(t_0)=a^{-1}$,
        \end{enumerate}
        or
        \begin{enumerate}
        \item
            $\pathlbl(q)=r_{j,2}$,
        \item
            $\pathlbl(s_i)=u_{j,2j+2+i}$ and
            $\pathlbl(s_i')=u_{j,2j+2+i}^{-1}$
            for every $i=1,\dots,2j+2$,
        \item
            $\pathlbl(t_1)=\dots=\pathlbl(t_{2j+2})=w_j$,
        \item
            $\pathlbl(t_0)=b^{-1}$,
        \end{enumerate}
    \item
        a c-path of $\Pi$ is selected if and only if
        it is a nontrivial subpath of one of the following
        $8j+8$ paths: $s_1^{\pm1}$, \dots, $s_{2j+2}^{\pm1}$,
        $s_1^{\prime\pm1}$, \dots, $s_{2j+2}^{\prime\pm1}$
    \end{enumerate}
    (hence $\kappa_\Delta(\Pi)=4\rank(\Pi)+4$);
\item
    for every face $\Pi$ of rank $0$,
    all c-pseudo-arcs of $\Pi$ are selected
    (in particular, $\kappa_\Delta(\Pi)=0$);
\item
    if two faces are congruent, then either
    they have the same rank, or
    at least one of them has rank~$0$.
\end{enumerate}
\end{definition}

(The last conditions in these two definitions may be redundant,
but are easy to satisfy, and 
they facilitate the proof of Lemma~\ref{lemma:condition_Y}.)

\begin{definition}
Faces of rank $0$ in a correct graded $S_1$-diagram of any sort
are called \emph{alien}, all the other faces are called \emph{native}.
Correct graded $S_1$-diagrams without alien faces are called
\emph{restricted}.
\end{definition}

For every diagram over \GP{\mathfrak A}{\mathcal R_n},
and for every diagram over \GP{\mathfrak A}{\mathcal R_\infty},
there is an essentially isomorphic diagram
that has a selection and a grading which turn it into a restricted correct
graded $S_1$-diagram of sort $(\mathrm{I}.n)$ or $(\mathrm{II})$,
respectively.

\begin{definition}
An $S$-diagram $\Delta$ of any sort is called \emph{correct\/} if
\begin{enumerate}
\item
    $\Delta$ is correct as a graded $S_1$-diagram;
\item
    an internal arc $u$ of $\Delta$ is exceptional if an only if
    there exist c-paths $s_1$, $s_{1-}$, $s_{10}$, $s_{1+}$,
    $s_2$, $s_{2-}$, $s_{20}$, $s_{2+}$ such that:
    \begin{enumerate}
    \item
        $s_1=s_{1-}s_{10}s_{1+}$ and $s_2=s_{2-}s_{20}s_{2+}$,
    \item
        $s_1$ and $s_2$ are maximal selected c-paths,
    \item
        $s_{10}$ and $s_{20}$ are distinct c-paths with a common
        image (in the $1$-skeleton of $\Delta$) which coincides with
        one of the oriented arcs associated with $u$, and
    \item
        $\pathlbl(s_1)=\pathlbl(s_2)$, $\abs{s_{1-}}=\abs{s_{2-}}$,
        $\abs{s_{1+}}=\abs{s_{2+}}$;
    \end{enumerate}
\item
    exceptional arcs are not incident to faces of rank~$0$
\end{enumerate}
(in particular, every internal exceptional arc of $\Delta$ is
a maximal selected arc, and there are no exceptional arcs of
rank~$0$).
\end{definition}

Every correct graded $S_1$-diagram of any sort has
a structure of a correct $S$-diagram that extends the given
structure of a graded $S_1$-diagram and is unique
up to choice of external exceptional arcs.
Every maximal selected internal arc of a correct $S$-diagram
either is exceptional, or does not overlap with any exceptional arc.

\begin{definition}
An internal exceptional arc $u$ of a correct $S$-diagram
is called \emph{non-extendible\/} if
it is the image of two maximal selected c-arcs;
otherwise $u$ is called \emph{extendible}.
\end{definition}

Every extendible exceptional arc of a correct $S$-diagram
of any sort can be ``extended'' to a longer
exceptional arc by a diamond move.
(Diamond moves are viewed here as operations on
correct $S$-diagrams of a given sort.)

\begin{definition}
An $S$-diagram $\Delta$ of any sort is called \emph{special\/} if
it is correct, weakly reduced, and has no
extendible internal exceptional arcs.
\end{definition}

Clearly, every $S$-subdiagram of every special $S$-diagram
is special.


\begin{lemma}
\label{lemma:special_diagrams}
Every correct\/ $S$-diagram of any sort
can be transformed by a series of
diamond moves into an\/ $S$-diagram each connected component of which
is either special and reduced\textup, or elementary spherical\textup.
In particular\textup, every reduced correct\/ $S$-diagram
can be transformed by diamond moves into a special\/ $S$-diagram\textup.
\end{lemma}


Note that even when the diagrams under consideration are graded,
the property of being reduced is the same as for non-graded ones,
unlike~\cite{Olshanskii:1989:gosg-rus,Olshanskii:1991:gdrg-eng}.


\begin{proof}[Proof of the lemma]
The number of connected components and the Euler
characteristic of any map that can be obtained from a given map
$\Gamma$ by diamond moves are bounded from above.
Indeed, the number of connected components is bounded by
$\norm{\Gamma(2)}+c_\Gamma$, and the Euler characteristic is bounded by
$2\norm{\Gamma(2)}+c_\Gamma$, as follows from
Lemma~\ref{lemma:maps__positive_Euler_characteristic}.
Here $c_\Gamma$ denotes the number of contours of~$\Gamma$.

Since diamond moves do not decrease the Euler characteristic,
and improper diamond moves increase it
(see Lemma~\ref{lemma:diamond_move}),
it follows that in any sequence
of diamond moves applied to a given diagram,
there is only bounded number of improper ones.

Consider an arbitrary correct $S$-diagram $\Delta$.
Assume without loss of generality that no improper
diamond move is applicable to $\Delta$, nor to any
(correct) $S$-diagram obtained from $\Delta$
by proper diamond moves.
In particular, neither the number of connected components of $\Delta$,
nor the Euler characteristic of $\Delta$ can be
increased by any sequence of diamond moves.
Then every connected component of $\Delta$ either is reduced or
otherwise can be turned into an elementary spherical diagram
by a sequence of proper diamond moves.
Thus it is left to show that every reduced connected component
of $\Delta$ can be made special by (proper) diamond moves.

Let $\Psi$ be a reduced connected component of $\Delta$.
Diamond moves allow one to ``extend'' all extendible
internal exceptional arcs one-by-one.
Any $S$-diagram obtained from $\Psi$ in this manner
will be special.
\end{proof}


For every $n$, denote $1/(4n+4)$ by~$\nu_n$.

If $\Delta$ is a restricted special $S$-diagram of sort $(\mathrm{I}.n)$,
then let $\lambda_\Delta$, $\mu_\Delta$, and $\nu_\Delta$
be the constant functions on $\Delta(2)$ defined as follows:
$$
\lambda_\Delta=\lambda_n,\quad
\mu_\Delta=\mu_n,\quad
\nu_\Delta=\nu_n=\frac{1}{4n+4}.
$$
Then, as follows from \thetag{\ref{display:main_inequality}},
\begin{equation}
\label{display:main_inequality_for_bsglcw.ineq}
2\lambda_\Delta+(2\kappa_\Delta+6n)\mu_\Delta+(2n+1)\nu_\Delta
<\frac{1}{2}.
\end{equation}

If $\Delta$ is a restricted special $S$-diagram of sort $(\mathrm{II})$,
then let $\lambda_\Delta$, $\mu_\Delta$, and $\nu_\Delta$
be the functions on $\Delta(2)$ defined as follows:
$$
\lambda_\Delta(\Pi)=\lambda_{\rank(\Pi)},\quad
\mu_\Delta(\Pi)=\mu_{\rank(\Pi)},\quad
\nu_\Delta(\Pi)=\nu_{\rank(\Pi)}.
$$
Then, as follows from \thetag{\ref{display:main_inequality}},
\begin{equation}
\label{display:main_inequality_for_sgicw.ineq}
2\lambda_\Delta+(2\kappa_\Delta+6\rank)\mu_\Delta
+(2\rank+1)\nu_\Delta
<\frac{1}{2}.
\end{equation}


\begin{lemma}
\label{lemma:condition_D}
Let\/ $\Delta$ be an arbitrary correct and weakly reduced\/
\textup(or special\/\textup) $S$-diagram\textup, and\/
$\Gamma$ be a restricted\/ $S$-subdiagram\textup.
Then\/ $\Delta$ satisfies the condition\/
$\mathsf D(\lambda_\Gamma,\mu_\Gamma,\nu_\Gamma)$
relative to\/ $\Gamma$\textup.
In the case\/ $\Delta$ is of sort\/ $(\mathrm{II})$\textup,
it satisfies\/ $\mathsf D'(\lambda_\Gamma,\mu_\Gamma,\nu_\Gamma)$
relative to\/~$\Gamma$\textup.
%
\end{lemma}


This lemma follows directly from the construction of the presentations
\GP{\mathfrak A}{\mathcal R_n}, $n\in\mathbb N$, and
\GP{\mathfrak A}{\mathcal R_\infty}
in Section~\ref{section:group_presentations}.


\begin{lemma}
\label{lemma:words_w}
Let\/ $j$ be a natural number and\/ $\Delta$ be a
restricted special disc\/ $S$-diagram satisfying the condition\/
$\mathsf Y$ and containing a face of rank at least\/ $j$\textup.
Let\/ $w=w_j$ where\/ $w_j$ is defined
in Subsection \/\textup{\ref{subsection:bsglcw.group_presentations}} or\/
\textup{\ref{subsection:sgicw.group_presentations}}
depending on the sort of\/ $\Delta$\textup.
Then\/ $\abs{\cntr{}\Delta}>\abs{w}$\textup.
\end{lemma}


\begin{proof}
By Lemma~\ref{lemma:condition_D}, the $S$-map $\Delta$ satisfies
the condition $\mathsf D(\lambda_\Delta,\mu_\Delta,\nu_\Delta)$
absolutely.

Let
$$
\gamma
=\lambda_\Delta+(3+\kappa_\Delta+\kappa_\Delta')\mu_\Delta+2\nu_\Delta.
$$
It follows from inequality
\thetag{\ref{display:main_inequality_for_bsglcw.ineq}} or
\thetag{\ref{display:main_inequality_for_sgicw.ineq}},
depending to the sort of $\Delta$, that
$\gamma(\Pi)<1/2-\lambda_\Delta(\Pi)$ for every face $\Pi$ of $\Delta$.

By Lemma~\ref{lemma:main_1},
$$
\abs{\cntr{}\Delta}
\ge\sum_{\Pi\in\Delta(2)}(1-2\gamma(\Pi))\abs{\cntr\Pi}
>\sum_{\Pi\in\Delta(2)}\!2\lambda_\Delta(\Pi)\abs{\cntr\Pi}.
$$

Let $\Pi$ be a face of $\Delta$ of rank $\ge j$.
It follows from one of the conditions imposed on the group presentations 
in Section~\ref{section:group_presentations} that
$\abs{w}\le\lambda_\Delta(\Pi)\abs{\cntr\Pi}$.
Hence $\abs{\cntr{}\Delta}>\abs{w}$.
\end{proof}


\begin{lemma}
\label{lemma:condition_Y}
Let\/ $\Delta$ be a connected restricted special\/ $S$-diagram\/
\textup(of any sort\/\textup) over
a group presentation\/ \GP{\mathfrak A}{\mathcal S}\textup.
Suppose that every proper\/ \textup(finite\/\textup) subpresentation
of\/ \GP{\mathfrak A}{\mathcal S} defines a torsion-free group\textup.
Then\/ $\Delta$ satisfies the condition\/~$\mathsf Y$\textup.
\end{lemma}


\begin{proof}
Without loss of generality, assume that every
connected special $S$-diagram over \GP{\mathfrak A}{\mathcal S}
whose set of face ranks is a proper subset of the
set of face ranks of $\Delta$, does satisfy $\mathsf Y$.
(Alternatively, one can induct on the number of different face ranks
of a diagram.)

Let $j$ be the rank of an arbitrary internal exceptional arc of $\Delta$.
Let $A$ be the set of all the internal exceptional arcs of $\Delta$
of rank $j$, and $B$ be the set of all the faces of $\Delta$ of rank~$j$.
Let $\Gamma$ be the $S$-subdiagram obtained
from $\Delta$ by removing all the faces and internal exceptional
arcs of rank $j$.

Let $\mathcal T$ be the minimal subset of $\mathcal S$ such that
$\Gamma$ is a diagram over \GP{\mathfrak A}{\mathcal T}. 
The set of face ranks of $\Gamma$ is a proper subset of that of $\Delta$;
therefore, $\mathcal T$ is a proper subset of $\mathcal S$.
Hence the groups presented by \GP{\mathfrak A}{\mathcal T}
is torsion-free.

Let $k=\kappa_\Delta(\Pi)$ for an arbitrary face $\Pi$ of rank $j$
($k=4n+4$ if the sort is $(\mathrm{I}.n)$,
and $k=4j+4$ if the sort is $(\mathrm{II})$).
Let $C$ be the set of corners of elements of $B$ chosen as follows:
If $\Pi$ is an arbitrary element of $B$,
let $s_1$, \dots, $s_{k/2}$, $s_1'$, \dots, $s_{k/2}'$,
$t_1$, \dots, $t_{k/2}$, and $t_0$ be the c-paths of $\Pi$
such as in the definition of correct graded $S_1$-diagrams
($s_1$, \dots, $s_{k/2}$, $s_1'$, \dots, $s_{k/2}'$ are maximal selected);
let $C$ contain the initial corners of the c-paths
$s_1t_1s_1'$, \dots, $s_{k/2}t_{k/2}s_{k/2}'$,
and no other corners of $\Pi$.
Thus $C$ contains exactly $k/2$ corners of each element of~$B$.

The image (under the attaching morphism) of each element of $C$
is a vertex of $\Gamma$.
As follows from Lemma~\ref{lemma:maps__positive_Euler_characteristic},
every connected component of $\Gamma$ of Euler characteristic $1$ is
a disc $S$-diagram.
It suffices to prove now that for every connected component $\Psi$
of $\Gamma$, if either $\Psi$ is a disc submap, or it contains the image
of a selected c-edge of a face of $\Delta$ of rank $j$, then
the number of elements of $C$ whose images
are in $\Psi$ is at least~$k/2$.

Since $\Delta$ is weakly reduced, it follows from the definition of
exceptional arcs in a correct $S$-diagram that for every connected
component $\Psi$ of $\Gamma$,
the number of elements of $C$ whose images are in $\Psi$
is divisible by $k/2$
(recall Lemma~\ref{lemma:submap}).

Consider an arbitrary connected component $\Psi$ of $\Gamma$ which
contains the image of a selected c-edge of an element of $B$.
Let $\Pi$ be an element of $B$ and $x$ be a selected c-edge of $\Pi$
such that the image of $x$ is in $\Psi$.
Let $s$ be the maximal selected c-path of $\Pi$ such that
$x$ lies on $s$, and the label of $s$ is of the form $u_{ji}$
($i\in\{1,\dots,2n+2\}$ if the diagrams are of sort $(\mathrm{I}.n)$,
and $i\in\{1,\dots,4j+4\}$ if the diagrams are of sort $(\mathrm{II})$).
Let $c$ be the element of $C$ which is the corner of $\Pi$
``closest'' to the initial vertex of $s$.
More precisely, $c$ is the only element of $C$ which is a corner
of $\Pi$ and either coincides with the initial vertex of $s$,
or can be connected to the initial vertex of $s$ by a path
(c-path of $\Pi$) without selected oriented c-edges.
Then the image of $c$ is in $\Psi$, and hence the number of elements
of $C$ mapped to $\Phi$ is at least $k/2$
(since it is not $0$ and is divisible by~$k/2$).

Now consider an arbitrary connected component $\Psi$ of $\Gamma$
endowed with the inherited structure of a special disc $S$-diagram.
Suppose $\Psi$ does not contain the image of any element of $C$.
It follows from Lemma~\ref{lemma:submap} and from $\Delta$'s
being weakly reduced that some cyclic shift of $(\cntr{}\Psi)^{\pm1}$
can be decomposed into the product of paths each of which
is labelled by $w_j$.
Therefore, $w_j$ represents a finite-order element in the
group presented by \GP{\mathfrak A}{\mathcal T}.
Since that group is torsion-free,
$w_j$ represents the identity element in it.

Let $\Phi$ be a restricted special disc $S$-diagram 
over \GP{\mathfrak A}{\mathcal T}
whose contour label is $w_j$, and whose set of face ranks is a subset
of that of $\Gamma$
(here Lemma~\ref{lemma:special_diagrams} is used to find such a
\emph{special\/} $\Phi$).
By the inductive assumption at the beginning of this proof,
$\Phi$ satisfies $\mathsf Y$.

By the construction of group presentations in
Section~\ref{section:group_presentations}, $w_j$ is not trivial
modulo \GP{\mathfrak A}{\mathcal R_n^{(j-1)}}
or \GP{\mathfrak A}{\mathcal R_\infty^{(j-1)}}, whichever corresponds
to the sort of the diagrams under consideration.
Therefore, $\Phi$ contains at least one face whose rank is $j$ or greater.
This contradicts Lemma~\ref{lemma:words_w}.

Thus $\Psi$ does contain the image of some element of $C$,
and hence $\Psi$ contains the images of at least $k/2$ element of~$C$.

On one hand, $\norm{C}=(k/2)\norm{B}$
(recall that $B$ is the set of all rank-$j$ faces).
On the other hand, all elements of $C$ can be distributed among
connected components of $\Gamma$ so that there are at least $k/2$
elements assigned to each component that either is disc, or
contains an external exceptional arc of $\Delta$ of rank $j$.
Therefore, the number of such components does not exceed the number
of faces of rank $j$.
The same is true for every $j$ such that $\Delta$ has an internal
exceptional arc of rank $j$.
Hence the condition~$\mathsf Y$.
\end{proof}


\begin{lemma}
\label{lemma:powers_among_relators}
No element of\/ $\bigcup_{n\in\mathbb N}\mathcal R_n\cup\mathcal R_\infty$
represents a proper power in the free group on\/ $\mathfrak A$\textup.
Distinct element of\/ $\mathcal R_n$\textup, $n\in\mathbb N$\textup, 
or of\/ $\mathcal R_\infty$ do not represent 
conjugate elements of the free group on\/ $\mathfrak A$\textup,
nor elements conjugate to each other's inverses\textup.
\end{lemma}


\begin{proof}
Most likely there is a straightforward way to prove these facts using
only the small-cancellation conditions imposed on the
(subwords of) defining relations
in Section~\ref{section:group_presentations},
or they can be obtained for free by imposing
additional restrictions on the defining relators 
of the constructed presentations.
Following is a proof which is more in the spirit of this paper.

Suppose $r\in\mathcal R_\infty$ or $r\in\mathcal R_n$ for some $n$,
and $r$ represents a proper power in the free group on $\mathfrak A$.
Then $r$ is freely conjugate to $s^m$ where $s$ is cyclically
reduced and $m>1$.
Let $\Phi$ be a special simple single-face disc $S$-diagram over
$\langle\,\mathfrak A\,\|\,\{r\}\,\rangle$
such that $\pathlbl(\cntr{}\Phi)=s^m$.

Let $\Delta$ be a special spherical $S$-diagram obtained from
two copies of $\Phi$ by attaching them to each other along their
contour cycles with a shift by $\abs{s}$ edges.
More precisely, let $\Phi_1$ and $\Phi_2$ be two copies of
the $S$-diagram $\Phi$.
Let $q_1=\cntr{}\Phi_1$, and
let $q_2$ be a cyclic shift of $\cntr{}\Phi_2$ by $\abs{s}$ edges
(in either direction).
Observe that $\pathlbl(q_1)=s^m=\pathlbl(q_2)$.
Let $\Delta$ be the correct spherical $S$-diagram obtained
by ``gluing'' $\Phi_1$ and $\Phi_2$ together along the pair of paths
$q_1$ and $q_2$.

Because of the shift in ``gluing'' the copies of $\Phi$,
the $S$-diagram $\Delta$ does not have any exceptional arcs,
and hence satisfies the condition $\mathsf Y$.
For the same reason, $\Delta$ satisfies the condition
$\mathsf D(\lambda_\Delta,\mu_\Delta,0)$ absolutely.

Let $\Pi_1$ and $\Pi_2$ be the two faces of $\Delta$.
Let $\hat\kappa=\kappa_\Delta(\Pi_1)=\kappa_\Delta(\Pi_2)$,
$\hat\lambda=\lambda_\Delta(\Pi_1)=\lambda_\Delta(\Pi_2)$, and
$\hat\mu=\mu_\Delta(\Pi_1)=\mu_\Delta(\Pi_2)$.

By Lemma~\ref{lemma:main_1} and inequalities
\thetag{\ref{display:main_inequality_for_bsglcw.ineq}}
and \thetag{\ref{display:main_inequality_for_sgicw.ineq}},
$$
0\ge\bigl(1-2\bigl(\hat\lambda+(3+2\hat\kappa)\hat\mu\bigr)\bigr)
\bigl(\abs{\cntr\Pi_1}+\abs{\cntr\Pi_2}\bigr)>0.
$$
This gives a contradiction.

Suppose two distinct defining relators of one of the constructed 
presentations represent conjugate elements of the free group
\GP{\mathfrak A}{\varnothing}.
This situation also gives rise to a special spherical $S$-diagram $\Delta$
without exceptional arcs and satisfying 
$\mathsf D(\lambda_\Delta,\mu_\Delta,0)$.
(Such $\Delta$ is also obtained by ``gluing'' together two single-face 
simple disc diagrams.)
This again leads to a contradiction with Lemma~\ref{lemma:main_1}.

The case of two relators conjugate to the inverses of each other
is dealt with similarly.
\end{proof}


\begin{lemma}
\label{lemma:condition_Y__singular_asphericity}
Let\/ \GP{\mathfrak A}{\mathcal S} be a finite subpresentation of\/
\GP{\mathfrak A}{\mathcal R_\infty} or of\/
\GP{\mathfrak A}{\mathcal R_n} for some\/ $n\in\mathbb N$\textup.
Then\/ \GP{\mathfrak A}{\mathcal S} is singularly aspherical\textup,
and every connected restricted special\/ $S$-diagram\/
\textup(of appropriate sort\/\textup) over\/ \GP{\mathfrak A}{\mathcal S}
satisfies the condition\/~$\mathsf Y$\textup.
\end{lemma}


\begin{proof}
Induction on $\mathcal S$:
if $\mathcal S=\varnothing$, then the conclusion is obvious;
assume $\mathcal S\ne\varnothing$, and the statement is true
for all proper subpresentations of \GP{\mathfrak A}{\mathcal S}.

By the inductive assumption and 
Lemma~\ref{lemma:singularly_aspherical_groups},
every proper subpresentation of \GP{\mathfrak A}{\mathcal S}
defines a torsion-free group.
By Lemma~\ref{lemma:condition_Y},
every connected restricted special $S$-diagram
over \GP{\mathfrak A}{\mathcal S} satisfies the condition~$\mathsf Y$.

It is left to show that \GP{\mathfrak A}{\mathcal S} is
singularly aspherical.
Suppose it is not.
Due to Lemma~\ref{lemma:powers_among_relators}, this means that
\GP{\mathfrak A}{\mathcal S} is not diagrammatically aspherical.

Let $\Delta_0$ be a restricted correct reduced spherical $S$-diagram
over \GP{\mathfrak A}{\mathcal S}
(it exists since the presentation is not diagrammatically aspherical).
Let $\Delta_1$ be a special $S$-diagram obtained from $\Delta_0$
by diamond moves (see Lemma~\ref{lemma:special_diagrams}).
Then, as follows from Lemmas \ref{lemma:diamond_move} and
\ref{lemma:maps__positive_Euler_characteristic},
every connected component of $\Delta_1$ is a reduced spherical diagram.
Let $\Delta$ be an arbitrary connected component of $\Delta_1$.
It is already shown that such $\Delta$
must satisfy the condition $\mathsf Y$.
By Lemma~\ref{lemma:condition_D}, $\Delta$ satisfies the condition
$\mathsf D(\lambda_\Delta,\mu_\Delta,\nu_\Delta)$.
By Lemma~\ref{lemma:main_1} and inequalities
\thetag{\ref{display:main_inequality_for_bsglcw.ineq}}
and \thetag{\ref{display:main_inequality_for_sgicw.ineq}},
$$
0\ge\sum_{\Pi\in\Delta(2)}
\bigl(1-2\bigl(\lambda_\Delta(\Pi)
+(3+2\kappa_\Delta(\Pi))\mu_\Delta(\Pi)
+2\nu_\Delta(\Pi)\bigl)\bigr)\abs{\cntr\Pi}>0,
$$
which gives a contradiction.
Thus \GP{\mathfrak A}{\mathcal S} is singularly aspherical.
\end{proof}


\begin{proposition}
\label{proposition:Gn_lcw}
For every\/ $n\in\mathbb N$\textup,
the group\/ $G_n$ constructed in
Subsection\/~\textup{\ref{subsection:bsglcw.group_presentations}}
is singularly aspherical\textup, torsion-free\textup, and
the elements\/ $[z_1]_{G_n}$ and $[z_2]_{G_n}$
freely generate a free subgroup\/ $H$ such that
$$
\bigl(\forall h\in H\setminus\{1\}\bigr)\,
\bigl(\forall m\ge2n\bigr)\,
\bigl(\cl_{G_n}(h^m)>n\bigr).
$$
\end{proposition}



\begin{proof}
By Lemma~\ref{lemma:condition_Y__singular_asphericity},
every finite subpresentation of \GP{\mathfrak A}{\mathcal R_n}
is singularly aspherical.
Therefore, \GP{\mathfrak A}{\mathcal R_n} itself is
singularly aspherical.
Therefore, by Lemma~\ref{lemma:singularly_aspherical_groups},
the group $G_n$ is torsion-free.

Let $w$ be an arbitrary nontrivial reduced product of
several copies of $z_1^{\pm1}$ and $z_2^{\pm1}$.
Let $m$ be an arbitrary integer such that $m\ge2n$.
Since, by Proposition~\ref{proposition:Gn_bs}, $G_n$ is simple or trivial,
the commutator length of $[w^m]$ in $G_n$ is defined.
To complete the proof, it is only left to show that
$\cl_{G_n}([w^m])>n$.

Suppose that on the contrary $\cl_{G_n}([w^m])\le n$.
By Lemma~\ref{lemma:diagrams}, there exists
a one-contour reduced diagram over \GP{\mathfrak A}{\mathcal R_n},
the underlying complex of whose closure
is a combinatorial handled sphere with $n$ or fewer handles,
and whose contour label is $w^m$.
Denote such a diagram by~$\Delta_0$.
Then $\chi_{\Delta_0}\ge1-2n$.

After cyclically shifting, if necessary, the c-contours of some of the
faces of $\Delta_0$,
endow $\Delta_0$ with the structure of a restricted correct $S$-diagram
(of sort $(\mathrm{I}.n)$)
without external exceptional arcs.
Transform $\Delta_0$ into a special $S$-diagram $\Delta_1$
by diamond moves.
This is possible by Lemma~\ref{lemma:special_diagrams} and
because $\Delta_0$ is reduced.
Let $\Delta$ be the connected component of $\Delta_1$ containing
$\cntr{}\Delta_1$.
Then a closure of $\Delta$ is a handled sphere.
Since diamond moves do not decrease the Euler characteristic,
the maximal possible Euler characteristic of a connected component
is $2$, and
every diamond move that increases the number of connected components
increases it by $1$ and increases the Euler characteristic by $2$,
it follows that $\chi_{\Delta}\ge\chi_{\Delta_0}$.
Note also that $\pathlbl(\cntr{}\Delta)=w^m$.

Case 1: $\Delta$ has no faces.
Let $F_{\mathfrak A}$ be the free group presented by
\GP{\mathfrak A}{\varnothing}.
Then, by Lemma~\ref{lemma:diagrams},
$[w^m]_{F_{\mathfrak A}}\in[F_{\mathfrak A},F_{\mathfrak A}]$
(hence $[w]_{F_{\mathfrak A}}\in[F_{\mathfrak A},F_{\mathfrak A}]$)
and $\cl_{F_{\mathfrak A}}[w^m]_{F_{\mathfrak A}}\le n$.
This contradicts with Corollary~5.2 in \cite{DuncanHowie:1991:gporplig}.
(That corollary implies, in particular, that
for every nontrivial element $x$ of the derived subgroup
of an arbitrary free group $F$, and for every $m\in\mathbb N$,
$\cl_F(x^m)\ge(m+1)/2$.)

Case 2: $\Delta$ has at least one face.
Let $\bar\Delta$ be a closure of $\Delta$.
Clearly, $\bar\Delta$ cannot be elementary spherical
(otherwise some cyclic shift of $w^{\pm m}$ would be a relator,
which is clearly not possible under the conditions imposed in
Subsection~\ref{subsection:bsglcw.group_presentations}).
The underlying complex of $\bar\Delta$ is a handled sphere
with at most $n$ handles since $\chi_{\bar\Delta}\ge2-2n$.
Let $\Theta$ be the face of $\bar\Delta$ that is not in $\Delta$
(the ``improper'' face).
Extend the existing structure of a (restricted special) $S$-diagram
on $\Delta$ to a structure of an (unrestricted special) $S$-diagram
on $\bar\Delta$, assigning to $\Theta$ rank $0$ and choosing all c-paths
of $\Theta$ as selected.
Then $\kappa_{\bar\Delta}(\Theta)=\kappa_{\bar\Delta}'(\Theta)=0$.

By Lemma~\ref{lemma:condition_D}, $\bar\Delta$ satisfies
$\mathsf D(\lambda_n,\mu_n,\nu_n)$ relative to $\Delta$.
By Lemma~\ref{lemma:condition_Y__singular_asphericity},
$\Delta$ satisfies $\mathsf Y$.
By induction and Inductive Lemma, using inequality
\thetag{\ref{display:main_inequality_for_bsglcw.ineq}},
obtain that $\bar\Delta$ satisfies $\mathsf Z(2)$ relative
to every simple disc subdiagram of $\Delta$.

Let $N$ be the sum of the lengths of all exceptional arcs
of $\bar\Delta$ (of $\Delta$),
$M$ be the sum of the lengths of all the non-exceptional maximal selected
arcs of $\bar\Delta$ that are incident to faces of $\Delta$
(recall that every non-exceptional maximal selected internal arc
of a correct $S$-diagram does not overlap with any exceptional arc), and
$L$ be the number of non-selected edges of $\bar\Delta$.
Note that every non-selected edge of $\bar\Delta$,
as well as every exceptional arc, is incident to a face of~$\Delta$.

Let $T$ be the set of all the edges of $\bar\Delta$ that
are incident to faces of $\Delta$.
Then $L+M+N=\norm{T}$.
To come to a contradiction, it is left to show that
$L+M+N\le(1/2)\sum_{\Pi\in\Delta(2)}\abs{\cntr\Pi}$,
because $\sum_{\Pi\in\Delta(2)}\abs{\cntr\Pi}<2\norm{T}$
(here it is used that $\Delta$ is non-degenerate).

The following upper estimate on $L$ follows from the condition
$\mathsf D_1(\lambda_\Delta)$:
$$
L\le\sum_{\Pi\in\Delta(2)}\!\lambda_n\abs{\cntr\Pi}.
$$

To estimate $M$, Estimating Lemma~\ref{estimating_lemma:selected_arcs}
and the condition $\mathsf D_2(\mu_n)$ shall be applied.
Let $A$ be the set of all the non-exceptional maximal
selected arcs that are incident to faces of $\Delta$.
Clearly, distinct elements of $A$ do not overlap and are not subarcs
of the same selected arc.
Let $B=\Delta(2)$, $C=\bar\Delta(2)$, $D=\{\Theta\}$.
By Estimating Lemma~\ref{estimating_lemma:selected_arcs} applied to
$\bar\Delta$, $A$, $B$, $C$, $D$,
there exist a subset $E$ of $A$ and a function
$h\!:A\setminus E\to\Delta(2)$ such that:
\begin{enumerate}
\item
    either $E$ is empty, or
    $$
    \norm{E}\le1+\kappa_{\bar\Delta}(\Theta)+\kappa_{\bar\Delta}'(\Theta)
    +2-3\chi_{\bar\Delta}\le3-3(2-2n)=6n-3;
    $$
\item
    for every $x\in A\setminus E$, the face $h(x)$ is incident to $x$;
\item
    for every face $\Pi$,
    the number of arcs mapped to $\Pi$ by $h$ is at most
    $3+\kappa_\Delta(\Pi)+\kappa_\Delta'(\Pi)\le3+2(4n+4)=8n+11$.
\end{enumerate}
Let $f\!:A\to\Delta(2)$ be an arbitrary extension of $h$ such that
for every $x\in A$, the face $f(x)$ is incident to $x$.
Then for every face $\Pi$ of $\Delta$,
the number of arcs mapped to $\Pi$ by $f$ is at most $14n+8$.
As follows from $\mathsf D_2(\mu_n)$,
$\abs{x}\le\mu_n\abs{\cntr{}(f(x))}$ for every $x\in A$.
Therefore,
$$
M\le\sum_{\Pi\in\Delta(2)}\!(14n+8)\mu_n\abs{\cntr\Pi}.
$$

It easily follows from Estimating
Lemma~\ref{estimating_lemma:exceptional_arcs}
applied to $\Delta$, that for every $j$,
the number of exceptional arcs of $\Delta$ of rank $j$
is at most $2n+1$ times the number of faces of $\Delta$ of rank $j$.
Then it follows from $\mathsf D_3(\nu_n)$ that
$$
N\le\sum_{\Pi\in\Delta(2)}\!(2n+1)\nu_n\abs{\cntr\Pi}.
$$

Thus, by inequality \thetag{\ref{display:main_inequality}}
and because $\Delta$ is non-degenerate,
$$
L+M+N<\sum_{\Pi\in\Delta(2)}\frac{1}{2}\abs{\cntr\Pi}.
$$
This leads to a contradiction in Case~2.
\end{proof}


\begin{proposition}
\label{proposition:G_icw}
The group\/ $G_\infty$ constructed in
Subsection\/~\textup{\ref{subsection:sgicw.group_presentations}}
is singularly aspherical\textup, torsion-free\textup, and
the elements\/ $[z_1]_{G_\infty}$ and $[z_2]_{G_\infty}$
freely generate a free subgroup\/ $H$ such that
$$
\bigl(\forall h\in H\setminus\{1\}\bigr)\,
\bigl(\lim_{n\to+\infty}\cl_{G_\infty}(h^n)=+\infty\bigr).
$$
\end{proposition}



\begin{proof}
The same way as in Proposition~\ref{proposition:Gn_lcw},
obtain that \GP{\mathfrak A}{\mathcal R_\infty} is singularly aspherical
and $G_\infty$ is torsion-free.

Let $w$ be an arbitrary nontrivial reduced product of
several copies of $z_1^{\pm1}$ and $z_2^{\pm1}$, and
$n$ be an arbitrary positive integer.
To complete this proof it is only left to show that
for every large enough $m$, $\cl_{G_\infty}([w^m])>n$.
Without loss of generality, assume that $w$ is cyclically reduced,
and that $\abs{w}\le\mu_j\abs{r_{j,1}}=\mu_j\abs{r_{j,2}}$
for every $j\ge n$.

Let $m$ be an arbitrary integer such that
$\abs{w^m}\ge\abs{r_{n,1}}$ and $1/m\le\mu_n$.
Suppose that $\cl_{G_\infty}([w^m])\le n$.
By the same argument as in the proof of
Proposition~\ref{proposition:Gn_lcw}, there exists
a one-contour reduced restricted special $S$-diagram over
\GP{\mathfrak A}{\mathcal R_\infty} (of sort $(\mathrm{II})$),
the underlying complex of whose closure
is a combinatorial handled sphere with at most $n$ handles,
and whose contour label is $w^m$.
Let $\Delta$ be such an $S$-diagram.

Let $\bar\Delta$ be a closure of $\Delta$.
The diagram $\bar\Delta$ cannot be elementary spherical
because of the conditions imposed in 
Subsection~\ref{subsection:sgicw.group_presentations}.
Let $\Theta$ be the face of $\bar\Delta$ that is not in $\Delta$.
Extend the existing structure of a (restricted special) $S$-diagram
on $\Delta$ to a structure of an (unrestricted special) $S$-diagram
on $\bar\Delta$, assigning to $\Theta$ rank $0$ and choosing all c-paths
of $\Theta$ as selected.
Then $\kappa_{\bar\Delta}(\Theta)=\kappa_{\bar\Delta}'(\Theta)=0$.

By Lemma~\ref{lemma:condition_D}, $\bar\Delta$ satisfies
$\mathsf D'(\lambda_\Delta,\mu_\Delta,\nu_\Delta)$
relative to $\Delta$.
By Lemma~\ref{lemma:condition_Y__singular_asphericity},
$\Delta$ satisfies $\mathsf Y$.
By induction and Inductive Lemma, using inequality
\thetag{\ref{display:main_inequality_for_sgicw.ineq}},
obtain that $\bar\Delta$ satisfies $\mathsf Z(2)$ relative
to every simple disc subdiagram of $\Delta$.

Let $N$ be the sum of the lengths of all exceptional arcs
of $\bar\Delta$,
$M$ be the sum of the lengths of all non-exceptional maximal selected
arcs of $\bar\Delta$, and
$L$ be the number of non-selected edges of $\bar\Delta$.
Then $L+M+N=\norm{\bar\Delta(1)}$.
To obtain a contradiction, it suffices to show that
$L+M+N<(1/2)\sum_{\Pi\in\bar\Delta(2)}\abs{\cntr\Pi}$.

The following upper estimate on $L$ follows from the condition
$\mathsf D_1(\lambda_\Delta)$ relative to $\Delta$, because every c-edge
of $\Theta$ is selected:
$$
L\le\sum_{\Pi\in\Delta(2)}\!\lambda_\Delta(\Pi)\abs{\cntr\Pi}.
$$

Let $A$ be the set of all maximal selected arcs of $\bar\Delta$.
Then $M+N=\sum_{x\in A}\abs{x}$.

Apply Estimating Lemma~\ref{estimating_lemma:selected_arcs}
to the $S$-diagram $\bar\Delta$ and the sets $A$,
$B=\Delta(2)$, $C=\bar\Delta(2)$, $D=\{\Theta\}$.
Let $E$ be a subset of $A$ and $h$ be a function
$A\setminus E\to\Delta(2)$ such that:
\begin{enumerate}
\item
    either $E$ is empty, or
    $$
    \norm{E}\le3-3\chi_{\bar\Delta}\le6n-3;
    $$
\item
    for every $x\in A\setminus E$, the face $h(x)$ is incident to $x$;
\item
    for every face $\Pi$,
    the number of arcs mapped to $\Pi$ by $h$ is at most
    $3+\kappa_\Delta(\Pi)+\kappa_\Delta'(\Pi)\le8\rank(\Pi)+11$.
\end{enumerate}
Let $M_1$ be the sum of the lengths of all non-exceptional
elements of $A\setminus E$,
and $M_2$ be the sum of the lengths of all non-exceptional
elements of $E$.
Then $M_1+M_2=M$, and, by the condition $\mathsf D_2(\mu_\Delta)$,
$$
M_1
\le\sum_{\Pi\in\Delta(2)}\!(8\rank(\Pi)+11)\mu_\Delta(\Pi)\abs{\cntr\Pi},
$$
while $M_2$ is less than or equal to $6n-3$ times
the maximal length of a non-exceptional maximal selected arc.

Apply Estimating Lemma~\ref{estimating_lemma:exceptional_arcs}
to $\Delta$.
Let $F$ be a set of exceptional arcs of $\Delta$ such that:
\begin{enumerate}
\item
    either $F$ is empty, or
    $\norm{F}\le-\chi_\Delta\le2n-1$, and
\item
    for every $j$,
    the number of exceptional arcs of rank $j$ that are not in $F$
    is at most twice the number of faces of $\Delta$ of rank $j$.
\end{enumerate}
Let $N_1$ be the sum of the lengths of all the exceptional arcs that
are not elements of $F$,
and $N_2$ be the sum of the lengths of all the elements of $F$.
Then $N_1+N_2=N$, and, by the condition $\mathsf D_3(\nu_\Delta)$,
$$
N_1\le\sum_{\Pi\in\Delta(2)}\!2\nu_\Delta(\Pi)\abs{\cntr\Pi},
$$
while $N_2$ is less than or equal to $2n-1$ times
the maximal length of an exceptional arc.

It is left to find suitable ``global'' estimates on the lengths of
non-exceptional maximal selected arc, and
on the lengths of exceptional ones.

Observe that the length of every arc of $\bar\Delta$ that is
incident to $\Theta$ and not incident to any other faces
is at most $\abs{w}-1<(1/m)\abs{\cntr\Theta}\le\mu_n\abs{\cntr\Theta}$.
Indeed, the label of each of the oriented arcs associated with such an arc
is a common subword of a power of $w$ and of a power of $w^{-1}$
(because $\bar\Delta$ is orientable).
Any such word of length $\abs{w}$ would be a cyclic
shift of $w$ and of $w^{-1}$ in the same time, but 
in a free group a nontrivial element is not conjugate to its own inverse
(if it was, it would commute with the square of the conjugating element,
and hence would commute with the conjugating element itself).

Case 1: the rank of every face of $\Delta$ is less than $n$.
If follows from conditions of 
Subsection~\ref{subsection:sgicw.group_presentations} and from the
inequality $\abs{\cntr\Theta}\ge\abs{r_{n,1}}$, that
the length of any non-exceptional maximal selected arc of
$\bar\Delta$ cannot be greater than $\mu_n\abs{\cntr\Theta}$,
the length of any exceptional arc of
$\bar\Delta$ cannot be greater than $\nu_n\abs{\cntr\Theta}$.
Therefore,
$$
M_2\le(6n-3)\mu_n\abs{\cntr\Theta}\quad\text{and}\quad
N_2\le(2n-1)\nu_n\abs{\cntr\Theta}.
$$
Thus, by inequalities \thetag{\ref{display:main_inequality}} and
\thetag{\ref{display:main_inequality_for_sgicw.ineq}}, obtain
a contradiction:
\begin{align*}
L+M+N
&\le\sum_{\Pi\in\Delta(2)}\!\lambda_\Delta(\Pi)\abs{\cntr\Pi}\\
&\qquad+\sum_{\Pi\in\Delta(2)}\!
(8\rank(\Pi)+11)\mu_\Delta(\Pi)\abs{\cntr\Pi}
+(6n-3)\mu_n\abs{\cntr\Theta}\\
&\qquad+\sum_{\Pi\in\Delta(2)}\!2\nu_\Delta(\Pi)\abs{\cntr\Pi}
+(2n-1)\nu_n\abs{\cntr\Theta}\\
&<\sum_{\Pi\in\bar\Delta(2)}\frac{1}{2}\abs{\cntr\Pi}.
\end{align*}

Case 2: $\Delta$ has a face of rank at least $n$.
Let $\hat\Pi$ be a face of $\Delta$ of maximal rank, 
$\rank(\hat\Pi)\ge n$.
By the condition $\mathsf D_2'(\mu_\Delta)$ and by inequality
$\abs{w}\le\mu_\Delta(\hat\Pi)\abs{\cntr\hat\Pi}$, the length of every
non-exceptional maximal selected arc of $\bar\Delta$ is at most
$\mu_\Delta(\hat\Pi)\abs{\cntr\hat\Pi}$.
By the condition $\mathsf D_3'(\nu_\Delta)$, the length of every
exceptional arc of $\bar\Delta$ is at most
$\nu_\Delta(\hat\Pi)\abs{\cntr\hat\Pi}$.
By inequality
\thetag{\ref{display:main_inequality_for_sgicw.ineq}}, obtain
a contradiction:
\begin{align*}
L+M+N
&\le\sum_{\Pi\in\Delta(2)}\!\lambda_\Delta(\Pi)\abs{\cntr\Pi}\\
&\qquad+\sum_{\Pi\in\Delta(2)}\!
(8\rank(\Pi)+11)\mu_\Delta(\Pi)\abs{\cntr\Pi}
+(6n-3)\mu_\Delta(\hat\Pi)\abs{\cntr\hat\Pi}\\
&\qquad+\sum_{\Pi\in\Delta(2)}\!2\nu_\Delta(\Pi)\abs{\cntr\Pi}
+(2n-1)\nu_\Delta(\hat\Pi)\abs{\cntr\hat\Pi}\\
&\le\sum_{\Pi\in\Delta(2)}\bigl(\lambda_\Delta(\Pi)
+(14\rank(\Pi)+8)\mu_\Delta(\Pi)\\
&\qquad+(2\rank(\Pi)+1)\nu_\Delta(\Pi)\bigr)
\abs{\cntr\Pi}\\
&<\sum_{\Pi\in\Delta(2)}\frac{1}{2}\abs{\cntr\Pi}.
\end{align*}
\end{proof}


It remains to show that the word and conjugacy problems
in the constructed groups are decidable.
Proving this fact could be facilitated by imposing additional restrictions
on the constructed presentations, but this is not necessary.

Observe that inequality \thetag{\ref{display:main_inequality}}
implies that for every $n\in\mathbb N$,
$$
\lambda_n+(8n+11)\mu_n+2\nu_n<\frac{19}{44}<0.45.
$$

The following lemma is helpful for solving the word and conjugacy
problems in the constructed groups.


\begin{lemma}
\label{lemma:minimal_disc_annular_maps}
Let\/ \GP{\mathfrak A}{\mathcal S} be a subpresentation of
one of the presentations\/
\GP{\mathfrak A}{\mathcal R_n}\textup, $n\in\mathbb N$\textup, or\/
\GP{\mathfrak A}{\mathcal R_\infty}\textup.
Let\/ $K$ be the group presented by\/ \GP{\mathfrak A}{\mathcal S}\textup.
Then
\begin{enumerate}
\item
    if\/ $w$ is a nontrivial group word over\/ $\mathfrak A$\textup,
    and\/ $\Delta$ is a minimal by the number of faces
    disc diagram over\/ \GP{\mathfrak A}{\mathcal S} such that\/
    $\pathlbl(\cntr{}\Delta)=w$\textup, then
    $$
    \abs{w}>\frac{1}{10}\sum_{\Pi\in\Delta(2)}\!\abs{\cntr\Pi}
    $$
    \textup(in particular\textup, $K$ is hyperbolic
    if\/ $\mathcal S$ is finite\/\textup{);}
\item
    if\/ $w_1$ and\/ $w_2$ are group words over\/ $\mathfrak A$
    such that\/ $[w_1]_K\ne1_K$\textup,
    and\/ $\Delta$ is a minimal by the number of faces
    contour-oriented annular diagram over\/ \GP{\mathfrak A}{\mathcal S}
    such that\/
    $\pathlbl(\cntr_1\Delta)=w_1$ and\/
    $\pathlbl(\cntr_2\Delta)^{-1}=w_2$\textup, then
    $$
    \abs{w_1}+\abs{w_2}>\frac{1}{10}
    \sum_{\Pi\in\Delta(2)}\!\abs{\cntr\Pi}.
    $$
\end{enumerate}
\end{lemma}


\begin{proof}
First, let $w$ be a group word over $\mathfrak A$,
and $\Delta$ be a minimal by the number of faces
disc diagram over \GP{\mathfrak A}{\mathcal S} such that
$\pathlbl(\cntr{}\Delta)=w$.
Because of minimality, $\Delta$ is reduced.
After cyclically shifting, if necessary, the c-contours of some of the
faces of $\Delta$,
endow $\Delta$ with a structure of a restricted correct $S$-diagram
of appropriate sort.
Transform $\Delta$ into a special $S$-diagram $\Gamma$ by diamond moves.
The connected component of $\Gamma$ containing $\cntr{}\Gamma$
is a disc diagram.
By the minimality of $\Delta$, this implies that $\Gamma$ is connected.
The $S$-diagram $\Gamma$ satisfies the conditions $\mathsf Y$
and $\mathsf D(\lambda_\Gamma,\mu_\Gamma,\nu_\Gamma)$.
Let
$$
\gamma=\lambda_\Gamma+(3+\kappa_\Gamma+\kappa_\Gamma')\mu_\Gamma
+2\nu_\Gamma.
$$
By inequalities
\thetag{\ref{display:main_inequality_for_bsglcw.ineq}} and
\thetag{\ref{display:main_inequality_for_sgicw.ineq}},
$\gamma(\Pi)<0.45$ for every $\Pi\in\Gamma$.
By Lemma~\ref{lemma:main_1},
$$
\abs{w}=\abs{\cntr{}\Gamma}
\ge\sum_{\Pi\in\Gamma(2)}\!(1-2\gamma(\Pi))\abs{\cntr\Pi}
\ge\frac{1}{10}\sum_{\Pi\in\Gamma(2)}\!\abs{\cntr\Pi}
=\frac{1}{10}\sum_{\Pi\in\Delta(2)}\!\abs{\cntr\Pi},
$$
and the equality in the both inequalities simultaneously
is not possible because $\abs{w}>0$.

Second, let $w_1$ and $w_2$ be group words over $\mathfrak A$
such that $[w_1]_G\ne1_G$,
and $\Delta$ be a minimal by the number of faces
contour-oriented restricted correct annular $S$-diagram
of appropriate sort over \GP{\mathfrak A}{\mathcal S}
such that
$\pathlbl(\cntr_1\Delta)=w_1$ and $\pathlbl(\cntr_2\Delta)^{-1}=w_2$.
Because of minimality, $\Delta$ is reduced.
Transform $\Delta$ into a special $S$-diagram $\Gamma$ by diamond moves.
The connected component of $\Gamma$ that contains $\cntr_1\Gamma$
is either disc or annular, and is contour-oriented.
By the minimality of $\Delta$, and because $[w_1]_G\ne1_G$,
this implies that $\Gamma$ is connected.
The $S$-diagram $\Gamma$ satisfies the conditions $\mathsf Y$
and $\mathsf D(\lambda_\Gamma,\mu_\Gamma,\nu_\Gamma)$.
Let $\gamma$ be as above.
Then $\gamma(\Pi)<0.45$ for every $\Pi\in\Gamma$.
By Lemma~\ref{lemma:main_1},
$$
\abs{w_1}+\abs{w_2}=\abs{\cntr_1\Gamma}+\abs{\cntr_2\Gamma}
\ge\sum_{\Pi\in\Gamma(2)}\!(1-2\gamma(\Pi))\abs{\cntr\Pi}
\ge\frac{1}{10}\sum_{\Pi\in\Delta(2)}\!\abs{\cntr\Pi},
$$
and the equality in the both inequalities is not possible simultaneously.
\end{proof}


\begin{proposition}
\label{proposition:Gn_G_swcp}
The groups\/ $G_n$\textup, $n\in\mathbb N$\textup, and\/ $G_\infty$
constructed in Section\/~\textup{\ref{section:group_presentations}}
have decidable word and conjugacy problems\textup.
\end{proposition}


\begin{proof}
Here follows a proof of decidability of the word an conjugacy problems
for the group $G_\infty$.
For the groups $G_n$, $n\in\mathbb N$, there is a completely
analogous proof which therefore shall not be given here.

For every $i\in\mathbb N$, let $w_i$, $r_{i,1}$ and $r_{i,2}$
be the same $w_i$, $r_{i,1}$ and $r_{i,2}$ as
in Subsection~\ref{subsection:sgicw.group_presentations}.

It is clear that the sequence
$r_{1,1}$, $r_{1,2}$, $r_{2,1}$, $r_{2,2}$, $r_{3,1}$, \dots
is recursive, and the sequence
$\abs{r_{1,1}}$, $\abs{r_{1,2}}$, $\abs{r_{2,1}}$, $\abs{r_{2,2}}$,
$\abs{r_{3,1}}$, \dots
is bounded from below by an increasing recursive sequence
(of rational numbers) tending to~$+\infty$.

Consider an arbitrary subset $\mathcal S$ of $\mathcal R_\infty$.
Let $K$ be the group presented by the (sub)presentation
\GP{\mathfrak A}{\mathcal S}.
It follows from Lemmas~\ref{lemma:diagrams} and
\ref{lemma:minimal_disc_annular_maps} that:
\begin{enumerate}
\item
    for every group word $x$ over $\mathfrak A$,
    $[x]_K=1_K$ (if and) only if there exists a disc diagram $\Delta$
    over \GP{\mathfrak A}{\mathcal S} such that
    $\pathlbl(\cntr{}\Delta)=x$ and
    $\norm{\Delta(1)}\le6\abs{x}$
    (because $(10+1)/2<6$);
\item
    for every group words $x$ and $y$ over $\mathfrak A$ such that
    $[x]_K\ne1_K$,
    $[x]_K$ and $[y]_K$ are conjugate in $K$ (if and) only if
    there exists a contour-oriented annular diagram $\Delta$
    over \GP{\mathfrak A}{\mathcal S} such that
    $\pathlbl(\cntr_1\Delta)=x$,
    $\pathlbl(\cntr_2\Delta)^{-1}=y$, and
    $\norm{\Delta(1)}<6(\abs{x}+\abs{y})$.
\end{enumerate}

The following algorithm decides the conjugacy problem for
\GP{\mathfrak A}{\mathcal R_\infty}:
Let $x$ and $y$ be arbitrary group words over $\mathfrak A$
given as an input.
Using the (effective) lower bound on $\abs{r_{i,j}}$, find a $k$ such that
for every $r\in\mathcal R_\infty$, if $(1/10)\abs{r}\le\abs{x}+\abs{y}$,
then $r\in\mathcal R_\infty^{(k)}$.
(Then, by Lemmas~\ref{lemma:diagrams} and
\ref{lemma:minimal_disc_annular_maps}, $x$ and $y$ represent conjugate
elements of $G_\infty$ if and only if
they represent conjugate elements of
\GP{\mathfrak A}{\mathcal R_\infty^{(k)}}.)
Determine the (finite) set $\mathcal R_\infty^{(k)}$ by finding
all the sets $\mathcal R_\infty^{(1)}$, $\mathcal R_\infty^{(2)}$, \dots,
$\mathcal R_\infty^{(k)}$ one-by-one in this order.
Do so in $k$ steps.
On the step number $i$, the set $\mathcal R_\infty^{(i-1)}$
is already determined.
To determine $\mathcal R_\infty^{(i)}$,
decide first whether $[w_i]_{\mathcal R_\infty^{(i-1)}}=1$
by checking if there exists a disc diagrams $\Delta$ over
\GP{\mathfrak A}{\mathcal R_\infty^{(i-1)}} with at most
$6\abs{w_i}$ edges and with the contour label $w_i$.
If $[w_i]_{\mathcal R_\infty^{(i-1)}}=1$, then
$\mathcal R_\infty^{(i)}=\mathcal R_\infty^{(i-1)}$, otherwise
$\mathcal R_\infty^{(i)}=\mathcal R_\infty^{(i-1)}\cup\{r_{i,1},r_{i,2}\}$.
After the set $\mathcal R_\infty^{(k)}$ is found,
decide whether $[x]_{G_\infty}=1_{G_\infty}$.
Do so by checking if there exists a disc diagram $\Delta$ over
\GP{\mathfrak A}{\mathcal R_\infty^{(k)}} with at most $6\abs{x}$
edges and with the contour label $x$.
If found that $[x]_{G_\infty}=1_{G_\infty}$,
then similarly decide whether $[y]_{G_\infty}=1_{G_\infty}$.
In the case $[x]_{G_\infty}=1_{G_\infty}=[y]_{G_\infty}$,
the elements $[x]_{G_\infty}$ and $[y]_{G_\infty}$
are conjugate in $G_\infty$,
and in the case $[x]_{G_\infty}=1_{G_\infty}\ne[y]_{G_\infty}$,
they are not.
If found that $[x]_{G_\infty}\ne1_{G_\infty}$, decide whether
$[x]_{G_\infty}$ and $[y]_{G_\infty}$ are conjugate
in $G_\infty$ by checking whether
there exists a contour-oriented annular diagram $\Delta$ over
\GP{\mathfrak A}{\mathcal R_\infty^{(k)}} with
less than $6(\abs{x}+\abs{y})$ edges and with the contour labels
$x$ and $y^{-1}$.

Thus the group $G_\infty$ has decidable word and conjugacy problems.
\end{proof}


Theorems \ref{theorem:bsglcw} and \ref{theorem:sgicw}
are direct corollaries of Propositions
\ref{proposition:Gn_bs}, \ref{proposition:G_s}, \ref{proposition:Gn_lcw},
\ref{proposition:G_icw}, \ref{proposition:Gn_G_swcp}.



\section*{Acknowledgements}

The author thanks Valerij Bardakov and Daniela Nikolova for
bringing questions about commutator width of simple groups
to author's attention.
The author is grateful to Alexander Ol'shanskii for helpful discussions.
Some of the recent developments in the area were 
pointed out to the author by Yves de Cornulier.


\nocite{Olshanskii:1989:gosg-rus}

\bibliographystyle{amsalpha}
\bibliography{\string~/Documents/bib}

\newcommand{\noopsort}[1]{}\newcommand{\singleletter}[1]{#1}
\providecommand{\bysame}{\leavevmode\hbox to3em{\hrulefill}\thinspace}
\providecommand{\MR}{\relax\ifhmode\unskip\space\fi MR }
\providecommand{\MRhref}[2]{%
  \href{http://www.ams.org/mathscinet-getitem?mr=#1}{#2}
}
\providecommand{\href}[2]{#2}
\begin{thebibliography}{CCH81}

\bibitem[BG92]{BargeGhys:1992:cEM-fr}
Jean Barge and {\'E}tienne Ghys, \emph{Cocycles d'{E}uler et de {M}aslov
  [{E}uler and {M}aslov cocycles]}, Math.\ Ann. \textbf{294} (1992), no.~2,
  235--265, in French.

\bibitem[Bro94]{Brown:1994:cg}
Kenneth~S. Brown, \emph{Cohomology of groups}, Springer-Verlag, 1994, corrected
  reprint of the 1982 original.

\bibitem[BRS76]{BuoncristianoRourkeSanderson:1976:gaht}
Sandro Buoncristiano, Colin~P. Rourke, and Brian~J. Sanderson, \emph{A
  geometric approach to homology theory}, London Mathematical Society Lecture
  Note Series, no.~18, Cambridge University Press, 1976.

\bibitem[CCH81]{ChiColHue:1981:agp}
Ian~M. Chiswell, Donald~J. Collins, and Johannes Huebschmann, \emph{Aspherical
  group presentations}, Math Z. \textbf{178} (1981), no.~1, 1--36.

\bibitem[DH91]{DuncanHowie:1991:gporplig}
Andrew~J. Duncan and James Howie, \emph{The genus problem for one-relator
  products of locally indicable groups}, Math.\ Z. \textbf{208} (1991), no.~2,
  225--237.

\bibitem[GG04]{GambaudoGhys:2004:cds}
Jean-Marc Gambaudo and {\'E}tienne Ghys, \emph{Commutators and diffeomorphisms
  of surfaces}, Ergod.\ Th.\ \& Dynam.\ Sys. \textbf{24} (2004), 1591--1617.

\bibitem[Hal35]{Hall:1935:ors}
Philip Hall, \emph{On representatives of subsets}, J.\ London Math.\ Soc.
  \textbf{10} (1935), 26--30.

\bibitem[Hue79]{Huebschmann:1979:ctagscg}
Johannes Huebschmann, \emph{Cohomology theory of aspherical groups and of small
  cancellation groups}, J.\ Pure Appl.\ Algebra \textbf{14} (1979), no.~2,
  137--143.

\bibitem[Hue80]{Huebschmann:1980:htcap}
\bysame, \emph{The homotopy type of a combinatorially aspherical presentation},
  Math.\ Z. \textbf{173} (1980), no.~2, 163--169.

\bibitem[Hue81]{Huebschmann:1981:a2cupW}
\bysame, \emph{Aspherical\/ $2$-complexes and an unsettled problem of {J.} {H.}
  {C.} {W}hitehead}, Math.\ Ann. \textbf{258} (1981), 17--37.

\bibitem[Isa77]{Isaacs:1977:ccs}
I.~Martin Isaacs, \emph{Commutators and the commutator subgroup}, Amer.\ Math.\
  Monthly \textbf{84} (1977), no.~9, 720--722.

\bibitem[LS01]{LyndonSchupp:2001:cgt}
Roger~C. Lyndon and Paul~E. Schupp, \emph{Combinatorial group theory},
  Springer-Verlag, 2001, reprint of the 1977 edition.

\bibitem[McC00]{McCammond:2000:gsct}
Jonathan~P. McCammond, \emph{A general small cancellation theory}, Internat.\
  J.\ Algebra Comput. \textbf{10} (2000), no.~1, 1--172.

\bibitem[MK99]{kn14:1999:knupgt-eng}
V.~D. Mazurov and E.~I. Khukhro (eds.), \emph{The {K}ourovka {N}otebook:
  unsolved problems in group theory}, 14th augmented ed., Russian Acad.\ of
  Sci.\ Siber.\ Div., Inst.\ Math., Novosibirsk, \noopsort{1999b}1999,
  translated from Russian.

\bibitem[Mur05]{Muranov:2005:dsmcbgbsg}
Alexey~Yu. Muranov, \emph{Diagrams with selection and method for constructing
  boundedly generated and boundedly simple groups}, Comm.\ Algebra \textbf{33}
  (2005), no.~4, 1217--1258, arXiv.org preprint:
  \href{http://arxiv.org/abs/math.GR/0404472}{math.GR/0404472}.

\bibitem[Mur07]{Muranov:2007:otfgfrfb}
\bysame, \emph{On torsion-free groups with finite regular file bases}, Trans.\
  Amer.\ Math.\ Soc. \textbf{359} (\noopsort{2007-01}2007), 3609--3645,
  arXiv.org preprint:
  \href{http://arxiv.org/abs/math.GR/0504438}{math.GR/0504438}.

\bibitem[Ol'89]{Olshanskii:1989:gosg-rus}
Alexander~Yu. Ol'shanskii, \emph{Geometrija opredeljajushchikh sootnoshenij v
  gruppakh [{G}eometry of defining relations in groups]}, Na\-u\-ka, Moscow,
  \noopsort{1989a}1989, in Russian.

\bibitem[Ol'91]{Olshanskii:1991:gdrg-eng}
\bysame, \emph{Geometry of defining relations in groups}, Kluwer Academic
  Publishers, Dordrecht, Boston, \noopsort{1989b}1991, translated from Russian.

\bibitem[Ore51]{Ore:1951:src}
Oystein Ore, \emph{Some remarks on commutators}, Proc.\ Amer.\ Math.\ Soc.
  \textbf{2} (1951), no.~2, 307--314.

\bibitem[Rou79]{Rourke:1979:ptg}
Colin~P. Rourke, \emph{Presentations and the trivial group}, Topology of
  low-dimensional manifolds ({P}roc. {S}econd {S}ussex {C}onf., {C}helwood
  {G}ate, 1977) (Berlin), Lecture Notes in Math., vol. 722, Springer, 1979,
  pp.~134--143.

\bibitem[Wil96]{Wilson:1996:fogt}
John~S. Wilson, \emph{First-order group theory}, Infinite groups '94
  (Berlin---New York) (Francesco de~Giovanni and Martin~L. Newell, eds.),
  Walter de Gruyter \& Co., 1996, proceedings of the international conference
  held in Ravello, Italy, May 23--27, 1994, pp.~301--314.

\end{thebibliography}


\end{document}